\documentclass[10pt, reqno]{amsart}
\usepackage{amssymb}
\usepackage{eucal}
\usepackage{amsmath}
\usepackage{amscd}
\usepackage[dvips]{color}
\usepackage{mathtools}
\usepackage{multicol}
\usepackage[all]{xy}           
\usepackage{graphicx}
\usepackage{color}
\usepackage{colordvi}
\usepackage{xspace}
\usepackage{ulem}
\usepackage{cancel}
\usepackage{tikz}
\usepackage{longtable}
\usepackage{multirow}
\usepackage{ifpdf}
\ifpdf
 \usepackage[colorlinks,final,backref=page,hyperindex]{hyperref}
\else
 \usepackage[colorlinks,final,backref=page,hyperindex,hypertex]{hyperref}
\fi

\begin{document}
\newcommand {\emptycomment}[1]{} 

\newcommand{\tabincell}[2]{\begin{tabular}{@{}#1@{}}#2\end{tabular}}

\newcommand{\nc}{\newcommand}
\newcommand{\delete}[1]{}

\nc{\mlabel}[1]{\label{#1}}  
\nc{\mcite}[1]{\cite{#1}}  
\nc{\mref}[1]{\ref{#1}}  
\nc{\meqref}[1]{~\eqref{#1}} 
\nc{\mbibitem}[1]{\bibitem{#1}} 

\delete{
\nc{\mlabel}[1]{\label{#1}  
{\hfill \hspace{1cm}{\bf{{\ }\hfill(#1)}}}}
\nc{\mcite}[1]{\cite{#1}{{\bf{{\ }(#1)}}}}  
\nc{\mref}[1]{\ref{#1}{{\bf{{\ }(#1)}}}}  
\nc{\meqref}[1]{~\eqref{#1}{{\bf{{\ }(#1)}}}} 
\nc{\mbibitem}[1]{\bibitem[\bf #1]{#1}} 
}

\newtheorem{thm}{Theorem}[section]
\newtheorem{lem}[thm]{Lemma}
\newtheorem{cor}[thm]{Corollary}
\newtheorem{pro}[thm]{Proposition}
\newtheorem{conj}[thm]{Conjecture}
\theoremstyle{definition}
\newtheorem{defi}[thm]{Definition}
\newtheorem{ex}[thm]{Example}
\newtheorem{rmk}[thm]{Remark}
\newtheorem{pdef}[thm]{Proposition-Definition}
\newtheorem{condition}[thm]{Condition}

\renewcommand{\labelenumi}{{\rm(\alph{enumi})}}
\renewcommand{\theenumi}{\alph{enumi}}
\renewcommand{\labelenumii}{{\rm(\roman{enumii})}}
\renewcommand{\theenumii}{\roman{enumii}}

\nc{\tred}[1]{\textcolor{red}{#1}}
\nc{\tblue}[1]{\textcolor{blue}{#1}}
\nc{\tgreen}[1]{\textcolor{green}{#1}}
\nc{\tpurple}[1]{\textcolor{purple}{#1}}
\nc{\btred}[1]{\textcolor{red}{\bf #1}}
\nc{\btblue}[1]{\textcolor{blue}{\bf #1}}
\nc{\btgreen}[1]{\textcolor{green}{\bf #1}}
\nc{\btpurple}[1]{\textcolor{purple}{\bf #1}}


\newcommand{\End}{\text{End}}

\nc{\calb}{\mathcal{B}}
\nc{\call}{\mathcal{L}}
\nc{\calo}{\mathcal{O}}
\nc{\frakg}{\mathfrak{g}}
\nc{\frakh}{\mathfrak{h}}
\nc{\ad}{\mathrm{ad}}

\nc{\ccred}[1]{\tred{\textcircled{#1}}}


\newcommand{\cm}[1]{\textcolor{purple}{\underline{CM:}#1 }}
\newcommand{\gl}[1]{\textcolor{blue}{\underline{GL:}#1 }}

\title[A new class of Manin triples of compatible Lie algebras]
{A new class of Manin triples of compatible Lie algebras}


\author{Guilai Liu}
\address{Chern Institute of Mathematics \& LPMC, Nankai University, Tianjin 300071, China
	\newline
	Institut de Recherche Math\'ematique Avanc\'ee,
	UMR 7501, Universit\'e de Strasbourg et CNRS,
	7 rue Ren\'e-Descartes, 67000 Strasbourg, France
}
\email{liugl@mail.nankai.edu.cn}

\begin{abstract}
	This paper introduces a new class of Manin triples of compatible Lie algebras and investigates the corresponding bialgebra structures.
	To this end, we first study an invariant-c$2$c form on a compatible Lie algebra, namely a nondegenerate symmetric bilinear form that is invariant with respect to one Lie bracket and is a commutative $2$-cocycle with respect to the other, which naturally arises from a quadratic differential Novikov algebra.
	Then we introduce the notion of compatible Lie-anti-pre-Lie algebras (compatible LAPL algebras in short), which serve as the underlying algebra structures of compatible Lie algebras endowed with invariant-c$2$c forms.
	Furthermore, we introduce the notion of a compatible LAPL bialgebra, which is shown to be equivalent to a Manin triple of compatible Lie algebras with respect to the invariant-c$2$c form.
	We investigate triangular compatible LAPL bialgebras, which are obtained from antisymmetric solutions of the compatible LAPL Yang-Baxter equation.
	We also use $\mathcal{O}$-operators of compatible LAPL algebras to construct such solutions.
	Finally, we show that every differential Novikov bialgebra gives rise to a compatible LAPL bialgebra and that this construction preserves triangularity.
\end{abstract}

\subjclass[2020]{
	    16T10, 
	17A36, 
	17B10, 
	17B38, 
	17B60, 
	17B62, 
	17D25. 
}

\keywords{Compatible Lie algebras, commutative $2$-cocycles, anti-pre-Lie algebras, compatible LAPL bialgebras, Yang-Baxter equations}

\maketitle


\tableofcontents

\allowdisplaybreaks

\section{Introduction}

Let $\mathcal{P}$ be an operad generated by a binary operation and let
$(\frak g,\circ_{\frak g}),(\frak g,\circ'_{\frak g})$ be $\mathcal{P}$-algebras.
If for all $\lambda,\mu\in \mathbb{K}$, the pair $(\frak g,\circ_{\lambda,\mu})$ with the multiplication $\circ_{\lambda,\mu}$ given by the linear combination
\begin{equation*}
	x\circ_{\lambda,\mu}y=\lambda x\circ_{\frak g}y+\mu x\circ'_{\frak g}y,\;\forall x,y\in\frak g
\end{equation*}
is again a $\mathcal{P}$-algebra, then the triple $(\frak g,\circ_{\frak g},\circ'_{\frak g})$ is called a {\bf compatible $\mathcal{P}$-algebra}.
This notion of compatible algebraic structures has attracted considerable attention in algebra, geometry and mathematical physics.
For example, compatible associative algebras have been extensively studied in connection with affine Dynkin diagrams, integrable matrix equations, infinitesimal bialgebras and quiver representations \cite{OS1,OS2,OS4}. Their operadic description, together with the corresponding free objects, was subsequently established in \cite{Dot}.
The notion of compatible Lie algebras was studied in detail due to their close relationships with compatible Poisson structures, bi-Hamiltonian systems and several deformation-theoretic constructions \cite{Gol1,Gol2,Gol3,LSB,OS3}.

A bialgebra is a vector space endowed with both an algebra structure and a coalgebra structure satisfying suitable compatibility conditions. Classical examples include Lie bialgebras \cite{CP1,Maj} which arise naturally in the study of Poisson-Lie groups and the infinitesimal theory of quantum groups, and antisymmetric infinitesimal (ASI) bialgebras \cite{Bai2010} which are closely related to symmetric Frobenius algebras.
These bialgebra structures admit a common characterization in terms of Manin triples with respect to suitable nondegenerate bilinear forms. More precisely, a Lie bialgebra is equivalent to a Manin triple of Lie algebras with respect to the nondegenerate symmetric invariant bilinear form, while an ASI bialgebra is equivalent to a double construction of a Frobenius algebra, or equivalently, a Manin triple of associative algebras with respect to the nondegenerate symmetric invariant bilinear form. This characterization has proved to be an effective framework for studying bialgebra structures associated with many other classes of algebras.

Different compatibility conditions on the nondegenerate bilinear forms in Manin triples can lead to different bialgebra structures.
For example, a Manin triple of Lie algebras with respect to the symplectic form (also called a parakähler Lie algebra) corresponds to a pre-Lie bialgebra \cite{Bai2008}, instead of a Lie bialgebra.
The notion of a compatible Lie bialgebra was introduced in \cite{WB}, which was shown to be equivalent to a Manin triple of compatible Lie algebras with respect to the nondegenerate symmetric bilinear form that is invariant with respect to both Lie brackets on the double space. It is quite natural to consider studying Manin triples of compatible Lie algebras with respect to other kinds of bilinear forms and the corresponding bialgebra structures.

Besides the symmetric invariant bilinear form, there is another type of symmetric bilinear form $\mathcal{B}$ on a Lie algebra
$(\frak g,[-,-]_{\frak g})$, namely the  commutative $2$-cocycle \cite{Dzh2010}, which satisfies the following condition:
\begin{equation*}
	\mathcal{B}([x,y]_{\frak g},z)+\mathcal{B}([y,z]_{\frak g},x)+\mathcal{B}([z,x]_{\frak g},y)=0,\;\forall x,y,z\in\frak g.
\end{equation*}
A nondegenerate commutative $2$-cocycle on a Lie algebra can be viewed as the symmetric analogue of a symplectic form on a Lie algebra.
The notion of an  anti-pre-Lie algebra was introduced in \cite{LB2022}, which serves as the underlying algebra structure of a Lie algebra with a nondegenerate commutative $2$-cocycle, in parallel with the correspondence between Lie algebras with symplectic forms and pre-Lie algebras \cite{Bai2021.2,Bur}.
Anti-pre-Lie algebras are shown to be closely related to many known types of nonassociative algebras, especially Novikov algebras.
The theory of anti-pre-Lie algebras has undergone rapid development in recent years.
In particular, the notion of an anti-pre-Lie bialgebra was introduced in \cite{TPA}, which was shown to be equivalent to a Manin triple of Lie algebras with respect to the nondegenerate commutative $2$-cocycle.
Anti-pre-Lie bialgebras serve as an important component of the bialgebra theories for transposed Poisson algebras \cite{TPA} and relative Poisson algebras \cite{LB2025}.

Up to now, compatible Lie algebras equipped with nondegenerate symmetric bilinear forms have been studied in two cases: the bilinear form that is invariant with respect to both Lie brackets \cite{WB}, and the bilinear form that is a commutative $2$-cocycle with respect to both Lie brackets. The latter naturally leads to the notion of a compatible anti-pre-Lie algebra \cite{Nor}.
Inspired by these two cases, in this paper we consider a mixed case in which a nondegenerate symmetric bilinear form on a compatible Lie algebra is invariant with respect to one Lie bracket and is a commutative $2$-cocycle with respect to the other. We call such a form an  invariant-c$2$c form.
We provide a natural construction of compatible Lie algebras with invariant-c$2$c forms from quadratic differential Novikov algebras.

To describe the underlying algebra structures of compatible Lie algebras with invariant-c$2$c forms, we introduce the notion of a compatible Lie-anti-pre-Lie algebra $(\frak g,[-,-]_{\frak g},\circ_{\frak g})$ (a compatible LAPL algebra in short), which combines a Lie algebra $(\frak g,[-,-]_{\frak g})$ and an anti-pre-Lie algebra $(\frak g,\circ_{\frak g})$ subject to suitable compatibility conditions.
A compatible Lie algebra with an invariant-c$2$c form gives rise to a compatible LAPL algebra.
Moreover, a compatible LAPL algebra $(\frak g,[-,-]_{\frak g},\circ_{\frak g})$ leads to a compatible Lie algebra on the double space $\frak g\oplus\frak g^{*}$ with an invariant-c$2$c form.
We also show that every differential Novikov algebra naturally gives rise to a compatible LAPL algebra.

Motivated by the above discussions, we study Manin triples of compatible Lie algebras with respect to the invariant-c$2$c form and the corresponding bialgebra structures. We introduce the notion of a compatible LAPL bialgebra,
which consists of a Lie bialgebra and an anti-pre-Lie bialgebra satisfying additional compatibility conditions.
Then we establish the equivalence between  a Manin triple of compatible Lie algebras with respect to the invariant-c$2$c form and a compatible LAPL bialgebra. 

We further investigate the construction theory of compatible LAPL bialgebras.
For this purpose, we introduce the compatible LAPL Yang--Baxter equation and show that every antisymmetric solution gives rise to a triangular compatible LAPL bialgebra.
We also introduce $\mathcal{O}$-operators of compatible LAPL algebras which provide a systematic approach to constructing such solutions.

Finally, we extend the construction from differential Novikov algebras to compatible LAPL algebras to the bialgebraic level.
We show that a Manin triple of differential Novikov algebras induces a Manin triple of compatible Lie algebras with respect to an invariant-c$2$c form.
Equivalently from the viewpoint of bialgebras, every differential Novikov bialgebra gives rise to a compatible LAPL bialgebra.
We further show that this construction preserves triangularity and establish the corresponding constructions at the levels of representations, $\mathcal{O}$-operators and pre-algebras.

Throughout this paper, unless otherwise specified, all the vector spaces and algebras are finite-dimensional over a  field $\mathbb{K}$ of characteristic $0$. Moreover, we fix the following notations.
\begin{enumerate}
	\item Let $\frak g$ be a vector space with a multiplication
	$\circ_{\frak g}:\frak g\otimes\frak g\rightarrow\frak g$.
	Define the linear maps
	$L_{\circ_{\frak g}},R_{\circ_{\frak g}}:\frak g\rightarrow\mathrm{End}_{\mathbb{K}}(\frak g)$
	by
	\begin{eqnarray*}
		L_{\circ_{\frak g}}(x)y=x\circ_{\frak g}y=R_{\circ_{\frak g}}(y)x,\;\forall x,y\in\frak g.
	\end{eqnarray*}
	In particular, for Lie algebras $(\frak g,[-,-]_{\frak g})$ and
	$(\frak g,\lbrace-,-\rbrace_{\frak g})$, we set the linear maps
	$\mathrm{ad}_{\frak g},\mathfrak{ad}_{\frak g}:\frak g\rightarrow\mathrm{End}_{\mathbb{K}}(\frak g)$
	by
	\begin{equation*}
		\mathrm{ad}_{\frak g}(x)y=[x,y]_{\frak g},\;
		\mathfrak{ad}_{\frak g}(x)y=\lbrace x,y\rbrace_{\frak g},\;\forall x,y\in\frak g.
	\end{equation*}
	
	\item For a vector space $\frak g$, define the linear maps
	$\tau:\frak g\otimes\frak g\rightarrow\frak g\otimes\frak g$ and
	$S:\frak g\otimes\frak g\otimes\frak g\rightarrow\frak g\otimes\frak g\otimes\frak g$
	by
	\begin{align*}
		\tau(x\otimes y)=y\otimes x,\;
		S(x\otimes y\otimes z)=x\otimes y\otimes z+y\otimes z\otimes x+z\otimes x\otimes y,\;\forall x,y,z\in\frak g.
	\end{align*}
	
	\item Let $\frak g$ be a vector space and $P:\frak g\rightarrow\frak g$
	be a linear map. Define a linear map
	$P^{t}:\frak g^{*}\rightarrow\frak g^{*}$ by
	\begin{eqnarray*}
		\langle P^{t}(a^{*}),x\rangle=\langle a^{*},P(x)\rangle,\;\forall x\in\frak g, a^{*}\in\frak g^{*},
	\end{eqnarray*}
	where $\langle-,-\rangle$ is the natural pairing between $\frak g$ and the dual space $\frak g^{*}$.
	
	\item Let $\frak g,V$ be vector spaces and
	$\mu:\frak g\rightarrow\mathrm{End}_{\mathbb{K}}(V)$ be a linear map.
	Define a linear map
	$\mu^{*}:\frak g\rightarrow\mathrm{End}_{\mathbb{K}}(V^{*})$ by
	\begin{eqnarray*}
		\langle\mu^{*}(x)u^{*},v\rangle=-\langle u^{*},\mu(x)v\rangle,\;\forall x\in\frak g, u^{*}\in V^{*}, v\in V.
	\end{eqnarray*}
	
	\item Let $\frak g$ be a vector space and $\mathcal{B}$ be a
	(nondegenerate) bilinear form on $\frak g$. Then there is an
	(invertible) linear map
	$\mathcal{B}^{\natural}:\frak g\rightarrow\frak g^{*}$ given by
	\begin{equation*}
		\langle\mathcal{B}^{\natural}(x),y\rangle=\mathcal{B}(x,y),\;\forall x,y\in\frak g.
	\end{equation*}
	
	\item Let $\frak g$ and $V$ be vector spaces and $r\in\frak g\otimes V$.
	Then $r$ is identified with a linear map from $\frak g^{*}$ to $V$
	(still denoted by $r$) by the following equation:
	\begin{equation*} 
		\langle r,a^{*}\otimes u^{*}\rangle=\langle r(a^{*}),u^{*}\rangle,\;\forall a^{*}\in\frak g^{*}, u^{*}\in V^{*}.
	\end{equation*}
 \item Let $V_{1},V_{2}$ be vector spaces and $P_{1}:V_{1}\rightarrow V_{1},P_{2}:V_{2}\rightarrow V_{2}$ be linear maps. Define a linear map $P_{1}+P_{2}:V_{1}\oplus V_{2}\rightarrow V_{1}\oplus V_{2}$ by
 \begin{equation*}
 	(P_{1}+P_{2})(v_{1}+v_{2})=P_{1}(v_{1})+P_{2}(v_{2}),\;\forall v_{1}\in V_{1},v_{2}\in V_{2}.
 \end{equation*}
\end{enumerate}

\section{Compatible Lie algebras with invariant-c$2$c forms and compatible LAPL algebras}\
In this section, we study compatible Lie algebras with invariant-c$2$c forms and compatible LAPL algebras.
First, we recall some basic foundations on anti-pre-Lie algebras, which serve as a necessary component of compatible LAPL algebras.
Then we present natural constructions of invariant-c$2$c forms on compatible Lie algebras.
Furthermore, we introduce the notion of compatible LAPL algebras, which naturally arise from compatible Lie algebras with invariant-c$2$c forms.

\subsection{Preliminaries on anti-pre-Lie algebras}\
We recall some basic foundations on anti-pre-Lie algebras. Most of the results can be found in \cite{LB2022} and \cite{TPA}.

\begin{defi}
	Let $\frak g$ be a vector space and $\circ_{\frak g}:\frak g\otimes\frak g\rightarrow\frak g$ be a multiplication.
	If the following equations hold:
	\begin{align*}
	& x\circ_{\frak g}(y\circ_{\frak g}z)-y\circ_{\frak g}(x\circ_{\frak g}z)=\lbrace y,x\rbrace_{\frak g}\circ_{\frak g}z,\\
	& \lbrace x,y\rbrace_{\frak g}\circ_{\frak g}z+\lbrace y,z\rbrace_{\frak g}\circ_{\frak g}x+\lbrace z,x\rbrace_{\frak g}\circ_{\frak g}y=0,\;\forall x,y,z\in\frak g,
	\end{align*}
where $\lbrace-,-\rbrace_{\frak g}:\frak g\otimes\frak g\rightarrow\frak g$ is given by
\begin{equation}\label{297}
\lbrace x,y\rbrace_{\frak g}=x\circ_{\frak g}y-y\circ_{\frak g}x,\;\forall x,y \in\frak g,
\end{equation}
then we say $(\frak g,\circ_{\frak g})$ is an {\bf anti-pre-Lie algebra}.
\end{defi}

Recall that a {\bf representation} of a Lie algebra $(\frak g,\lbrace-,-\rbrace_{\frak g})$ is a pair $(\mu,V)$ such that $V$ is a vector space, $\mu:\frak g\rightarrow\mathrm{End}_{\mathbb{K}}(V)$ is a linear map and the following equation holds:
\begin{equation*}
	\mu(\lbrace x,y\rbrace_{\frak g})v=\mu(x)\mu(y)v-\mu(y)\mu(x)v,\;\forall x,y\in \frak g, v\in V.
\end{equation*}
In  fact, the pair $(\mu,V)$ is a representation of $(\frak g,\lbrace-,-\rbrace_{\frak g})$ if and only if there is a Lie algebra structure $(\frak g\oplus V,\lbrace-,-\rbrace_{d})=\frak g\ltimes_{\mu}V$ on the direct sum $\frak g\oplus V$ of vector spaces
given by 
\begin{equation*}
	\lbrace x+u,y+v\rbrace_{d}=\lbrace x,y\rbrace_{\frak g}+\mu(x)v-\mu(y)u,\;\forall x,y\in\frak g,u,v\in V.
\end{equation*}
Moreover, if $(\mu,V)$ is a representation of $(\frak g,\lbrace-,-\rbrace_{\frak g})$, then $(\mu^{*},V^{*})$ is also a representation of $(\frak g,\lbrace-,-\rbrace_{\frak g})$.

\begin{pro}\label{317}
	Let $\frak g$ be a vector space, $\circ_{\frak g}:\frak g\otimes\frak g\rightarrow\frak g$ be a multiplication and $\lbrace-,-\rbrace_{\frak g}:\frak g\otimes\frak g\rightarrow\frak g$ be given by \eqref{297}.
	Then $(\frak g,\circ_{\frak g})$ is an anti-pre-Lie algebra if and only if the pair $(\frak g,\lbrace-,-\rbrace_{\frak g})$ is a Lie algebra and the negative multiplication operator $-L_{\circ_{\frak g}}:\frak g\rightarrow\mathrm{End}_{\mathbb{K}}(\frak g)$ induces a representation $(-L_{\circ_{\frak g}},\frak g)$ of $(\frak g,\lbrace-,-\rbrace_{\frak g})$.
	In this case, $(\frak g,\lbrace-,-\rbrace_{\frak g})$ is called the {\bf sub-adjacent Lie algebra} of $(\frak g,\circ_{\frak g})$ and $(-L^{*}_{\circ_{\frak g}},\frak g^{*})$ is also a representation of $(\frak g,\lbrace-,-\rbrace_{\frak g})$.
\end{pro}

Let $(\frak g,\lbrace-,-\rbrace_{\frak g})$ be a Lie algebra and  $(\mu,V)$ and $(\mu',V')$  be two representations of $(\frak g,\lbrace-,-\rbrace_{\frak g})$.
If there exists a linear isomorphism $\phi:V\rightarrow V'$  such that the following equation holds:
\begin{equation}\label{324}
\phi\big( \mu(x)v\big)=\mu'(x)\phi(v),\;\forall x\in\frak g, v\in V,
\end{equation}
then we say $(\mu,V)$ and $(\mu',V')$ are {\bf equivalent} through $\phi$.

\begin{pro}\label{322}
	Suppose that $\mathcal{B}$ is a nondegenerate commutative $2$-cocycle on a Lie algebra $(\frak g,\lbrace-,-\rbrace_{\frak g})$.
	Then there is an anti-pre-Lie algebra $(\frak g,\circ_{\frak g})$ with $\circ_{\frak g}:\frak g\otimes\frak g\rightarrow\frak g$ given by
	\begin{equation}\label{332-}
		\mathcal{B}(x\circ_{\frak g}y,z)=\mathcal{B}(y,\lbrace x,z\rbrace_{\frak g}),\;\forall x,y,z\in\frak g
	\end{equation}
whose sub-adjacent Lie algebra is 
 $(\frak g, \lbrace-,-\rbrace_{\frak g})$.
 Moreover, the adjoint representation $(\mathfrak{ad}_{\frak g}, \frak g)$
of $ (\frak g,\lbrace-,-\rbrace_{\frak g}) $
is equivalent to $(-L^{*}_{\circ_{\frak g}},\frak g^{*})$ through $\mathcal{B}^{\natural}$, that is, the following equation holds:
\begin{equation}\label{344}
\mathcal{B}^{\natural}\big( \mathfrak{ad}_{\frak g}(x)y  \big)=-L^{*}_{\circ_{\frak g}}(x)\mathcal{B}^{\natural}(y),\;\forall x,y\in\frak g.
\end{equation}
\end{pro}

The following result concerns the converse side of Proposition \ref{322}.

\begin{pro}\label{364}
Let $(\frak g,\circ_{\frak g})$	be an anti-pre-Lie algebra and $(\frak g,\lbrace-,-\rbrace_{\frak g})$ be the sub-adjacent Lie algebra of $(\frak g,\circ_{\frak g})$.
Then $\frak g\ltimes_{-L^{*}_{\circ_{\frak g}}}\frak g^{*}$ is a Lie algebra and the natural nondegenerate symmetric bilinear form $\mathcal{B}_{d}$ on $\frak g\oplus\frak g^{*}$ given by
\begin{equation}\label{367}
	\mathcal{B}_{d}(x+a^{*},y+b^{*})=\langle x, b^{*}\rangle+\langle a^{*},y\rangle,\;\forall x,y\in\frak g, a^{*},b^{*}\in\frak g^{*}
\end{equation}
is a commutative $2$-cocycle on $\frak g\ltimes_{-L^{*}_{\circ_{\frak g}}}\frak g^{*}$.
\end{pro}

Recall \cite{Bai2021.2,Bur} that a {\bf pre-Lie algebra} is a pair $(\frak g,\star_{\frak g})$ such that $\frak g$ is a vector space and $\star_{\frak g}:\frak g\otimes\frak g\rightarrow\frak g$ is a multiplication satisfying the following equation:
\begin{equation}\label{388}
		(x\star_{\frak g}y)\star_{\frak g}z-x\star_{\frak g}(y\star_{\frak g}z)=(y\star_{\frak g}x)\star_{\frak g}z-y\star_{\frak g}(x\star_{\frak g}z) , \;\forall x,y,z\in\frak g.
\end{equation}
A {\bf Novikov algebra} \cite{Bal,Gel} is a pre-Lie algebra $(\frak g,\star_{\frak g})$ which moreover satisfies the following equation:
\begin{equation}\label{392}
	(x\star_{\frak g}y)\star_{\frak g}z=(x\star_{\frak g}z)\star_{\frak g}y, \;\forall x,y,z\in\frak g.
\end{equation}

\begin{pro}\label{474}
Suppose that $(\frak g,\star_{\frak g})$ is a pre-Lie algebra.
Let $(\frak g,\circ_{\frak g})$ be the $2$-algebra of $(\frak g,\star_{\frak g})$, that is, the multiplication $\circ_{\frak g}:\frak g\otimes\frak g\rightarrow\frak g$ is given by
\begin{equation}\label{408}
	x\circ_{\frak g} y=x\star_{\frak g} y+2y\star_{\frak g} x,\;\forall x,y\in \frak g.
\end{equation}
Then $(\frak g,\circ_{\frak g})$ is an anti-pre-Lie algebra if and only if
$(\frak g,\star_{\frak g})$ is further a Novikov algebra. 
\end{pro}

\subsection{Compatible Lie algebras with invariant-c$2$c forms}\

\begin{defi}\cite{Gol1,Gol2,OS3}
	Let $(\frak g,[-,-]_{\frak g})$ and $(\frak g,\lbrace-,-\rbrace_{\frak g})$ be Lie algebras. If for all $\lambda,\mu\in \mathbb{K}$, the pair $(\frak g,[-,-]_{\lambda,\mu})$ with the multiplication $[-,-]_{\lambda,\mu}$ given by
	\begin{equation*}
	[x,y]_{\lambda,\mu}=\lambda  [x,y]_{\frak g}+\mu \lbrace x,y\rbrace_{\frak g},\;\forall x,y\in\frak g
	\end{equation*}
	is again a Lie algebra, then  the triple $(\frak g,[-,-]_{\frak g},\lbrace-,-\rbrace_{\frak g})$ is called a {\bf compatible Lie algebra}. 
\end{defi}
 
 \begin{pro}\cite{Str}\label{442}
 	Let $(\frak g,[-,-]_{\frak g})$ and $(\frak g,\lbrace-,-\rbrace_{\frak g})$ be Lie algebras. Then $(\frak g,[-,-]_{\frak g},\lbrace-,-\rbrace_{\frak g})$ is  a  compatible Lie algebra if and only if the following equation holds for all $x,y,z\in\frak g$:
 	\begin{equation}\label{292}
 	 \lbrace[x,y]_{\frak g},z\rbrace_{\frak g}+ \lbrace[y,z]_{\frak g},x\rbrace_{\frak g}+ \lbrace[z,x]_{\frak g},y\rbrace_{\frak g}
 	 +[\lbrace x,y\rbrace_{\frak g},z]_{\frak g}+ [\lbrace y,z\rbrace_{\frak g},x]_{\frak g}+ [\lbrace z,x\rbrace_{\frak g},y]_{\frak g}=0.
 	\end{equation}
 \end{pro}
  
  \begin{defi}
  	Let $(\frak g,[-,-]_{\frak g},\lbrace-,-\rbrace_{\frak g})$  be a  compatible Lie algebra and let $\mathcal{B}$ be a nondegenerate symmetric bilinear form on $\frak g$. We say $\mathcal{B}$ is an {\bf invariant-c$2$c form} on $(\frak g,[-,-]_{\frak g},\lbrace-,-\rbrace_{\frak g})$ 
  	if $\mathcal{B}$ is invariant on $(\frak g,[-,-]_{\frak g})$ in the sense that 
  	\begin{equation}\label{455}
  		\mathcal{B}([x,y]_{\frak g},z)=\mathcal{B}(x,[y,z]_{\frak g}),\;\forall x,y,z\in\frak g
  	\end{equation}
  	and $\mathcal{B}$ is a commutative $2$-cocycle on $(\frak g, \lbrace-,-\rbrace_{\frak g})$.
   \end{defi}

 Recall  \cite{HBG}  that a {\bf quadratic Novikov algebra} is a triple $(\frak g,\star_{\frak g},\mathcal{B})$ in which $(\frak g,\star_{\frak g})$ is a Novikov algebra and $\mathcal{B}$ is a nondegenerate symmetric bilinear form on $\frak g$ which is invariant in the sense that
   \begin{equation}\label{319}
   	\mathcal{B}(x\star_{\frak g}y,z)=-\mathcal{B}(y,x\star_{\frak g}z+z\star_{\frak g} x),\;\forall x,y,z\in\frak g.
   \end{equation}
A {\bf differential Novikov algebra} is a Novikov algebra $(\frak g,\star_{\frak g} )$ together with a linear map $P:\frak g\rightarrow\frak g$ which is a derivation in the sense that
\begin{equation*}
	P(x\star_{\frak g}y)=P(x)\star_{\frak g}y+x\star_{\frak g}P(y),\;\forall x,y\in\frak g.
\end{equation*}
By a {\bf quadratic differential Novikov algebra} we mean a quadruple $(\frak g,\star_{\frak g},\mathcal{B},P)$ in which $(\frak g,\star_{\frak g},\mathcal{B})$ is a quadratic Novikov algebra, $P:\frak g\rightarrow\frak g$ is a derivation of $(\frak g,\star_{\frak g})$ and $P$ is skew-symmetric with respect to $\mathcal{B}$ in the sense that 
\begin{equation}\label{484}
	\mathcal{B}\big(P(x),y\big)+\mathcal{B}\big(x,P(y)\big)=0,
	\;\forall x,y\in\frak g.
\end{equation}
 We show that compatible Lie algebras with invariant-c$2$c forms naturally arise from quadratic differential Novikov algebras.

 \begin{pro}\label{505}
 	Suppose that $(\frak g,\star_{\frak g},P)$ is a differential Novikov algebra.
 	Then there is a compatible Lie algebra $(\frak g,[-,-]_{\frak g},\lbrace-,-\rbrace_{\frak g})$ in which the Lie brackets $[-,-]_{\frak g}$ and $\lbrace-,-\rbrace_{\frak g}$ are respectively given by
 	\begin{align}
 		[x,y]_{\frak g}&=P(x)\star_{\frak g}y-P(y)\star_{\frak g}x,\label{331}\\
 		\lbrace x,y\rbrace_{\frak g}&=y\star_{\frak g}x-x\star_{\frak g}y,\;\forall x,y\in\frak g.\label{332}
 	\end{align}
 	In particular, if $(\frak g,\star_{\frak g},\mathcal{B},P)$ is a quadratic differential Novikov algebra, 
 	then $\mathcal{B}$ is an invariant-c$2$c form on the compatible Lie algebra $(\frak g,[-,-]_{\frak g},\lbrace-,-\rbrace_{\frak g})$ obtained in the above way.
 \end{pro}
 \begin{proof}
 	It is known that $(\frak g,\lbrace-,-\rbrace_{\frak g})$ is a Lie algebra.
 	For all $x,y,z\in\frak g$, we have
 	\begin{align*}
 		[[x,y]_{\frak g},z]_{\frak g}
 		&=[P(x)\star_{\frak g}y-P(y)\star_{\frak g}x,z]_{\frak g}\\
 		&=P\big(P(x)\star_{\frak g}y-P(y)\star_{\frak g}x\big)\star_{\frak g}z
 		-P(z)\star_{\frak g}\big(P(x)\star_{\frak g}y-P(y)\star_{\frak g}x\big)\\
 		&=\big(P^{2}(x)\star_{\frak g}y\big)\star_{\frak g}z
 		+\big(P(x)\star_{\frak g}P(y)\big)\star_{\frak g}z
 		-\big(P^{2}(y)\star_{\frak g}x\big)\star_{\frak g}z
 		-\big(P(y)\star_{\frak g}P(x)\big)\star_{\frak g}z\\
 		&\quad
 		-P(z)\star_{\frak g}\big(P(x)\star_{\frak g}y\big)
 		+P(z)\star_{\frak g}\big(P(y)\star_{\frak g}x\big),\\
 		\ [[y,z]_{\frak g},x]_{\frak g}
 		&=\big(P^{2}(y)\star_{\frak g}z\big)\star_{\frak g}x
 		+\big(P(y)\star_{\frak g}P(z)\big)\star_{\frak g}x
 		-\big(P^{2}(z)\star_{\frak g}y\big)\star_{\frak g}x
 		-\big(P(z)\star_{\frak g}P(y)\big)\star_{\frak g}x\\
 		&\quad
 		-P(x)\star_{\frak g}\big(P(y)\star_{\frak g}z\big)
 		+P(x)\star_{\frak g}\big(P(z)\star_{\frak g}y\big),\\
 		\ [[z,x]_{\frak g},y]_{\frak g}
 		&=\big(P^{2}(z)\star_{\frak g}x\big)\star_{\frak g}y
 		+\big(P(z)\star_{\frak g}P(x)\big)\star_{\frak g}y
 		-\big(P^{2}(x)\star_{\frak g}z\big)\star_{\frak g}y
 		-\big(P(x)\star_{\frak g}P(z)\big)\star_{\frak g}y\\
 		&\quad
 		-P(y)\star_{\frak g}\big(P(z)\star_{\frak g}x\big)
 		+P(y)\star_{\frak g}\big(P(x)\star_{\frak g}z\big).
 	\end{align*}
 	By \eqref{388} and \eqref{392}, we have
 	\begin{equation*}
 		[[x,y]_{\frak g},z]_{\frak g}
 		+[[y,z]_{\frak g},x]_{\frak g}
 		+[[z,x]_{\frak g},y]_{\frak g}=0,
 	\end{equation*}
 	and hence $(\frak g,[-,-]_{\frak g})$ is a Lie algebra. Moreover, we have
 	\begin{align*}
 		\lbrace[x,y]_{\frak g},z\rbrace_{\frak g}
 		&=\lbrace P(x)\star_{\frak g}y-P(y)\star_{\frak g}x,z\rbrace_{\frak g}\\
 		&=-\big(P(x)\star_{\frak g}y\big)\star_{\frak g}z
 		+\big(P(y)\star_{\frak g}x\big)\star_{\frak g}z
 		+z\star_{\frak g}\big(P(x)\star_{\frak g}y\big)
 		-z\star_{\frak g}\big(P(y)\star_{\frak g}x\big),\\
 		\ \lbrace[y,z]_{\frak g},x\rbrace_{\frak g}
 		&=-\big(P(y)\star_{\frak g}z\big)\star_{\frak g}x
 		+\big(P(z)\star_{\frak g}y\big)\star_{\frak g}x
 		+x\star_{\frak g}\big(P(y)\star_{\frak g}z\big)
 		-x\star_{\frak g}\big(P(z)\star_{\frak g}y\big),\\
 		\ \lbrace[z,x]_{\frak g},y\rbrace_{\frak g}
 		&=-\big(P(z)\star_{\frak g}x\big)\star_{\frak g}y
 		+\big(P(x)\star_{\frak g}z\big)\star_{\frak g}y
 		+y\star_{\frak g}\big(P(z)\star_{\frak g}x\big)
 		-y\star_{\frak g}\big(P(x)\star_{\frak g}z\big),\\
 		\ [\lbrace x,y\rbrace_{\frak g},z]_{\frak g}
 		&=[y\star_{\frak g}x-x\star_{\frak g}y,z]_{\frak g}\\
 		&=P\big(y\star_{\frak g}x-x\star_{\frak g}y\big)\star_{\frak g}z
 		-P(z)\star_{\frak g}\big(y\star_{\frak g}x-x\star_{\frak g}y\big)\\
 		&=-\big(P(x)\star_{\frak g}y\big)\star_{\frak g}z
 		-\big(x\star_{\frak g}P(y)\big)\star_{\frak g}z
 		+\big(P(y)\star_{\frak g}x\big)\star_{\frak g}z
 		+\big(y\star_{\frak g}P(x)\big)\star_{\frak g}z\\
 		&\quad
 		+P(z)\star_{\frak g}(x\star_{\frak g}y)
 		-P(z)\star_{\frak g}(y\star_{\frak g}x),\\
 		\ [\lbrace y,z\rbrace_{\frak g},x]_{\frak g}
 		&=-\big(P(y)\star_{\frak g}z\big)\star_{\frak g}x
 		-\big(y\star_{\frak g}P(z)\big)\star_{\frak g}x
 		+\big(P(z)\star_{\frak g}y\big)\star_{\frak g}x
 		+\big(z\star_{\frak g}P(y)\big)\star_{\frak g}x\\
 		&\quad
 		+P(x)\star_{\frak g}(y\star_{\frak g}z)
 		-P(x)\star_{\frak g}(z\star_{\frak g}y),\\
 		\ [\lbrace z,x\rbrace_{\frak g},y]_{\frak g}
 		&=-\big(P(z)\star_{\frak g}x\big)\star_{\frak g}y
 		-\big(z\star_{\frak g}P(x)\big)\star_{\frak g}y
 		+\big(P(x)\star_{\frak g}z\big)\star_{\frak g}y
 		+\big(x\star_{\frak g}P(z)\big)\star_{\frak g}y\\
 		&\quad
 		+P(y)\star_{\frak g}(z\star_{\frak g}x)
 		-P(y)\star_{\frak g}(x\star_{\frak g}z).
 	\end{align*}
 	By \eqref{388} and \eqref{392} again, we conclude that \eqref{292} holds, and hence $(\frak g,[-,-]_{\frak g},\lbrace-,-\rbrace_{\frak g})$ is a compatible Lie algebra.
 	By \cite{LB2022}, $\mathcal{B}$ is a commutative $2$-cocycle on $(\frak g,\lbrace-,-\rbrace_{\frak g})$.
 	Moreover, we have
 	\begin{align*}
 		\mathcal{B}([x,y]_{\frak g},z)
 		&=\mathcal{B}\big(P(x)\star_{\frak g}y-P(y)\star_{\frak g}x,z\big)
 		\overset{\eqref{319}}{=}
 		\mathcal{B}\big(P(x),z\star_{\frak g}y\big)
 		+\mathcal{B}\big(x,P(y)\star_{\frak g}z+z\star_{\frak g}P(y)\big)\\
 		&=-\mathcal{B}\big(x,P(z\star_{\frak g}y)\big)
 		+\mathcal{B}\big(x,P(y)\star_{\frak g}z+z\star_{\frak g}P(y)\big)
 		=\mathcal{B}\big(x,P(y)\star_{\frak g}z-P(z)\star_{\frak g}y\big)\\
 		&=\mathcal{B}(x,[y,z]_{\frak g}).
 	\end{align*}
 	Therefore $\mathcal{B}$ is invariant on $(\frak g,[-,-]_{\frak g})$. Hence $\mathcal{B}$ is an invariant-c$2$c form on
 	 $(\frak g,[-,-]_{\frak g},\lbrace-,-\rbrace_{\frak g})$.
 \end{proof}
  
 \begin{rmk}
 The bracket $\lbrace-,-\rbrace_{\frak g}$ in \eqref{332}
 is the opposite of the usual commutator bracket of
 $(\frak g,\star_{\frak g})$.
 This sign convention is chosen so that,
 	when $(\frak g,\star_{\frak g},\mathcal{B},P)$ is a
 	quadratic differential Novikov algebra,
 	the anti-pre-Lie algebra induced by $\mathcal{B}$
 	through Proposition \ref{322} is exactly the $2$-algebra
 	of $(\frak g,\star_{\frak g})$, see Proposition \ref{751}.
 \end{rmk}
  
  \begin{ex}\label{584}
  	Let $\frak g$ be a vector space with a multiplication
  	$\star_{\frak g}$.
  	Define a multiplication $[-,-]_{\frak h}$ on
  	$\frak h=\frak g\otimes\mathbb{K}[t,t^{-1}]$ by
  	\begin{equation*}
  		[x\otimes t^{i},y\otimes t^{j}]_{\frak h}
  		=i(x\star_{\frak g}y)\otimes t^{i+j-1}
  		-j(y\star_{\frak g}x)\otimes t^{i+j-1},
  		\;\forall x,y\in\frak g,\;i,j\in\mathbb{Z}.
  	\end{equation*}
  	By \cite{Bal}, $(\frak h,[-,-]_{\frak h})$ is a Lie algebra
  	if and only if $(\frak g,\star_{\frak g})$ is a Novikov algebra.
  	
  	Suppose further that $(\frak g,\star_{\frak g},\mathcal{B})$
  	is a quadratic Novikov algebra.
  	Define a multiplication $\star_{\frak h}$, a bilinear form
  	$\overline{\mathcal{B}}$ and a linear map
  	$P:\frak h\rightarrow\frak h$ by
  	\begin{align*}
  		(x\otimes t^{i})\star_{\frak h}(y\otimes t^{j})
  		&=(x\star_{\frak g}y)\otimes t^{i+j},\\
  		\overline{\mathcal{B}}(x\otimes t^{i},y\otimes t^{j})
  		&=\delta_{i+j,-1}\mathcal{B}(x,y),\\
  		P(x\otimes t^{i})
  		&=i x\otimes t^{i-1},
  	\end{align*}
  	for all $x,y\in\frak g$ and $i,j\in\mathbb{Z}$,
  	where $\delta_{m,n}$ is the Kronecker delta.
  	Then it is straightforward to check that
  	$(\frak h,\star_{\frak h},\overline{\mathcal{B}},P)$
  	is an (infinite-dimensional) quadratic differential Novikov algebra.
  	Consequently, there is a compatible Lie algebra
  	$(\frak h,[-,-]_{\frak h},\lbrace-,-\rbrace_{\frak h})$,
  	where $[-,-]_{\frak h}$ is given above and
  	\begin{equation*}
  		\lbrace x\otimes t^i,y\otimes t^j\rbrace_{\frak h}
  		=(y\star_{\frak g}x-x\star_{\frak g}y)\otimes t^{i+j},
  		\;\forall x,y\in\frak g,\;i,j\in\mathbb{Z}.
  	\end{equation*}
  	Furthermore, $\overline{\mathcal{B}}$ is an invariant-c$2$c
  	form on
  	$(\frak h,[-,-]_{\frak h},\lbrace-,-\rbrace_{\frak h})$.
  \end{ex}
  
  We next give an example of a compatible Lie algebra with an invariant-c$2$c form on $\frak s\frak l(2,\mathbb{K})$.
  
\begin{ex}\label{631}
		Let $(\frak g,\lbrace-,-\rbrace_{\frak g})$ be the $3$-dimensional simple Lie algebra
	$\frak s\frak l(2,\mathbb{K})$ with a basis 
	\begin{equation*}
		\{
		x=\left(\begin{matrix} 0&1\cr
			0&0\cr\end{matrix}\right),h=\left(\begin{matrix} 1&0\cr
			0&-1\cr\end{matrix}\right),y=\left(\begin{matrix} 0&0\cr
			1&0\cr\end{matrix}\right)\}
	\end{equation*}
	whose products are given by
	\begin{equation*}
		\lbrace h,x \rbrace_{\frak g}=2x,\;
		\lbrace h,y \rbrace_{\frak g}=-2y,\;
		\lbrace x,y\rbrace_{\frak g}=h.
	\end{equation*}
By \cite{LB2022}, the bilinear form $\mathcal B$
	with non-zero values given by
	\begin{equation*}
		\mathcal{B}(x,y)=\mathcal{B}(y,x)=-1,\;\mathcal{B}(h,h)=4
	\end{equation*}
	is a nondegenerate commutative $2$-cocycle on $(\frak g,\lbrace-,-\rbrace_{\frak g})$.
	By Proposition \ref{322}, there is an anti-pre-Lie algebra $(\frak g,\circ_{\frak g})$
	given by the following non-zero products:
	\begin{equation*}
		x\circ_{\frak g} h=-4x,\; x\circ_{\frak g} y=\dfrac{1}{2}h,\; h\circ_{\frak g} x=-2x,\; h\circ_{\frak g} y=2y,\; y\circ_{\frak g} x=-\dfrac{1}{2}h,\; y\circ_{\frak g} h=4y 
	\end{equation*}
	such that \eqref{297} holds. 
Define another multiplication  $[-,-]_{\frak g}:\frak g\otimes\frak g\rightarrow\frak g$ by the following nonzero products:
\begin{equation*}
[h,x]_{\frak g}=-[x,h]_{\frak g}=-4x,\;
[h,y]_{\frak g}=-[y,h]_{\frak g}=4y,\;
[x,y]_{\frak g}=-[y,x]_{\frak g}=h.
\end{equation*}
Then $(\frak g,[-,-]_{\frak g})$ is a Lie algebra which is compatible with $(\frak g ,\lbrace-,-\rbrace_{\frak g})$.
Moreover, it is straightforward to check that $\mathcal B$ is a nondegenerate symmetric invariant bilinear form on  $(\frak g,[-,-]_{\frak g})$.
Therefore $\mathcal B$ is an invariant-c$2$c form on the compatible Lie algebra $(\frak g ,[-,-]_{\frak g},\lbrace-,-\rbrace_{\frak g})$.
\end{ex}
  
  \subsection{Compatible LAPL algebras}\
  Recall the representations of compatible Lie algebras.
  
  \begin{defi}\cite{WB}
  	Let $(\frak g,[-,-]_{\frak g},\lbrace-,-\rbrace_{\frak g})$ be a compatible Lie algebra.
  	Suppose that $(\rho,V)$ is a representation of $(\frak g,[-,-]_{\frak g} )$ and $(\mu,V)$ is a representation of $(\frak g, \lbrace-,-\rbrace_{\frak g})$. If the following equation holds:
  	\begin{equation}\label{652}
  		\mu([x,y]_{\frak g})v-\mu(x)\rho(y)v+\mu(y)\rho(x)v+\rho(\lbrace x,y\rbrace_{\frak g})v-\rho(x)\mu(y)v+\rho(y)\mu(x)v=0,\;\forall x,y\in\frak g, v\in V,
  	\end{equation}
  then we say $(\rho,\mu,V)$ is a {\bf representation} of $(\frak g,[-,-]_{\frak g},\lbrace-,-\rbrace_{\frak g})$.
  \end{defi}
  
  \begin{ex}
  Let $(\frak g,[-,-]_{\frak g},\lbrace-,-\rbrace_{\frak g})$ be a compatible Lie algebra.
  Then $(\mathrm{ad}_{\frak g},\mathfrak{ad}_{\frak g},\frak g)$ is a representation of $(\frak g,[-,-]_{\frak g},\lbrace-,-\rbrace_{\frak g})$, which is called the {\bf adjoint representation} of $(\frak g,[-,-]_{\frak g},\lbrace-,-\rbrace_{\frak g})$.
  \end{ex}
  
  \begin{pro}\cite{WB}
 Let $(\frak g,[-,-]_{\frak g},\lbrace-,-\rbrace_{\frak g})$ be a compatible Lie algebra, $V$ be a vector space and $\rho,\mu:\frak g\rightarrow\mathrm{End}_{\mathbb{K}}(V)$ be linear maps.
 Define multiplications $[-,-]_{d},\lbrace-,-\rbrace_{d}$ on $\frak g\oplus V$ by
 \begin{align}
   [x+u,y+v]_{d}&=[x,y]_{\frak g}+\rho(x)v-\rho(y)u,\label{667}\\
   \lbrace x+u,y+v\rbrace_{d}&=\lbrace x,y\rbrace_{\frak g}+ \mu(x)v-\mu(y)u,\;\forall x,y\in\frak g, u,v\in V.
\notag \end{align}
Then $(\rho,\mu,V)$ is a  representation of $(\frak g,[-,-]_{\frak g},\lbrace-,-\rbrace_{\frak g})$ if and only if $(\frak g\oplus V,[-,-]_{d},\lbrace-,-\rbrace_{d})=\frak g \ltimes_{\rho,\mu}V$ is a compatible Lie algebra.
  \end{pro}

 \begin{pro}
 Let $(\frak g,[-,-]_{\frak g},\lbrace-,-\rbrace_{\frak g})$ be a compatible Lie algebra and $(\rho,\mu,V)$ be a  representation of $(\frak g,[-,-]_{\frak g},\lbrace-,-\rbrace_{\frak g})$.
 Then  $(\rho^{*},\mu^{*},V^{*})$ is also a representation of $(\frak g,[-,-]_{\frak g},\lbrace-,-\rbrace_{\frak g})$.
 \end{pro}
  
 \begin{defi}
 	Let $(\frak g,[-,-]_{\frak g})$ be a Lie algebra, $(\frak g,\circ_{\frak g})$ be an anti-pre-Lie algebra and $(\frak g, \lbrace-,-\rbrace_{\frak g})$ be the sub-adjacent Lie algebra of $(\frak g,\circ_{\frak g})$.
 	If $(\frak g,[-,-]_{\frak g},\lbrace-,-\rbrace_{\frak g})$ is a compatible Lie algebra and the following equation holds for all $x,y,z\in\frak g$:
 	\begin{equation}\label{646}
 	[x,y]_{\frak g}\circ_{\frak g}z-x\circ_{\frak g}[y,z]_{\frak g}+y\circ_{\frak g}[x,z]_{\frak g}
 	-[\lbrace x,y\rbrace_{\frak g},z]_{\frak g}-[x,y\circ_{\frak g} z]_{\frak g}+[y,x\circ_{\frak g}z]_{\frak g}=0,
 	\end{equation}
 then we say $(\frak g,[-,-]_{\frak g},\circ_{\frak g})$ is a {\bf compatible Lie-anti-pre-Lie algebra} or a {\bf compatible LAPL algebra} for short. Moreover, we say $(\frak g,[-,-]_{\frak g},\lbrace-,-\rbrace_{\frak g})$ is the {\bf sub-adjacent compatible Lie algebra} of $(\frak g,[-,-]_{\frak g},\circ_{\frak g})$.
 \end{defi}
  
  \begin{pro}
 Let $(\frak g,[-,-]_{\frak g})$ be a Lie algebra, $(\frak g,\circ_{\frak g})$ be an anti-pre-Lie algebra and $(\frak g, \lbrace-,-\rbrace_{\frak g})$ be the sub-adjacent Lie algebra of $(\frak g,\circ_{\frak g})$. Then the following conditions are equivalent:
 \begin{enumerate}
 	\item $(\frak g,[-,-]_{\frak g},\circ_{\frak g})$ is a  compatible LAPL algebra.
 	\item $(\frak g,[-,-]_{\frak g},\lbrace-,-\rbrace_{\frak g})$ is a compatible Lie algebra with a representation $(\mathrm{ad}_{\frak g},-L_{\circ_{\frak g}},\frak g)$.
 	\item $(\frak g,[-,-]_{\frak g},\lbrace-,-\rbrace_{\frak g})$ is a compatible Lie algebra with a representation $(\mathrm{ad}^{*}_{\frak g},-L^{*}_{\circ_{\frak g}},\frak g^{*})$.
 \end{enumerate}
  \end{pro}
  \begin{proof}
  	It follows from a straightforward check.
  \end{proof}
  
  Let $(\frak g,[-,-]_{\frak g},\lbrace-,-\rbrace_{\frak g})$ be a compatible Lie algebra and  $(\rho,\mu,V)$ and $(\rho',\mu',V')$  be two representations of $(\frak g,[-,-]_{\frak g}$,
  $\lbrace-,-\rbrace_{\frak g})$.
  If there exists a linear isomorphism $\phi:V\rightarrow V'$  such that \eqref{324} and the following equation hold:
  \begin{equation*}
  	\phi\big( \rho(x)v\big)=\rho'(x)\phi(v),\;\forall x\in\frak g, v\in V,
  \end{equation*}
  then we say $(\rho,\mu,V)$ and $(\rho',\mu',V')$ are {\bf equivalent} through $\phi$.
  
  \begin{pro}\label{673}
  	Let $(\frak g,[-,-]_{\frak g},\lbrace-,-\rbrace_{\frak g})$ be a compatible Lie algebra and 
  	 let $\mathcal{B}$ be an invariant-c$2$c form on  $(\frak g,[-,-]_{\frak g},\lbrace-,-\rbrace_{\frak g})$.
  	Then there is a compatible LAPL algebra $(\frak g,[-,-]_{\frak g},\circ_{\frak g})$
    with $\circ_{\frak g}:\frak g\otimes\frak g\rightarrow\frak g$ given by \eqref{332-}
  	such that \eqref{297} holds,
  	that is, $(\frak g,[-,-]_{\frak g},\lbrace-,-\rbrace_{\frak g})$ is the sub-adjacent compatible Lie algebra of $(\frak g,[-,-]_{\frak g},\circ_{\frak g})$.
  	 Moreover, the adjoint representation $(\mathrm{ad}_{\frak g},\mathfrak{ad}_{\frak g}, \frak g)$
  	of $ (\frak g,[-,-]_{\frak g},\lbrace-,-\rbrace_{\frak g}) $
  	is equivalent to $(\mathrm{ad}^{*}_{\frak g},-L^{*}_{\circ_{\frak g}},\frak g^{*})$ through $\mathcal{B}^{\natural}$.
  \end{pro}
  \begin{proof}
  	By Proposition \ref{322}, there is an anti-pre-Lie algebra $(\frak g, \circ_{\frak g})$
  	with $\circ_{\frak g}:\frak g\otimes\frak g\rightarrow\frak g$ given by \eqref{332-}
  	such that \eqref{297} holds.
  	Moreover, for all $x,y,z,t\in\frak g$, we have
  	\begin{align*}
&\mathcal{B}([x,y]_{\frak g}\circ_{\frak g}z,t)   \overset{\eqref{332-}} {=} \mathcal{B}(z,\lbrace [x,y]_{\frak g}, t\rbrace_{\frak g}),\\
-&\mathcal{B}(x\circ_{\frak g}[y,z]_{\frak g},t)   \overset{\eqref{332-}} {=} -\mathcal{B}([y,z]_{\frak g},\lbrace x,t\rbrace_{\frak g}) \overset{\eqref{455}}{=}\mathcal{B}(z,[y,\lbrace x,t\rbrace_{\frak g}]_{\frak g}),\\
&\mathcal{B}(y\circ_{\frak g}[x,z]_{\frak g},t)   = -\mathcal{B}(z,[x,\lbrace y,t\rbrace_{\frak g}]_{\frak g}),\\
-&\mathcal{B}( [\lbrace x,y\rbrace_{\frak g},z]_{\frak g},t  )   \overset{\eqref{455}} {=} \mathcal{B}(z,[\lbrace x,y\rbrace_{\frak g},t]_{\frak g}),\\
-&\mathcal{B}([x,y\circ_{\frak g}z]_{\frak g},t) \overset{\eqref{455}}{=} \mathcal{B}(y\circ_{\frak g}z, [x,t]_{\frak g})\overset{\eqref{332-}}{=}\mathcal{B}(z,\lbrace y,[x,t]_{\frak g}\rbrace_{\frak g}),\\
 &\mathcal{B}([y,x\circ_{\frak g}z]_{\frak g},t) = -\mathcal{B}(z,\lbrace x,[y,t]_{\frak g}\rbrace_{\frak g}),
  	\end{align*}
  Therefore we have
  \begin{align*}
  	& \mathcal{B}([x,y]_{\frak g}\circ_{\frak g}z-x\circ_{\frak g}[y,z]_{\frak g}+y\circ_{\frak g}[x,z]_{\frak g}
  	-[\lbrace x,y\rbrace_{\frak g},z]_{\frak g}-[x,y\circ_{\frak g} z]_{\frak g}+[y,x\circ_{\frak g}z]_{\frak g},t)\\
  	&= \mathcal{B}(z,\lbrace [x,y]_{\frak g}, t\rbrace_{\frak g}+[y,\lbrace x,t\rbrace_{\frak g}]_{\frak g}-[x,\lbrace y,t\rbrace_{\frak g}]_{\frak g}+[\lbrace x,y\rbrace_{\frak g},t]_{\frak g}+\lbrace y,[x,t]_{\frak g}\rbrace_{\frak g}-\lbrace x,[y,t]_{\frak g}\rbrace_{\frak g}  )\\
  	&\overset{\eqref{292}}{=} 0.
  \end{align*}
By the nondegeneracy of $\mathcal{B}$, \eqref{646} holds and hence $(\frak g,[-,-]_{\frak g},\circ_{\frak g})$ is a compatible LAPL algebra. Moreover, by Proposition \ref{322}, \eqref{344} holds. By \eqref{455}, the following equation holds:
\begin{equation*}
\mathcal{B}^{\natural}\big( \mathrm{ad}_{\frak g}(x)y  \big)=\mathrm{ad}^{*}_{\frak g}(x) \mathcal{B}^{\natural}(y),\;\forall x,y\in\frak g.
\end{equation*}
Therefore the representations $(\mathrm{ad}_{\frak g},\mathfrak{ad}_{\frak g}, \frak g)$
and  $(\mathrm{ad}^{*}_{\frak g},-L^{*}_{\circ_{\frak g}},\frak g^{*})$ of $(\frak g,[-,-]_{\frak g},\lbrace-,-\rbrace_{\frak g})$ are equivalent through $\mathcal{B}^{\natural}$.
  \end{proof}
  
  \begin{pro}\label{751}
  	Let $(\frak g,\star_{\frak g},P)$ be a differential Novikov algebra and $(\frak g,[-,-]_{\frak g})$ be the Lie algebra with $[-,-]_{\frak g}$ given by \eqref{331}.
  	Let $\circ_{\frak g}:\frak g\otimes\frak g\rightarrow\frak g$ be the multiplication given by \eqref{408}.
  	Then $(\frak g,[-,-]_{\frak g},\circ_{\frak g})$ is a compatible LAPL algebra.
  	Furthermore, suppose that there exists a nondegenerate symmetric invariant bilinear form $\mathcal{B}$ on $(\frak g,\star_{\frak g})$ satisfying \eqref{484}, which is consequently an invariant-c$2$c form on the compatible Lie algebra $(\frak g,[-,-]_{\frak g},\lbrace-,-\rbrace_{\frak g})$ with $\lbrace-,-\rbrace_{\frak g}$ given by \eqref{332}.
  	Then $(\frak g,[-,-]_{\frak g},\circ_{\frak g})$ is exactly the compatible LAPL algebra induced from the invariant-c$2$c form $\mathcal{B}$ on the compatible Lie algebra
  	$(\frak g,[-,-]_{\frak g},\lbrace-,-\rbrace_{\frak g})$
  	given in Proposition \ref{505}
  	through Proposition \ref{673}.
  \end{pro}
  \begin{proof}
  	By Proposition \ref{474}, $(\frak g,\circ_{\frak g})$ is an anti-pre-Lie algebra.
  	Moreover, it is straightforward to check that \eqref{646} holds.
  	Therefore
  	$(\frak g,[-,-]_{\frak g},\circ_{\frak g})$ is a compatible LAPL algebra.
  	Furthermore, the symmetry of $\mathcal{B}$ and \eqref{319} imply
  	\begin{equation*}
  		\mathcal{B}(x\circ_{\frak g}y,z)
  		=\mathcal{B}(x\star_{\frak g}y+2y\star_{\frak g}x,z)
  		=\mathcal{B}(y,z\star_{\frak g}x-x\star_{\frak g}z)
  		=\mathcal{B}(y,\lbrace x,z\rbrace_{\frak g}),
  		\;\forall x,y,z\in\frak g.
  	\end{equation*}
  	Hence the conclusion follows from Proposition \ref{673}.
  \end{proof}
  
  \begin{ex}
  	Let $(\frak h,\star_{\frak h},\overline{\mathcal{B}},P)$
  	be the quadratic differential Novikov algebra given in
  	Example \ref{584}.
  	By Proposition \ref{751}, there is a compatible LAPL algebra
  	$(\frak h,[-,-]_{\frak h},\circ_{\frak h})$,
  	where $[-,-]_{\frak h}$ is given in Example \ref{584}
  	and $\circ_{\frak h}$ is given by
  	\begin{equation*}
  		\begin{aligned}
  			(x\otimes t^{i})\circ_{\frak h}(y\otimes t^{j})
  			&=(x\otimes t^{i})\star_{\frak h}(y\otimes t^{j})
  			+2(y\otimes t^{j})\star_{\frak h}(x\otimes t^{i})\\
  			&=(x\star_{\frak g}y+2y\star_{\frak g}x)\otimes t^{i+j},
  		\end{aligned}
  	\end{equation*}
  for all $x,y\in\frak g,\;i,j\in\mathbb{Z}$.
  	Furthermore, $(\frak h,[-,-]_{\frak h},\circ_{\frak h})$
  	is exactly the compatible LAPL algebra induced from
  	the invariant-c$2$c form $\overline{\mathcal{B}}$ on
  	the compatible Lie algebra
  	$(\frak h,[-,-]_{\frak h},\lbrace-,-\rbrace_{\frak h})$
  	given in Example \ref{584}.
  \end{ex}

  \begin{ex}
Keep the notations in Example \ref{631}.
By Proposition \ref{673}, $(\frak g,[-,-]_{\frak g},\circ_{\frak g})$ is a compatible LAPL algebra, which is induced from the invariant-c$2$c form $\mathcal{B}$ on
$(\frak g ,[-,-]_{\frak g},\lbrace-,-\rbrace_{\frak g})$.
  \end{ex}
  
  The following result concerns the converse side of Proposition \ref{673}.

  	\begin{pro}\label{794}
  		Let $(\frak g,[-,-]_{\frak g},\circ_{\frak g})$	be a compatible LAPL algebra and $(\frak g,[-,-]_{\frak g},\lbrace-,-\rbrace_{\frak g})$ be the sub-adjacent compatible Lie algebra of $(\frak g,[-,-]_{\frak g},\circ_{\frak g})$.
  		Then $\frak g\ltimes_{\mathrm{ad}_{\frak g}^{*},-L^{*}_{\circ_{\frak g}}}\frak g^{*}$ is a compatible Lie algebra and the natural nondegenerate symmetric bilinear form $\mathcal{B}_{d}$ on $\frak g\oplus\frak g^{*}$ given by \eqref{367}
  		is an invariant-c$2$c form on $\frak g\ltimes_{\mathrm{ad}_{\frak g}^{*},-L^{*}_{\circ_{\frak g}}}\frak g^{*}$.
  \end{pro}
  \begin{proof}
  	It follows from a  straightforward check.
  \end{proof}
   
  Recall the representation theory of anti-pre-Lie algebras.
  
  \begin{defi}\cite{TPA} 
  	A {\bf representation} of an anti-pre-Lie algebra $(\frak g,\circ_{\frak g})$ is a triple $(l_{\circ},r_{\circ},V)$ in which $V$ is a vector space and $l_{\circ},r_{\circ}:\frak g\rightarrow\mathrm{End}_{\mathbb{K}}(V)$ are linear maps satisfying the following equations:
  	\begin{align*}
  	 	r_{\circ} (x\circ_{\frak g} y)v&=l_{\circ} (x)r_{\circ} (y)v+r_{\circ} (y)l_{\circ} (x)v-r_{\circ} (y)r_{\circ} (x)v,\\
  	 	l_{\circ} (y\circ_{\frak g} x)v-l_{\circ} (x\circ_{\frak g} y)v&=l_{\circ} (x)l_{\circ} (y)v-l_{\circ} (y)l_{\circ} (x)v,\\
  		 l_{\circ} (y\circ_{\frak g} x)v-l_{\circ} (x\circ_{\frak g} y)v&=r_{\circ} (x)l_{\circ} (y)v-r_{\circ} (y)l_{\circ} (x)v-r_{\circ} (x)r_{\circ} (y)v+r_{\circ} (y)r_{\circ} (x)v,\;\forall x,y\in\frak g, v\in V.
  	\end{align*}
  \end{defi}
  
  \begin{pro}\cite{TPA} 
  	Suppose that $(\frak g,\circ_{\frak g})$ is an anti-pre-Lie algebra,
  	 $V$ is a vector space and $l_{\circ},r_{\circ}:\frak g\rightarrow\mathrm{End}_{\mathbb{K}}(V)$ are linear maps.
  	 Then $(l_{\circ},r_{\circ}, V)$ is a  representation  of $(\frak g,\circ_{\frak g} )$ if and only if there is an anti-pre-Lie algebra structure $(\frak g\oplus V,\circ_{d} )=\frak g\ltimes_{l_{\circ},r_{\circ}}V$ on $\frak g\oplus V$ with $\circ_{d}$
  	 given by
  	 \begin{eqnarray*}
  	 	(x+u)\circ_{d}(y+v)=x\circ_{\frak g}y+l_{\circ}(x)v+r_{\circ}(y)u,\;\forall x,y\in\frak g,u,v\in V.
  	 \end{eqnarray*}
  \end{pro}

\begin{lem}\cite{TPA} 
	Let $(\frak g,\circ_{\frak g} )$ be an anti-pre-Lie algebra and $(\frak g,\lbrace-,-\rbrace_{\frak g})$ be the sub-adjacent Lie algebra of $(\frak g,\circ_{\frak g} )$.
	Suppose that $(l_{\circ},r_{\circ},V)$ is a representation of $(\frak g,\circ_{\frak g})$.
	Then $( r_{\circ}^{*}-l_{\circ}^{*}, r_{\circ}^{*},V^{*})$ is a representation of $(\frak g,\circ_{\frak g})$.
	Moreover,  $(-l_{\circ},V),(-l_{\circ}^{*},V^{*}), (l_{\circ}-r_{\circ},V)$ and $(l_{\circ}^{*}-r_{\circ}^{*},V^{*})$ are representations of $(\frak g,\lbrace-,-\rbrace_{\frak g})$.
\end{lem}
 
  Now we shall study the representation theory of compatible LAPL algebras.
  
  \begin{defi}
  	Let $(\frak g,[-,-]_{\frak g},\circ_{\frak g})$ be a compatible LAPL algebra and $(\frak g,[-,-]_{\frak g},\lbrace-,-\rbrace_{\frak g})$ be the sub-adjacent compatible Lie algebra.
  	Suppose that $(\rho,V)$ is a representation of $(\frak g,[-,-]_{\frak g})$
  	and $(l_{\circ},r_{\circ},V)$ is a representation of $(\frak g, \circ_{\frak g})$.
  	Assume that $(\rho,l_{\circ}-r_{\circ},V)$ is a representation of $(\frak g,[-,-]_{\frak g},\lbrace-,-\rbrace_{\frak g})$, that is, \eqref{652} holds for $\mu=l_{\circ}-r_{\circ}$.
  	If the following equations hold for all $x,y\in\frak g, v\in V$:
  	\begin{align}
 l_{\circ}([x,y]_{\frak g})v-l_{\circ}(x)\rho(y)v+l_{\circ}(y)\rho(x)v-\rho(\lbrace x,y\rbrace_{\frak g})v-\rho(x)l_{\circ}(y)v+\rho(y)l_{\circ}(x)v&=0,\label{838}\\
  r_{\circ}(y)\rho(x)v+l_{\circ}(x)\rho(y)v+r_{\circ}([x,y]_{\frak g})v+\rho(y)(l_{\circ}-r_{\circ})(x)v-\rho(x)r_{\circ}(y)v-\rho(x\circ_{\frak g}y)v&=0, \label{839}
  	\end{align}
  then we say $(\rho,l_{\circ},r_{\circ},V)$ is a {\bf representation} of $(\frak g,[-,-]_{\frak g},\circ_{\frak g})$.
   \end{defi}
  
  \begin{pro}
  	Let $(\frak g,[-,-]_{\frak g},\circ_{\frak g})$ be a compatible LAPL algebra, $V$ be a vector space and $\rho,l_{\circ},r_{\circ}:\frak g\rightarrow\mathrm{End}_{\mathbb{K}}(V)$ be linear maps. Define multiplications $[-,-]_{d}$ and $\circ_{d}$ respectively by \eqref{667} and the following equation:
  	\begin{equation*}
  		(x+u)\circ_{d}(y+v)=x\circ_{\frak g}y+l_{\circ}(x)v+r_{\circ}(y)u,\;\forall x,y\in\frak g, u,v\in V.
  	\end{equation*}
  Then $(\rho,l_{\circ},r_{\circ},V)$ is a  representation of $(\frak g,[-,-]_{\frak g},\circ_{\frak g})$ if and only if $(\frak g\oplus V,[-,-]_{d},\circ_{d})=\frak g\ltimes_{\rho,l_{\circ},r_{\circ}}V$ is a compatible LAPL algebra.
  \end{pro}
  \begin{proof}
  	It follows from a straightforward computation.
  \end{proof}
  
  \begin{ex}
	Let $(\frak g,[-,-]_{\frak g},\circ_{\frak g})$ be a compatible LAPL algebra.
	Then $( \mathrm{ad} _{\frak g}, L _{\circ_{\frak g}},R _{\circ_{\frak g}},\frak g  )$
	 is a representation of $(\frak g,[-,-]_{\frak g},\circ_{\frak g})$, which is called the {\bf adjoint representation} of $(\frak g,[-,-]_{\frak g},\circ_{\frak g})$.
  \end{ex}

 \begin{pro}
	Let $(\frak g,[-,-]_{\frak g},\circ_{\frak g})$ be a compatible LAPL algebra and $(\rho,l_{\circ},r_{\circ},V)$ be a representation of $(\frak g,[-,-]_{\frak g},\circ_{\frak g})$.
	Then $(\rho^{*},r_{\circ}^{*}-l_{\circ}^{*},r_{\circ}^{*},V^{*})$ is a representation of $(\frak g,[-,-]_{\frak g},\circ_{\frak g})$.
	In particular, $(\mathrm{ad}^{*}_{\frak g},- \mathfrak{ad}^{*}_{\frak g},R^{*}_{\circ_{\frak g}},\frak g^{*})$ is a representation of $(\frak g,[-,-]_{\frak g},\circ_{\frak g})$.
 \end{pro}
  \begin{proof}
  	We already know that $(\rho^{*}, V^{*})$ is a representation of $(\frak g,[-,-]_{\frak g} )$ and $( r_{\circ}^{*}-l_{\circ}^{*},r_{\circ}^{*},V^{*})$ is a representation of $(\frak g, \circ_{\frak g})$.
  	For all $x,y\in\frak g, u^{*}\in V^{*}, v\in V$, we have
\begin{align*}
& \langle -l_{\circ}^{*}([x,y]_{\frak g})u^{*}+l_{\circ}^{*}(x)\rho^{*}(y)u^{*}-l_{\circ}^{*}(y)\rho^{*}(x)u^{*}
+\rho^{*}(\lbrace x,y\rbrace_{\frak g})u^{*}+\rho^{*}(x)l_{\circ}^{*}(y)u^{*}-\rho^{*}(y)l_{\circ}^{*}(x)u^{*}, v\rangle\\
&= \langle u^{*}, l_{\circ}([x,y]_{\frak g})v+\rho(y)l_{\circ}(x)v-\rho(x)l_{\circ}(y)v-\rho(\lbrace x,y\rbrace_{\frak g})v+l_{\circ}(y)\rho(x)v-l_{\circ}(x)\rho(y)v\rangle\\
&\overset{\eqref{838}}{=} 0.
\end{align*}
Therefore, \eqref{652} holds for the triple $(\rho^{*}, r_{\circ}^{*}-l_{\circ}^{*}-r_{\circ}^{*}=-l_{\circ}^{*}, V^{*})$.
Similarly, \eqref{838} and \eqref{839} hold for the quadruple 
$(\rho^{*}, r_{\circ}^{*}-l_{\circ}^{*}, r_{\circ}^{*}, V^{*})$. 	Hence $(\rho^{*},r_{\circ}^{*}-l_{\circ}^{*},r_{\circ}^{*},V^{*})$ is a representation of $(\frak g,[-,-]_{\frak g},\circ_{\frak g})$. The remaining result is obvious.
  \end{proof}
  
  At the end of this subsection, we shall introduce an identity that holds on compatible LAPL algebras, which we shall use later.
  
  \begin{pro}
  	Let $(\frak g,[-,-]_{\frak g},\circ_{\frak g})$ be a compatible LAPL algebra.
  	Then the following identity holds for all $x,y,z\in\frak g$:
  	\begin{equation}\label{886}
  		[x\circ_{\frak g}y,z]_{\frak g}-[x\circ_{\frak g}z,y]_{\frak g}+x\circ_{\frak g}[z,y]_{\frak g}+[x,y]_{\frak g}\circ_{\frak g}z-[x,z]_{\frak g}\circ_{\frak g}y+2[x,\lbrace z,y\rbrace_{\frak g}]_{\frak g}=0.
  	\end{equation}
  \end{pro}
  \begin{proof}
  	For all $x,y,z\in\frak g$, we have
  	\begin{align*}
  			&[x\circ_{\frak g}y,z]_{\frak g}-[x\circ_{\frak g}z,y]_{\frak g}+x\circ_{\frak g}[z,y]_{\frak g}+[x,y]_{\frak g}\circ_{\frak g}z-[x,z]_{\frak g}\circ_{\frak g}y+2[x,\lbrace z,y\rbrace_{\frak g}]_{\frak g}\\
  		&	=
  	\lbrace[x,y]_{\frak g},z\rbrace_{\frak g}+ \lbrace[y,z]_{\frak g},x\rbrace_{\frak g}+ \lbrace[z,x]_{\frak g},y\rbrace_{\frak g}
  	+[\lbrace x,y\rbrace_{\frak g},z]_{\frak g}+ [\lbrace y,z\rbrace_{\frak g},x]_{\frak g}+ [\lbrace z,x\rbrace_{\frak g},y]_{\frak g}\\
  	&\quad-[y,z]_{\frak g}\circ_{\frak g}x+y\circ_{\frak g}[z,x]_{\frak g}-z\circ_{\frak g}[y,x]_{\frak g}
  	+[\lbrace y,z\rbrace_{\frak g},x]_{\frak g}+[y,z\circ_{\frak g} x]_{\frak g}-[z,y\circ_{\frak g}x]_{\frak g}\\
  	&\overset{\eqref{292},\eqref{646}}{=}0.
  	\end{align*}
  Hence the conclusion holds.
  \end{proof}

  \section{Manin triples of compatible Lie algebras with respect to the invariant-c$2$c forms and compatible LAPL bialgebras}
  
  In this section, we introduce the notions of Manin triples of compatible Lie algebras with respect to the invariant-c$2$c forms and compatible LAPL bialgebras, and show that they are equivalent structures.
  
  \subsection{Manin triples of compatible Lie algebras with respect to the invariant-c$2$c forms}
  
  \begin{defi}\label{931}
  	Let $(\frak g\oplus \frak g^{*},[-,-]_{d},\lbrace-,-\rbrace_{d})$ be a compatible Lie algebra containing $(\frak g ,[-,-]_{\frak g},\lbrace-,-\rbrace_{\frak g})$ and $( \frak g^{*},[-,-]_{\frak g^{*}},\lbrace-,-\rbrace_{\frak g^{*}})$ as compatible Lie subalgebras. If the natural nondegenerate symmetric bilinear form $\mathcal{B}_{d}$   given by \eqref{367} is an invariant-c$2$c form on $(\frak g\oplus \frak g^{*},[-,-]_{d},\lbrace-,-\rbrace_{d})$, then we say $$\big( (\frak g\oplus\frak g^{*},[-,-]_{d},\lbrace-,-\rbrace_{d},\mathcal{B}_{d}),   \frak g,  \frak g^{*} \big)$$
  	is a {\bf Manin triple of compatible Lie algebras with respect to the invariant-c$2$c form}.
  \end{defi}
  
In Definition \ref{931},  it is clear that $	\big( (\frak g\oplus\frak g^{*},[-,-]_{d}, \mathcal{B}_{d}),   \frak g,  \frak g^{*} \big)$
 is a  Manin triple of  Lie algebras with respect to the nondegenerate symmetric invariant bilinear form \cite{CP1} and 
  $\big( (\frak g\oplus\frak g^{*}, \lbrace-,-\rbrace_{d},\mathcal{B}_{d}),   \frak g,  \frak g^{*} \big)$
  is a  Manin triple of Lie algebras with respect to the nondegenerate commutative $2$-cocycle \cite{TPA}.
  
  \begin{lem}\label{971}\cite{CP1,Maj}
  	Assume that $(\frak g\oplus\frak g^{*},[-,-]_{d})$ is a Lie algebra containing $(\frak g,[-,-]_{\frak g})$ and $(\frak g^{*},[-,-]_{\frak g^{*}})$ as Lie subalgebras.
  	Then the natural nondegenerate symmetric bilinear form $\mathcal{B}_{d}$ given by \eqref{367} is invariant on  $(\frak g\oplus\frak g^{*},[-,-]_{d})$  if and only if the multiplication $[-,-]_{d}$ on $\frak g\oplus\frak g^{*}$ is explicitly given by
  	\begin{equation}\label{974}
  		[x+a^{*},y+b^{*}]_{d}=[x,y]_{\frak g}+\mathrm{ad}^{*}_{\frak g^{*}}(a^{*})y-\mathrm{ad}^{*}_{\frak g^{*}}(b^{*})x
  		+[a^{*},b^{*}]_{\frak g^{*}}+\mathrm{ad}^{*}_{\frak g}(x)b^{*}-\mathrm{ad}^{*}_{\frak g}(y)a^{*}, 
  	\end{equation}
  for all $x,y\in\frak g,a^{*},b^{*}\in\frak g^{*}$.
  \end{lem}
  
  \begin{lem}\label{981}\cite{TPA}
  Suppose that $(\frak g,\circ_{\frak g})$ and $(\frak g^{*},\circ_{\frak g^{*}})$  are anti-pre-Lie algebras, whose sub-adjacent Lie algebras are $(\frak g,\lbrace-,-\rbrace_{\frak g})$ and $(\frak g^{*},\lbrace-,-\rbrace_{\frak g^{*}})$ respectively.
  Assume that $(\frak g\oplus\frak g^{*},\lbrace-,-\rbrace_{d})$ is a Lie algebra containing $(\frak g,\lbrace-,-\rbrace_{\frak g})$ and $(\frak g^{*},\lbrace-,-\rbrace_{\frak g^{*}})$ as Lie subalgebras.
  Then the natural nondegenerate symmetric bilinear form $\mathcal{B}_{d}$ given by \eqref{367} is a commutative $2$-cocycle on   $(\frak g\oplus\frak g^{*},\lbrace-,-\rbrace_{d})$ such that the anti-pre-Lie algebra $(\frak g\oplus\frak g^{*},\circ_{d})$ induced from  $\mathcal{B}_{d}$ through Proposition \ref{322}
  contains $(\frak g,\circ_{\frak g})$ and $(\frak g^{*},\circ_{\frak g^{*}})$ as anti-pre-Lie subalgebras
   if and only if the multiplication $\lbrace-,-\rbrace_{d}$ on $\frak g\oplus\frak g^{*}$ is explicitly given by
  \begin{eqnarray}\label{985}
  	\lbrace x+a^{*},y+b^{*}\rbrace_{d}
  	=\lbrace x,y\rbrace_{\frak g}
  	-L^{*}_{\circ_{\frak g^{*}}}(a^{*})y
  	+L^{*}_{\circ_{\frak g^{*}}}(b^{*})x
  	+\lbrace a^{*},b^{*}\rbrace_{\frak g^{*}}
  	-L^{*}_{\circ_{\frak g}}(x)b^{*}
  	+L^{*}_{\circ_{\frak g}}(y)a^{*}, 
  \end{eqnarray}
  for all $x,y\in\frak g,a^{*},b^{*}\in\frak g^{*}$.
  In this case the multiplication $\circ_{d}$ in the anti-pre-Lie algebra $(\frak g\oplus\frak g^{*},\circ_{d})$ is explicitly given by 
   \begin{eqnarray}\label{993}
  	 (x+a^{*})\circ_{d}(y+b^{*}) =x\circ_{\frak g}y-\mathfrak{ad}^{*}_{\frak g^{*}}(a^{*})y+R^{*}_{\circ_{\frak g^{*}}}(b^{*})x
  	+a^{*}\circ_{\frak g^{*}}b^{*}-\mathfrak{ad}^{*}_{\frak g}(x)b^{*}+R^{*}_{\circ_{\frak g}}(y)a^{*},
  \end{eqnarray}
for all $x,y\in\frak g,a^{*},b^{*}\in\frak g^{*}$.
  \end{lem}
  
  \begin{pro}\label{1005}
  	 Suppose that $(\frak g,[-,-]_{\frak g},\circ_{\frak g})$ and $(\frak g^{*},[-,-]_{\frak g^{*}},\circ_{\frak g^{*}})$ are compatible LAPL algebras, whose sub-adjacent compatible Lie algebras are $(\frak g,[-,-]_{\frak g},\lbrace-,-\rbrace_{\frak g})$ and $(\frak g^{*},[-,-]_{\frak g^{*}},\lbrace-,-\rbrace_{\frak g^{*}})$ respectively. Then the following conditions are equivalent:
  	 \begin{enumerate}
  	 	\item\label{1008} There is a Manin triple of compatible Lie algebras with respect to the invariant-c$2$c form $\big( (\frak g\oplus\frak g^{*}  ,[-,-]_{d},\lbrace-,-\rbrace_{ d},\mathcal{B}_{d}),\frak g, \frak g^{*}\big)$  such that the compatible LAPL algebra $(\frak g\oplus\frak g^{*},[-,-]_{d},\circ_{d})$ induced from $\mathcal{B}_{d}$ through Proposition \ref{673} contains $(\frak g,[-,-]_{\frak g},\circ_{\frak g})$ and $(\frak g^{*},[-,-]_{\frak g^{*}},\circ_{\frak g^{*}})$ as compatible LAPL subalgebras.
  	 	\item\label{1009} There is a compatible Lie algebra $(\frak g\oplus\frak g^{*},[-,-]_{d},\lbrace-,-\rbrace_{ d})$ in which 
  	 	$[-,-]_{d}$ is given by \eqref{974} and $\lbrace-,-\rbrace_{d}$ is given by \eqref{985}.
  	 	\item\label{1011} There is a compatible LAPL algebra $(\frak g\oplus\frak g^{*},[-,-]_{d},\circ_{ d})$ in which 
  	 	$[-,-]_{d}$ is given by \eqref{974} and $\circ_{d}$ is given by \eqref{993}.
  	 \end{enumerate}
  \end{pro}
  \begin{proof}
  	Suppose that Condition \eqref{1008} holds. Then by Lemmas \ref{971} and \ref{981}, the multiplications $[-,-]_{d}$ and $\circ_{d}$ in the compatible LAPL algebra $(\frak g\oplus\frak g^{*},[-,-]_{d},\circ_{d})$ induced from $\mathcal{B}_{d}$ through Proposition \ref{673} are given by \eqref{974} and \eqref{993} respectively.
  	Therefore we have \eqref{1008}$\Rightarrow$\eqref{1011}.
  	
  	Suppose that Condition \eqref{1011} holds. Then it is straightforward to check that the sub-adjacent compatible Lie algebra $(\frak g\oplus\frak g^{*},[-,-]_{d},\lbrace-,-\rbrace_{ d})$ of $(\frak g\oplus\frak g^{*},[-,-]_{d},\circ_{ d})$ satisfies \eqref{974} and \eqref{985}.
  	Therefore we have \eqref{1011}$\Rightarrow$\eqref{1009}.
  	
  		Suppose that Condition \eqref{1009} holds. Then by Lemmas \ref{971} and \ref{981} again, $\big( (\frak g\oplus\frak g^{*}  ,[-,-]_{d},\lbrace-,-\rbrace_{ d},\mathcal{B}_{d} ),\frak g$,
  		$\frak g^{*}\big)$ is a Manin triple of compatible Lie algebras with respect to the invariant-c$2$c form  such that the compatible LAPL algebra $(\frak g\oplus\frak g^{*},[-,-]_{d}$,
  		$\circ_{d})$ induced from $\mathcal{B}_{d}$ through Proposition \ref{673} contains $(\frak g,[-,-]_{\frak g},\circ_{\frak g})$ and $(\frak g^{*},[-,-]_{\frak g^{*}},\circ_{\frak g^{*}})$ as compatible LAPL subalgebras.
  		Therefore we have \eqref{1009}$\Rightarrow$\eqref{1008}.
   \end{proof}

  \subsection{Compatible LAPL bialgebras} Recall the notions of Lie bialgebras and anti-pre-Lie bialgebras.
  
   \begin{defi}\cite{CP1}
  	A {\bf Lie coalgebra} $(\frak g,\delta )$ is a vector space $\frak g$ together with a comultiplication $\delta :\frak g\rightarrow \frak g\otimes \frak g$ satisfying the following equations:
  	\begin{align*}
  	 \delta(x)+\tau\delta(x)&=0,\\
  		 S(\mathrm{id}_{\frak g}\otimes\delta)\delta(x)&=0,\;\forall x\in \frak g.
  	\end{align*}
  	A {\bf Lie bialgebra} is a triple $(\frak g,[-,-]_{\frak g},\delta)$ consisting of a Lie algebra $(\frak g,[-,-]_{\frak g})$ and a Lie coalgebra $(\frak g,\delta)$ such that $\delta$ is a $1$-cocycle of $(\frak g,[-,-]_{\frak g})$ with values in $\frak g\otimes \frak g$, that is,
  	\begin{eqnarray*}
  		\delta([x,y]_{\frak g})=\big( \mathrm{ad}_{\frak g}(x)\otimes\mathrm{id}_{\frak g}+\mathrm{id}_{\frak g}\otimes\mathrm{ad}_{\frak g}(x) \big)\delta(y)-\big( \mathrm{ad}_{\frak g}(y)\otimes\mathrm{id}_{\frak g}+\mathrm{id}_{\frak g}\otimes\mathrm{ad}_{\frak g}(y) \big)\delta(x), \;\forall x,y\in\frak g.
  	\end{eqnarray*}
  \end{defi}
  
  \begin{pro}\label{1043}\cite{CP1,Maj}
  	Let $(\frak g,[-,-]_{\frak g})$ and $(\frak g^{*},[-,-]_{\frak g^{*}})$ be Lie algebras. Let $\delta:\frak g\rightarrow\frak g\otimes\frak g$ be the linear dual of $[-,-]_{\frak g^{*}}$, that is,
  	 \begin{equation}\label{1154}
  		\langle \delta(x),a^{*}\otimes b^{*}\rangle=\langle x,[a^{*},b^{*}]_{\frak g^{*}}\rangle,\;\forall x\in\frak g, a^{*},b^{*}\in\frak g^{*}.
  	\end{equation}
  Then there is a Lie algebra $(\frak g\oplus\frak g^{*}, [-,-]_{d})$ with $[-,-]_{d}$ given by \eqref{974} if and only if $(\frak g,[-,-]_{\frak g},\delta)$ is a Lie bialgebra.
  \end{pro}
  
  \begin{defi}\cite{TPA}
  	An {\bf anti-pre-Lie coalgebra} $(\frak g,\theta )$ is a vector space $\frak g$ together with a comultiplication $\theta :\frak g\rightarrow \frak g\otimes \frak g$ satisfying the following equations:
  	\begin{align*}
  		 (\mathrm{id}_{\frak g}^{\otimes 3}-\tau\otimes\mathrm{id}_{\frak g})\big(  (\mathrm{id}_{\frak g}\otimes\theta)\theta+(\theta\otimes\mathrm{id}_{\frak g})\theta \big)(x)&=0,\\
  		 S( \tau\otimes\mathrm{id}_{\frak g}-  \mathrm{id}_{\frak g}^{\otimes 3})(\theta\otimes\mathrm{id}_{\frak g})\theta(x)&=0,\;\forall x\in\frak g.
  	\end{align*}
 An {\bf anti-pre-Lie bialgebra} is a triple $(\frak g,\circ_{\frak g},\theta)$ in which $(\frak g,\circ_{\frak g})$ is an anti-pre-Lie algebra, $(\frak g,\theta)$ is an anti-pre-Lie coalgebra and the following equations hold:
\begin{eqnarray*}
	&&\Delta(x\circ_{\frak g} y)=
	(\mathrm{id}_{\frak g}^{\otimes 2}-\tau)\Big(   
	\big( L_{\circ_{\frak g}}(x)\otimes\mathrm{id}_{\frak g}\big)\Delta(y)
	-\big(  \mathrm{id}_{\frak g}\otimes R_{\circ_{\frak g}}(y) \big)\theta(x)
	\Big),
	\\
	&&\theta( \lbrace x,y\rbrace_{\frak g})=\big(\mathrm{id}_{\frak g}\otimes \mathfrak{ad}_{\frak g}(x)-{L}_{\circ_{\frak g}}(x)\otimes \mathrm{id}_{\frak g}\big)\theta(y)
	-\big(\mathrm{id}_{\frak g}\otimes \mathfrak{ad}_{\frak g} (y)-{L}_{\circ_{\frak g}}(y)\otimes \mathrm{id}_{\frak g}\big)\theta(x),\;\forall x,y\in\frak g,\;\;
\end{eqnarray*}
in which $\lbrace -,-\rbrace_{\frak g}$ is given by \eqref{297} and $\Delta:\frak g\rightarrow\frak g\otimes\frak g$ is given by
\begin{equation}\label{1071}
\Delta(x)=\theta(x)-\tau\theta(x),\;\forall x\in\frak g.
\end{equation}
  \end{defi}
  
  \begin{pro} \cite{TPA}\label{1112}
  	 Suppose that $(\frak g,\circ_{\frak g})$ and $(\frak g^{*},\circ_{\frak g^{*}})$  are anti-pre-Lie algebras, whose sub-adjacent Lie algebras are $(\frak g,\lbrace-,-\rbrace_{\frak g})$ and $(\frak g^{*},\lbrace-,-\rbrace_{\frak g^{*}})$ respectively.
  	 Let $\theta:\frak g\rightarrow\frak g\otimes\frak g$ be the linear dual of $\circ_{\frak g^{*}}$, that is,
  	 \begin{equation}\label{1186}
  	 	\langle \theta(x),a^{*}\otimes b^{*}\rangle=\langle x, a^{*}\circ_{\frak g^{*}}b^{*} \rangle,\;\forall x\in\frak g, a^{*},b^{*}\in\frak g^{*}.
  	 \end{equation}
   Then there is a Lie algebra $(\frak g\oplus\frak g^{*}, \lbrace-,-\rbrace_{d})$ with $\lbrace-,-\rbrace_{d}$ given by \eqref{985}
   if and only if $(\frak g,\circ_{\frak g},\theta)$ is an anti-pre-Lie bialgebra.
  \end{pro}
  
  Recall \cite{WB} that a {\bf compatible Lie coalgebra} is a triple $(\frak g,\delta,\Delta)$ such that both $(\frak g,\delta)$ and $(\frak g,\Delta)$ are Lie coalgebras and the following equation holds:
  \begin{equation}\label{1115}
  	S\big( (\delta\otimes\mathrm{id}_{\frak g})\Delta+ (\Delta\otimes\mathrm{id}_{\frak g})\delta \big)(x)=0,\;\forall x \in\frak g.
  \end{equation}
  Now we introduce the notion of a compatible LAPL coalgebra.

  \begin{defi}
  	A {\bf compatible LAPL coalgebra} is a triple $(\frak g,\delta,\theta)$ such that the following conditions hold:
  	\begin{enumerate}
  		\item $(\frak g,\theta)$ is an anti-pre-Lie coalgebra.
  		\item $(\frak g,\delta,\Delta)$ with $\Delta$ given by \eqref{1071} is a compatible Lie coalgebra.
  		\item the following equation holds:
  		\begin{equation}\label{1098}
  			(\delta\otimes\mathrm{id}_{\frak g})\theta(x)+(\tau\otimes\mathrm{id}_{\frak g}-\mathrm{id}^{\otimes 3}_{\frak g})\big( (\mathrm{id}_{\frak g}\otimes\delta)\theta+(\theta\otimes\mathrm{id}_{\frak g}+\mathrm{id}_{\frak g}\otimes\theta)\delta  \big)(x)=0,\;\forall x\in\frak g.
  		\end{equation}
  	\end{enumerate}
\end{defi}

\begin{pro}
Let $\frak g$ be a vector space and $\delta,\theta:\frak g\rightarrow\frak g\otimes\frak g$ be comultiplications.
Then $(\frak g,\delta,\theta)$ is a compatible LAPL coalgebra if and only if $(\frak g^{*},[-,-]_{\frak g^{*}},\circ_{\frak g^{*}})$ 
with $[-,-]_{\frak g^{*}}$ given by \eqref{1154} and $\circ_{\frak g^{*}}$ given by \eqref{1186} is a compatible LAPL  algebra.
\end{pro}
  \begin{proof}
  	The conclusion follows from \cite{TPA,WB} and observing that
  	\eqref{1098} holds if and only if \eqref{646} holds on $\frak g^{*}$
  	under the natural nondegenerate tensor pairing.
  \end{proof}
  
  \begin{defi}
Suppose that $(\frak g,[-,-]_{\frak g},\circ_{\frak g})$ is a compatible LAPL algebra and $(\frak g,\delta,\theta)$ is a compatible LAPL coalgebra, such that $(\frak g,[-,-]_{\frak g},\delta)$ is a Lie bialgebra and $(\frak g,\circ_{\frak g},\theta)$ is an anti-pre-Lie bialgebra. 
Let $\lbrace-,-\rbrace_{\frak g}$ be given by \eqref{297} and $\Delta$ be given by \eqref{1071}.
If the following equations hold:
\begin{align}
 \theta([x,y]_{\frak g})
+\big(\mathrm{ad}_{\frak g}(y)\otimes\mathrm{id}_{\frak g}+\mathrm{id}_{\frak g}\otimes \mathrm{ad}_{\frak g}(y)\big)\theta(x)
-\big(\mathrm{ad}_{\frak g}(x)\otimes\mathrm{id}_{\frak g}+\mathrm{id}_{\frak g}\otimes \mathrm{ad}_{\frak g}(x)\big)\theta(y)&\nonumber\\
 -\delta(\lbrace x,y\rbrace_{\frak g})
+\big( \mathrm{id}_{\frak g}\otimes\mathfrak{ad}_{\frak g}(x)-L_{\circ_{\frak g}}(x)\otimes\mathrm{id}_{\frak g}  \big)\delta(y)
-\big( \mathrm{id}_{\frak g}\otimes\mathfrak{ad}_{\frak g}(y)-L_{\circ_{\frak g}}(y)\otimes\mathrm{id}_{\frak g}  \big)\delta(x)& =0,\label{1179}\\
 \delta(x\circ_{\frak g}y)
-\big( L_{\circ_{\frak g}}(x)\otimes\mathrm{id}_{\frak g}+\mathrm{id}_{\frak g}\otimes L_{\circ_{\frak g}}(x)  \big)\delta(y)
-\big( R_{\circ_{\frak g}}(y)\otimes\mathrm{id}_{\frak g}+\mathrm{id}_{\frak g}\otimes R_{\circ_{\frak g}}(y)  \big)\delta(x)&\nonumber\\
 -\Delta([x,y]_{\frak g})+\big( \mathrm{id}_{\frak g}\otimes\mathrm{ad}_{\frak g}(x)+ \mathrm{ad}_{\frak g}(x)\otimes \mathrm{id}_{\frak g} \big)\Delta(y)+(\mathrm{id}^{\otimes 2}_{\frak g}-\tau)\big(\mathrm{id}_{\frak g}\otimes\mathrm{ad}_{\frak g}(y)\big)\theta(x)&=0,\label{1183}
\end{align}
for all $x,y\in\frak g$, then we say $(\frak g,[-,-]_{\frak g},\circ_{\frak g},\delta,\theta)$ is a {\bf compatible LAPL bialgebra}.
  \end{defi}

  \begin{lem}\label{1188} \cite{WB}
  	Let $(\frak g,[-,-]_{\frak g},\lbrace-,-\rbrace_{\frak g}),(\frak h,[-,-]_{\frak h},\lbrace-,-\rbrace_{\frak h})$
  	be compatible Lie algebras and
  	$\rho,\mu:\frak g\rightarrow\mathrm{End}_{\mathbb{K}}(\frak h)$,
  	$s,t:\frak h\rightarrow\mathrm{End}_{\mathbb{K}}(\frak g)$
  	be linear maps.
  	Suppose that $(\rho,\mu,\frak h)$ is a representation of
  	$(\frak g,[-,-]_{\frak g},\lbrace-,-\rbrace_{\frak g})$
  	and $(s,t,\frak g)$ is a representation of
  	$(\frak h,[-,-]_{\frak h},\lbrace-,-\rbrace_{\frak h})$.
  	Define multiplications $[-,-]_{d},\lbrace-,-\rbrace_{d}$
  	on $\frak g\oplus\frak h$ by
  	\begin{align*}
  		[x+a,y+b]_{d}
  		&=[x,y]_{\frak g}+s(a)y-s(b)x
  		+[a,b]_{\frak h}+\rho(x)b-\rho(y)a,\\
  		\lbrace x+a,y+b\rbrace_{d}
  		&=\lbrace x,y\rbrace_{\frak g}+t(a)y-t(b)x
  		+\lbrace a,b\rbrace_{\frak h}+\mu(x)b-\mu(y)a,
  		\;\forall x,y\in\frak g,\ a,b\in\frak h.
  	\end{align*}
  	Then $(\frak g\oplus\frak h,[-,-]_{d},\lbrace-,-\rbrace_{d})$
  	is a compatible Lie algebra if and only if
  	$(\frak g\oplus\frak h,[-,-]_{d})$ and
  	$(\frak g\oplus\frak h,\lbrace-,-\rbrace_{d})$
  	are Lie algebras and the following equations hold for all $x,y\in\frak g,\ a,b\in\frak h$:
  	\begin{align}
  		 t(a)([x,y]_{\mathfrak g})
  		-[x,t(a)y]_{\mathfrak g}
  		+[y,t(a)x]_{\mathfrak g}
  		+t(\rho(x)a)y-t(\rho(y)a)x&\nonumber\\
  		 \quad+s(a)(\lbrace x,y\rbrace_{\mathfrak g})
  		-\lbrace x,s(a)y\rbrace_{\mathfrak g}
  		+\lbrace y,s(a)x\rbrace_{\mathfrak g}
  		+s(\mu(x)a)y-s(\mu(y)a)x&
  		=0,\label{1210}\\
  		 \mu(x)([a,b]_{\mathfrak h})
  		-[a,\mu(x)b]_{\mathfrak h}
  		+[b,\mu(x)a]_{\mathfrak h}
  		+\mu(s(a)x)b-\mu(s(b)x)a&\nonumber\\
  	 \quad+\rho(x)(\lbrace a,b\rbrace_{\mathfrak h})
  		-\lbrace a,\rho(x)b\rbrace_{\mathfrak h}
  		+\lbrace b,\rho(x)a\rbrace_{\mathfrak h}
  		+\rho(t(a)x)b-\rho(t(b)x)a&
  		=0.\label{1222}
  	\end{align}
  \end{lem}

  \begin{pro}\label{1226}
  	 Suppose that $(\frak g,[-,-]_{\frak g},\circ_{\frak g})$ and $(\frak g^{*},[-,-]_{\frak g^{*}},\circ_{\frak g^{*}})$ are compatible LAPL algebras, whose sub-adjacent compatible Lie algebras are $(\frak g,[-,-]_{\frak g},\lbrace-,-\rbrace_{\frak g})$ and $(\frak g^{*},[-,-]_{\frak g^{*}},\lbrace-,-\rbrace_{\frak g^{*}})$ respectively. Then there is a compatible Lie algebra $(\frak g\oplus\frak g^{*},[-,-]_{d},\lbrace-,-\rbrace_{ d})$ in which 
  		$[-,-]_{d}$ is given by \eqref{974} and $\lbrace-,-\rbrace_{d}$ is given by \eqref{985} if and only if $(\frak g,[-,-]_{\frak g},\circ_{\frak g},\delta,\theta)$ is a compatible LAPL bialgebra, in which $\delta$ is the linear dual of $[-,-]_{\frak g^{*}}$ and $\theta$ is the linear dual of $\circ_{\frak g^{*}}$.
  \end{pro}
  \begin{proof}
  	By Proposition \ref{1043}, $(\frak g\oplus\frak g^{*},[-,-]_{d})$ with $[-,-]_{d}$ given by \eqref{974} is a Lie algebra if and only if $(\frak g,[-,-]_{\frak g} ,\delta )$ is a Lie bialgebra.
  	By Proposition \ref{1112}, $(\frak g\oplus\frak g^{*}, \lbrace-,-\rbrace_{ d})$ with $\lbrace-,-\rbrace_{d}$ given by \eqref{985} is a Lie algebra if and only if $(\frak g, \circ_{\frak g} ,\theta)$ is an anti-pre-Lie bialgebra.
  	Moreover, for all $x,y\in\frak g, a^{*},b^{*}\in\frak g^{*}$, we have
  	\begin{align*}
-\langle L^{*}_{\circ_{\frak g^{*}}}(a^{*})([x,y]_{\frak g}), b^{*}\rangle
&= \langle [x,y]_{\frak g}, a^{*}\circ_{\frak g^{*}}b^{*}\rangle=\langle \theta([x,y]_{\frak g}), a^{*}\otimes b^{*}\rangle,\\
\langle [x,L^{*}_{\circ_{\frak g^{*}}}(a^{*})y]_{\frak g}, b^{*}\rangle&= -\langle L^{*}_{\circ_{\frak g^{*}}}(a^{*})y, \mathrm{ad}^{*}_{\frak g}(x)b^{*}\rangle=\langle y, a^{*}\circ_{\frak g^{*}} \mathrm{ad}^{*}_{\frak g}(x)b^{*}\rangle\\
&= -\langle \big( \mathrm{id}_{\frak g}\otimes\mathrm{ad}_{\frak g}(x)  \big)\theta(y), a^{*}\otimes b^{*}\rangle,\\
-\langle [y,L^{*}_{\circ_{\frak g^{*}}}(a^{*})x]_{\frak g}, b^{*}\rangle&= 
\langle \big( \mathrm{id}_{\frak g}\otimes\mathrm{ad}_{\frak g}(y)  \big)\theta(x), a^{*}\otimes b^{*}\rangle,\\
-\langle L^{*}_{\circ_{\frak g^{*}}}\big( \mathrm{ad}^{*}_{\frak g} (x)a^{*} \big)y, b^{*}\rangle&= \langle y, \mathrm{ad}^{*}_{\frak g}(x)a^{*}\circ_{\frak g^{*}}b^{*}\rangle=-\langle \big( \mathrm{ad}_{\frak g}(x)\otimes\mathrm{id}_{\frak g}  \big)\theta(y), a^{*}\otimes b^{*}\rangle,\\
\langle L^{*}_{\circ_{\frak g^{*}}}\big( \mathrm{ad}^{*}_{\frak g} (y)a^{*} \big)x, b^{*}\rangle&=  \langle \big( \mathrm{ad}_{\frak g}(y)\otimes\mathrm{id}_{\frak g}  \big)\theta(x), a^{*}\otimes b^{*}\rangle,\\
\langle \mathrm{ad}^{*}_{\frak g^{*}}(a^{*})(\lbrace x,y\rbrace_{\frak g}), b^{*}\rangle&= -\langle \lbrace x,y\rbrace_{\frak g}, [a^{*},b^{*}]_{\frak g^{*}} \rangle=-\langle \delta ( \lbrace x,y\rbrace_{\frak g} ), a^{*}\otimes b^{*}\rangle,\\
-\langle \lbrace x,\mathrm{ad}^{*}_{\frak g^{*}}(a^{*})y\rbrace_{\frak g}, b^{*}\rangle&= \langle \mathrm{ad}^{*}_{\frak g^{*}}(a^{*})y, \mathfrak{ad}^{*}_{\frak g}(x)b^{*}\rangle=-\langle y, [a^{*}, \mathfrak{ad}^{*}_{\frak g}(x)b^{*}]_{\frak g^{*}}\rangle\\
&= \langle \big( \mathrm{id}_{\frak g}\otimes\mathfrak{ad}_{\frak g}(x)  \big) \delta(y), a^{*}\otimes b^{*}   \rangle,\\
\langle \lbrace y,\mathrm{ad}^{*}_{\frak g^{*}}(a^{*})x\rbrace_{\frak g}, b^{*}\rangle&=  -\langle \big( \mathrm{id}_{\frak g}\otimes\mathfrak{ad}_{\frak g}(y)  \big) \delta(x), a^{*}\otimes b^{*}   \rangle,\\
-\langle \mathrm{ad}^{*}_{\frak g^{*}}\big( L^{*}_{\circ_{\frak g}}(x)a^{*}  \big)y, b^{*}\rangle&= \langle y, [ L^{*}_{\circ_{\frak g}}(x)a^{*}, b^{*} ]_{\frak g^{*}}\rangle=-\langle \big( L_{\circ_{\frak g}}(x)\otimes\mathrm{id}_{\frak g}  \big)\delta(y), a^{*}\otimes b ^{*}\rangle,\\
\langle \mathrm{ad}^{*}_{\frak g^{*}}\big( L^{*}_{\circ_{\frak g}}(y)a^{*}  \big)x, b^{*}\rangle&=  \langle \big( L_{\circ_{\frak g}}(y)\otimes\mathrm{id}_{\frak g}  \big)\delta(x), a^{*}\otimes b ^{*}\rangle.
  	\end{align*}
Therefore \eqref{1210} holds for the assumption 
\begin{equation}\label{1252}
	\frak h=\frak g^{*},\;
	\rho =\mathrm{ad}^{*}_{\frak g},\;
	\mu=-L^{*}_{\circ_{\frak g}},\;
	s=\mathrm{ad}^{*}_{\frak g^{*}},\;
	t=-L^{*}_{\circ_{\frak g^{*}}}
\end{equation}
if and only if \eqref{1179} holds.
Similarly, \eqref{1222} holds for the assumption \eqref{1252} if and only if \eqref{1183} holds. 
Since
$(\frak g,[-,-]_{\frak g},\circ_{\frak g})$ and
$(\frak g^{*},[-,-]_{\frak g^{*}},\circ_{\frak g^{*}})$
are compatible LAPL algebras, $(\mathrm{ad}^{*}_{\frak g}, - {L}^{*}_{\circ_{\frak g}},\frak g^{*})$ and $(\mathrm{ad}^{*}_{\frak g^{*}}, - {L}^{*}_{\circ_{\frak g^{*}}},\frak g )$
are representations of
$(\frak g,[-,-]_{\frak g},\lbrace-,-\rbrace_{\frak g})$
and
$(\frak g^{*},[-,-]_{\frak g^{*}},\lbrace-,-\rbrace_{\frak g^{*}})$ 
respectively.
Therefore by Lemma \ref{1188}, $(\frak g\oplus\frak g^{*},[-,-]_{d},\lbrace-,-\rbrace_{ d})$  is a compatible Lie algebra if and only if $(\frak g,[-,-]_{\frak g},\circ_{\frak g},\delta,\theta)$ is a compatible LAPL bialgebra.
  \end{proof}
  
  Combining Propositions \ref{1005} and \ref{1226}, we have the following result.
  
  \begin{cor}\label{1265}
  	Suppose that $(\frak g,[-,-]_{\frak g},\circ_{\frak g})$ and $(\frak g^{*},[-,-]_{\frak g^{*}},\circ_{\frak g^{*}})$ are compatible LAPL algebras, whose sub-adjacent compatible Lie algebras are $(\frak g,[-,-]_{\frak g},\lbrace-,-\rbrace_{\frak g})$ and $(\frak g^{*},[-,-]_{\frak g^{*}},\lbrace-,-\rbrace_{\frak g^{*}})$ respectively. Then the following conditions are equivalent:
  	\begin{enumerate}
  		\item\label{1268} There is a  Manin triple of compatible Lie algebras with respect to the invariant-c$2$c form $\big( (\frak g\oplus\frak g^{*},[-,-]_{d},\lbrace-,-\rbrace_{d},\mathcal{B}_{d}),\frak g,\frak g^{*}\big)$  such that the compatible LAPL algebra $(\frak g\oplus\frak g^{*},[-,-]_{d},\circ_{d})$ induced from $\mathcal{B}_{d}$ through Proposition \ref{673} contains $(\frak g,[-,-]_{\frak g},\circ_{\frak g})$ and $(\frak g^{*},[-,-]_{\frak g^{*}},\circ_{\frak g^{*}})$ as compatible LAPL subalgebras.
  		\item\label{1270} There is a compatible Lie algebra $(\frak g\oplus\frak g^{*},[-,-]_{d},\lbrace-,-\rbrace_{ d})$ in which 
  		$[-,-]_{d}$ is given by \eqref{974} and $\lbrace-,-\rbrace_{d}$ is given by \eqref{985}.
  		\item\label{1272} There is a compatible LAPL algebra $(\frak g\oplus\frak g^{*},[-,-]_{d},\circ_{ d})$ in which 
  		$[-,-]_{d}$ is given by \eqref{974} and $\circ_{d}$ is given by \eqref{993}.
  		\item $(\frak g,[-,-]_{\frak g},\circ_{\frak g},\delta,\theta)$ is a compatible LAPL bialgebra, in which $\delta$ is the linear dual of $[-,-]_{\frak g^{*}}$ and $\theta$ is the linear dual of $\circ_{\frak g^{*}}$.
  	\end{enumerate}
  \end{cor}

  \section{Triangular compatible LAPL bialgebras}
  
  In this section, we introduce the notion of the compatible LAPL Yang-Baxter equation in a compatible LAPL algebra, and show that an antisymmetric solution of this equation gives rise to a compatible LAPL bialgebra that is called triangular.
  Then we interpret solutions of the compatible LAPL Yang-Baxter equation in terms of $\mathcal{O}$-operators of compatible LAPL algebras.
  Finally, we introduce pre-compatible LAPL algebras as a means of constructing such $\mathcal{O}$-operators.
  
  \subsection{The compatible LAPL Yang-Baxter equation}\
  
  Let $(\frak g,[-,-]_{\frak g})$ be a Lie algebra and $r=\sum\limits_{i} a_{i}\otimes b_{i}\in\frak g\otimes\frak g$. Define $\textbf{C}(r)\in \frak g\otimes\frak g\otimes\frak g$ by 
   \begin{equation*}
  	\textbf{C}(r)=[r_{12},r_{13}]_{\frak g}+[r_{12},r_{23}]_{\frak g}+[r_{13},r_{23}]_{\frak g},
  \end{equation*}
 where
  \begin{equation*} 
  	[r_{12},r_{13}]= \sum_{i,j}[a_{i},a_{j}]_{\frak g}\otimes b_{i}\otimes b_{j},\; [r_{12} ,r_{23}]= \sum_{i,j}a_{i}\otimes [b_{i}, a_{j}]_{\frak g}\otimes b_{j},\; [r_{13},r_{23}]=\sum_{i,j}a_{i}\otimes a_{j}\otimes[b_{i},b_{j}]_{\frak g}.
  \end{equation*}
  We say $r$ is a solution of the \textbf{classical
  	Yang-Baxter equation} (CYBE in short) in $(\frak g,[-,-]_{\frak g})$ if the following equation holds:
  \begin{equation}\label{1292}
  	\textbf{C}(r)=0.
  \end{equation}
  
  \begin{pro}\cite{CP1}\label{1329}
  	 Let $(\frak g,[-,-]_{\frak g})$ be a Lie algebra, $r=\sum\limits_{i} a_{i}\otimes b_{i}\in\frak g\otimes\frak g$ and $\delta_{r}:\frak g\rightarrow
  	\frak g\otimes\frak g$ be a co-multiplication defined by 
  	\begin{equation}\label{1299}
  		 \delta_{r}(x) =\big(\mathrm{ad}_{\frak g}(x)\otimes \mathrm{id}_{\frak g}+\mathrm{id}_{\frak g}\otimes\mathrm{ad}_{\frak g}(x)\big)r, \;\;\forall x\in \frak g.
  	\end{equation}
  Then  $(\frak g,[-,-]_{\frak g},\delta_{r})$ is a Lie bialgebra if and only if 
  	the following equations hold:
  	\begin{align}
  		\big(\mathrm{ad}_{\frak g}(x)\otimes \mathrm{id}_{\frak g}+\mathrm{id}_{\frak g}\otimes\mathrm{ad}_{\frak g}(x)\big)\big(r+\tau(r)\big)&=0,&\label{eq:CYBE1}\\
  		\big (\mathrm{ad}_{\frak g}(x)\otimes \mathrm{id}_{\frak g}\otimes \mathrm{id}_{\frak g}+\mathrm{id}_{\frak g}\otimes\mathrm{ad}_{\frak g}(x)\otimes \mathrm{id}_{\frak g}+\mathrm{id}_{\frak g}\otimes \mathrm{id}_{\frak g}\otimes\mathrm{ad}_{\frak g}(x)\big)\textbf{C}(r)&=0,\;\forall x\in\frak g.\label{eq:CYBE2}
  	\end{align}
  In particular, if $r$ is an antisymmetric solution of the CYBE in $(\frak g,[-,-]_{\frak g})$, then $(\frak g,[-,-]_{\frak g},\delta_{r})$ is a Lie bialgebra, which is called a {\bf triangular Lie bialgebra}.
  \end{pro}
  
Let $(\frak g,\circ_{\frak g})$ be an anti-pre-Lie algebra and $r=\sum\limits_i a_i\otimes b_i\in\frak g\otimes\frak g$. Define $\textbf{T}(r)\in\frak g\otimes\frak g\otimes\frak g$ by
\begin{equation}\label{1318}
	\textbf{T}(r)=r_{12}\circ_{\frak g}r_{13}
	+r_{12}\circ_{\frak g}r_{23}
	-\lbrace r_{13},r_{23}\rbrace_{\frak g},
\end{equation}
where
\begin{equation*}
		r_{12}\circ_{\frak g}r_{13}
	 =\sum_{i,j}a_i\circ_{\frak g}a_j\otimes b_i\otimes b_j,\;
		r_{12}\circ_{\frak g}r_{23}
	 =\sum_{i,j}a_i\otimes b_i\circ_{\frak g}a_j\otimes b_j,\;
		\lbrace r_{13},r_{23}\rbrace_{\frak g}
	 =\sum_{i,j}a_i\otimes a_j\otimes
		\lbrace b_i,b_j\rbrace_{\frak g}.
\end{equation*}
We say that $r$ is a solution of the \textbf{anti-pre-Lie Yang--Baxter equation} (APL-YBE in short) in $(\frak g,\circ_{\frak g})$ if the following equation holds:
\begin{equation}\label{eq:APL-YBE}
	\textbf{T}(r)=0.
\end{equation}

\begin{pro}\cite{TPA}\label{1367}
	Let $(\frak g,\circ_{\frak g})$ be an anti-pre-Lie algebra, $r=\sum\limits_i a_i\otimes b_i\in\frak g\otimes\frak g$, and $\theta_r:\frak g\rightarrow\frak g\otimes\frak g$ be a co-multiplication defined by
	\begin{equation}\label{eq:aplCob}
	\theta_r(x) =
	\big(L_{\circ_{\frak g}}(x)\otimes\mathrm{id}_{\frak g}
	-\mathrm{id}_{\frak g}\otimes\mathfrak{ad}_{\frak g}(x)\big)r,
	\; \forall x\in\frak g.
\end{equation}
	Then $(\frak g,\circ_{\frak g},\theta_{r})$ is an anti-pre-Lie bialgebra if and only if the following equations hold:
{\small
	\begin{align}
		\Big( {L}_{\circ_{\frak g}}(x)\otimes\mathrm{id}_{\frak g}\otimes\mathrm{id}_{\frak g}-(\tau\otimes\mathrm{id}_{\frak g})\big( {L}_{\circ_{\frak g} }(x)\otimes\mathrm{id}_{\frak g}\otimes\mathrm{id}_{\frak g}\big)-\mathrm{id}_{\frak g}\otimes\mathrm{id}_{\frak g}\otimes\mathfrak{ad}_{\frak g}(x)\Big)\textbf{T}(r)&\nonumber\\
		+\sum_{j}\big(\mathrm{id}_{\frak g}\otimes {L}_{\circ_{\frak g}}(a_{j})\otimes\mathfrak{ad}_{\frak g}(x)-\mathfrak{ad}_{\frak g}(a_{j})\otimes\mathrm{id}_{\frak g}\otimes\mathfrak{ad}_{\frak g}(x)\big)\Big(\big(r+\tau(r)\big)\otimes b_{j}\Big)&\nonumber\\
		+(\mathrm{id}_{\frak g}^{\otimes 3}-\tau\otimes\mathrm{id}_{\frak g})\sum_{j}\big( {L}_{\circ_{\frak g}}(x\circ_{\frak g} a_{j})\otimes\mathrm{id}_{\frak g}\otimes\mathrm{id}_{\frak g}- {L}_{\circ_{\frak g}}(x) {R}_{\circ_{\frak g}}(a_{j})\otimes\mathrm{id}_{\frak g}\otimes\mathrm{id}_{\frak g}\big)\Big(\big(r+\tau(r)\big)\otimes b_{j}\Big)&=0,\label{eq:Apl 1}\\
	S\bigg(
	-\big(\mathrm{id}_{\frak g}\otimes\mathrm{id}_{\frak g}
	\otimes L_{\circ_{\frak g}}(x)\big)\textbf{T}(r)
	+\sum_{j}
	\big(\mathrm{id}_{\frak g}\otimes\mathrm{id}_{\frak g}
	\otimes L_{\circ_{\frak g}}(\lbrace x,b_j\rbrace_{\frak g})\big)
	(\mathrm{id}_{\frak g}^{\otimes3}-\tau\otimes\mathrm{id}_{\frak g})
	\Big(a_j\otimes\big(r+\tau(r)\big)\Big)
	&\nonumber\\
	+\bigg(
	\big(L_{\circ_{\frak g}}(a_j)\otimes\mathrm{id}_{\frak g}
	\otimes L_{\circ_{\frak g}}(x)\big)
	(\tau\otimes\mathrm{id}_{\frak g})
	-\big(\mathrm{id}_{\frak g}\otimes R_{\circ_{\frak g}}(a_j)
	\otimes L_{\circ_{\frak g}}(x)\big)
	\bigg)\Big(b_j\otimes\big(r+\tau(r)\big)\Big)
	&\nonumber\\
	+\big(\mathrm{id}_{\frak g}\otimes L_{\circ_{\frak g}}(a_j)
	\otimes L_{\circ_{\frak g}}(x)\big)
	\Big(\big(r+\tau(r)\big)\otimes b_j\Big)
	&\nonumber\\
	+\big(\mathrm{id}_{\frak g}\otimes\mathfrak{ad}_{\frak g}(b_j)
	\otimes L_{\circ_{\frak g}}(x)\big)
	\Big(a_j\otimes\big(r+\tau(r)\big)\Big)
	+\big(\mathfrak{ad}_{\frak g}(b_j)\otimes\mathrm{id}_{\frak g}
	\otimes L_{\circ_{\frak g}}(x)\big)
	\Big(\big(r+\tau(r)\big)\otimes a_j\Big)
	\bigg)
	&=0,\label{eq:Apl 2}\\
		\big(\mathrm{id}_{\frak g}\otimes {L}_{\circ_{\frak g}}(x\circ_{\frak g} y)-\mathrm{id}_{\frak g}\otimes {L}_{\circ_{\frak g}}(x) {L}_{\circ_{\frak g}}(y)&\nonumber\\
		+ {L}_{\circ_{\frak g}}(x) {L}_{\circ_{\frak g}}(y)\otimes \mathrm{id}_{\frak g}- {L}_{\circ_{\frak g}}(x\circ_{\frak g} y)\otimes \mathrm{id}_{\frak g}+ {L}_{\circ_{\frak g}}(y)\otimes {L}_{\circ_{\frak g}}(x)- {L}_{\circ_{\frak g}}(x)\otimes {L}_{\circ_{\frak g}}(y)\big)\big(r+\tau(r)\big)&=0,\label{1399}
	\end{align}
}
	for all $x,y\in\frak g$.
	In particular, if $r$ is an antisymmetric solution of the APL-YBE in $(\frak g,\circ_{\frak g})$, then $(\frak g,\circ_{\frak g},\theta_r)$ is an anti-pre-Lie bialgebra, which is called a \textbf{triangular anti-pre-Lie bialgebra}.
\end{pro}

Now we are motivated to introduce the notion of   the  compatible LAPL Yang-Baxter equation in a compatible LAPL algebra.

\begin{defi}
	Let $(\frak g,[-,-]_{\frak g},\circ_{\frak g})$ be a compatible LAPL algebra and $r\in\frak g\otimes\frak g$.
	If \eqref{1292} and \eqref{eq:APL-YBE} hold, then we say $r$ is a solution of the {\bf compatible LAPL Yang-Baxter equation} (CLAPL-YBE  in short) in $(\frak g,[-,-]_{\frak g},\circ_{\frak g})$.
\end{defi}

\begin{pro}\label{1413}
	Let $(\frak g,[-,-]_{\frak g},\circ_{\frak g})$ be a compatible LAPL algebra and $r=\sum\limits_{i}a_{i}\otimes b_{i}\in\frak g\otimes\frak g$. Let $\delta_{r}:\frak g\rightarrow
	\frak g\otimes\frak g$ be a co-multiplication defined by 
	\eqref{1299} and $\theta_r:\frak g\rightarrow\frak g\otimes\frak g$ be a co-multiplication defined by \eqref{eq:aplCob}.
	Then \eqref{1115} holds if and only if the following equation holds for all $x\in\frak g$:
	\begin{align}
	S\Bigg[ (\mathrm{id}_{\frak g}\otimes\tau-\mathrm{id}^{\otimes 3}_{\frak g})\big(\mathrm{ad}_{\frak g}(x)\otimes\mathrm{id}_{\frak g}\otimes\mathrm{id}_{\frak g}\big){\bf T}(r)-(R_{\circ_{\frak g}}(x)\otimes\mathrm{id}_{\frak g}\otimes\mathrm{id}_{\frak g}) {\bf C}(r)&\nonumber\\
	 +\sum_{j}\Bigg(\big( \mathfrak{ad}_{\frak g}(x)\otimes\mathrm{ad}_{\frak g}(a_{j})\otimes\mathrm{id}_{\frak g} \big)\big((r+\tau(r)\big)\otimes b_{j}\big)
	+\big(\mathrm{id}_{\frak g}\otimes\mathrm{ad}_{\frak g}(b_{j})\otimes\mathfrak{ad}_{\frak g}(x)\big) \big( a_{j}\otimes \big(r+\tau(r)\big) \big)&\nonumber\\
 -(\mathrm{id}_{\frak g}\otimes \mathrm{id}_{\frak g}\otimes \mathrm{ad}_{\frak g}(a_{j})\mathfrak{ad}_{\frak g}(x) ) \big( b_{j}\otimes (r+\tau(r)) \big) -\big( \mathrm{ad}_{\frak g}(b_{j})\mathfrak{ad}_{\frak g}(x)\otimes\mathrm{id}_{\frak g}\otimes\mathrm{id}_{\frak g} \big)\big( \big(r+\tau(r)\big) \otimes a_{j}\big)&\nonumber\\
	 +\big( \mathrm{id}_{\frak g}\otimes\mathrm{id}_{\frak g}\otimes\mathrm{ad}_{\frak g}(a_{j}\circ_{\frak g}x)  \big)\big( b_{j}\otimes \big(r+\tau(r)\big) \big)
	-\big( \mathrm{id}_{\frak g}\otimes\mathrm{id}_{\frak g}\otimes\mathfrak{ad}_{\frak g}([x,a_{j}]_{\frak g})\big)\big( b_{j}\otimes (r+\tau(r)) \big)&\nonumber\\
	 +\big( \mathfrak{ad}_{\frak g}([x,a_{j}]_{\frak g})\otimes\mathrm{id}_{\frak g} \otimes\mathrm{id}_{\frak g}\big)\big(  (r+\tau(r))\otimes b_{j} \big)&\nonumber\\
 +(\mathrm{id}_{\frak g}^{\otimes 3}-\mathrm{id}_{\frak g}\otimes \tau)\big( \mathrm{ad}_{\frak g}(x)\otimes R_{\circ_{\frak g}}(a_{j})\otimes\mathrm{id}_{\frak g} \big)\big(  (r+\tau(r))\otimes b_{j} \big)&\nonumber\\
	 \quad+(\mathrm{id}_{\frak g}^{\otimes 3}-\mathrm{id}_{\frak g}\otimes \tau)\big(\mathrm{id}_{\frak g}\otimes\mathfrak{ad}_{\frak g}(b_{j})\otimes\mathrm{ad}_{\frak g}(x)\big)\big( a_{j}\otimes (r+\tau(r)) \big)\Bigg)\Bigg]&=0.\label{1421}
	\end{align}
Moreover, \eqref{1098} holds if and only if the following equation holds  for all $x\in\frak g$: 
	\begin{align}
		\Big(\mathrm{id}_{\frak g}\otimes\mathrm{id}_{\frak g}
		\otimes\mathfrak{ad}_{\frak g}(x)-
		L_{\circ_{\frak g}}(x)\otimes\mathrm{id}_{\frak g}
		\otimes\mathrm{id}_{\frak g}
		-\mathrm{id}_{\frak g}\otimes L_{\circ_{\frak g}}(x)
		\otimes\mathrm{id}_{\frak g} 
		\Big)\textbf{C}(r)&\nonumber\\
		-\Big(
		\mathrm{ad}_{\frak g}(x)\otimes\mathrm{id}_{\frak g}
		\otimes\mathrm{id}_{\frak g}
		+\mathrm{id}_{\frak g}\otimes\mathrm{ad}_{\frak g}(x)
		\otimes\mathrm{id}_{\frak g} 
		+\mathrm{id}_{\frak g}\otimes\mathrm{id}_{\frak g}
		\otimes\mathrm{ad}_{\frak g}(x)
		\Big)\textbf{T}(r)&\nonumber\\
		+\sum_j\Bigg(
		\Big(
		\mathrm{ad}_{\frak g}(x\circ_{\frak g}a_j)
		-L_{\circ_{\frak g}}([x,a_j]_{\frak g})
		-L_{\circ_{\frak g}}(x)\mathrm{ad}_{\frak g}(a_j)
		+\mathrm{ad}_{\frak g}(x)R_{\circ_{\frak g}}(a_j)
		\Big)
		\otimes\mathrm{id}_{\frak g}\otimes\mathrm{id}_{\frak g}&\nonumber\\
		-\mathrm{ad}_{\frak g}(a_j)\otimes
		L_{\circ_{\frak g}}(x)\otimes\mathrm{id}_{\frak g}
		+\mathrm{id}_{\frak g}\otimes
		L_{\circ_{\frak g}}([x,a_j]_{\frak g})
		\otimes\mathrm{id}_{\frak g}
		+R_{\circ_{\frak g}}(a_j)\otimes
		\mathrm{ad}_{\frak g}(x)\otimes\mathrm{id}_{\frak g}
		&\nonumber\\
		-\mathfrak{ad}_{\frak g}(a_j)\otimes
		\mathrm{id}_{\frak g}\otimes\mathrm{ad}_{\frak g}(x)
		+\mathrm{id}_{\frak g}\otimes
		L_{\circ_{\frak g}}(a_j)\otimes\mathrm{ad}_{\frak g}(x)
		\Bigg)
		\Big((r+\tau(r))\otimes b_j\Big)&=0.\label{1706}
	\end{align}
\end{pro}
\begin{proof}
	For all $x\in\frak g$, we have
	\begin{equation*}
		\Delta_{r}(x)=\theta_{r}(x)-\tau\theta_{r}(x)=\sum_{i} x\circ_{\frak g} a_{i}\otimes b_{i}-a_{i}\otimes \lbrace x, b_{i}\rbrace_{\frak g}-b_{i}\otimes x\circ_{\frak g}a_{i}+\lbrace x, b_{i}\rbrace_{\frak g}\otimes a_{i}.
	\end{equation*}	
	Therefore, we have
	\begin{align*}
		(\delta_{r}\otimes\mathrm{id}_{\frak g})\Delta_{r}(x)&=	(\delta_{r}\otimes\mathrm{id}_{\frak g})(\sum_{i} x\circ_{\frak g} a_{i}\otimes b_{i}-a_{i}\otimes \lbrace x, b_{i}\rbrace_{\frak g}-b_{i}\otimes x\circ_{\frak g}a_{i}+\lbrace x, b_{i}\rbrace_{\frak g}\otimes a_{i})\\
		&=\sum_{i} \delta_{r}(x\circ_{\frak g} a_{i})\otimes b_{i}-\delta_{r}(a_{i})\otimes \lbrace x, b_{i}\rbrace_{\frak g}-\delta_{r}(b_{i})\otimes x\circ_{\frak g}a_{i}+\delta_{r}(\lbrace x, b_{i}\rbrace_{\frak g})\otimes a_{i}\\
		&=\sum_{i,j} [x\circ_{\frak g} a_{i}, a_{j}]_{\frak g}\otimes b_{j}\otimes b_{i}+a_{j}\otimes [x\circ_{\frak g}a_{i}, b_{j}]_{\frak g}\otimes b_{i}\\
		&\quad -[a_{i}, a_{j}]_{\frak g}\otimes b_{j}\otimes \lbrace x, b_{i}\rbrace_{\frak g}-a_{j}\otimes [a_{i}, b_{j}]_{\frak g}\otimes \lbrace x, b_{i}\rbrace_{\frak g}\\
		&\quad-[b_{i}, a_{j}]_{\frak g}\otimes b_{j}\otimes x\circ_{\frak g}a_{i}-a_{j}\otimes [b_{i}, b_{j}]_{\frak g}\otimes x\circ_{\frak g} a_{i}\\
		&\quad +[\lbrace x, b_{i}\rbrace_{\frak g}, a_{j}]_{\frak g}\otimes b_{j}\otimes a_{i}+a_{j}\otimes [\lbrace x, b_{i}\rbrace_{\frak g}, b_{j}]_{\frak g}\otimes a_{i},
	\end{align*}
	and similarly
	\begin{align*}
		(\Delta_{r}\otimes\mathrm{id}_{\frak g})\delta_{r}(x)&=\sum_{i,j} [x, a_{i}]_{\frak g}\circ_{\frak g} a_{j}\otimes b_{j}\otimes b_{i}-a_{j}\otimes \lbrace [x,a_{i}]_{\frak g}, b_{j}\rbrace_{\frak g}\otimes b_{i}\\
		&\quad -b_{j}\otimes [x, a_{i}]_{\frak g}\circ_{\frak g}a_{j}\otimes b_{i}+\lbrace [x,a_{i}]_{\frak g}, b_{j}\rbrace_{\frak g}\otimes a_{j}\otimes b_{i}\\
		&\quad +a_{i}\circ_{\frak g}a_{j}\otimes b_{j}\otimes [x, b_{i}]_{\frak g}- a_{j}\otimes \lbrace a_{i}, b_{j}\rbrace_{\frak g}\otimes [x, b_{i}]_{\frak g}\\
		&\quad -b_{j}\otimes a_{i}\circ_{\frak g} a_{j}\otimes [x, b_{i}]_{\frak g}+\lbrace a_{i}, b_{j}\rbrace_{\frak g}\otimes a_{j}\otimes [x, b_{i}]_{\frak g}.
	\end{align*}
	Thus we have
	\begin{align*}
	 S\big(
		(\delta_{r}\otimes\mathrm{id}_{\frak g})\Delta_{r}
		+(\Delta_{r}\otimes\mathrm{id}_{\frak g})\delta_{r}
		\big)(x) 
		&=S\sum_{i,j}
		\begin{aligned}[t]
			\Bigl(&[x\circ_{\frak g}a_j,a_i]_{\frak g}\otimes b_i\otimes b_j
			+[x\circ_{\frak g}a_i,b_j]_{\frak g}\otimes b_i\otimes a_j\\
			&-\lbrace x,b_i\rbrace_{\frak g}\otimes[a_i,a_j]_{\frak g}\otimes b_j
			-\lbrace x,b_i\rbrace_{\frak g}\otimes a_j\otimes[a_i,b_j]_{\frak g}\\
			&-x\circ_{\frak g}a_i\otimes[b_i,a_j]_{\frak g}\otimes b_j
			-x\circ_{\frak g}a_i\otimes a_j\otimes[b_i,b_j]_{\frak g}\\
			&+[\lbrace x,b_i\rbrace_{\frak g},a_j]_{\frak g}\otimes b_j\otimes a_i
			+[\lbrace x,b_i\rbrace_{\frak g},b_j]_{\frak g}\otimes a_i\otimes a_j\\
			&+[x,a_j]_{\frak g}\circ_{\frak g}a_i\otimes b_i\otimes b_j
			-\lbrace[x,a_i]_{\frak g},b_j\rbrace_{\frak g}\otimes b_i\otimes a_j\\
			&-[x,a_i]_{\frak g}\circ_{\frak g}a_j\otimes b_i\otimes b_j
			+\lbrace[x,a_i]_{\frak g},b_j\rbrace_{\frak g}\otimes a_j\otimes b_i\\
			&+[x,b_i]_{\frak g}\otimes a_i\circ_{\frak g}a_j\otimes b_j
			-[x,b_i]_{\frak g}\otimes a_j\otimes\lbrace a_i,b_j\rbrace_{\frak g}\\
			&-[x,b_i]_{\frak g}\otimes b_j\otimes a_i\circ_{\frak g}a_j
			+[x,b_i]_{\frak g}\otimes\lbrace a_i,b_j\rbrace_{\frak g}\otimes a_j
			\Bigr).
		\end{aligned}
	\end{align*}
Noticing that
\begin{align*}
	& S (\sum_{i,j} [x\circ_{\frak g} a_i,b_j]_{\frak g}
	\otimes b_i\otimes a_j)\\
	=& S \Bigl(
	\sum_j
	(\mathrm{id}_{\frak g}\otimes\mathrm{id}_{\frak g}
	\otimes\mathrm{ad}_{\frak g}(x\circ_{\frak g}a_j))
	\bigl(b_j\otimes(r+\tau(r))\bigr)
	-\sum_{i,j}[x\circ_{\frak g}a_i,a_j]_{\frak g}
	\otimes b_i\otimes b_j
	\Bigr),\\
	& S (\sum_{i,j}
	-\lbrace x,b_i\rbrace_{\frak g}\otimes[a_i,a_j]_{\frak g}\otimes b_j
	-\lbrace x,b_i\rbrace_{\frak g}\otimes a_j\otimes[a_i,b_j]_{\frak g})\\
	=& S \Bigl(
	\sum_j\Bigl(
	\bigl(\mathfrak{ad}_{\frak g}(x)\otimes\mathrm{ad}_{\frak g}(a_j)
	\otimes\mathrm{id}_{\frak g}\bigr)
	\bigl((r+\tau(r))\otimes b_j\bigr)
	+\bigl(\mathrm{id}_{\frak g}\otimes\mathrm{ad}_{\frak g}(b_j)
	\otimes\mathfrak{ad}_{\frak g}(x)\bigr)
	\bigl(a_j\otimes(r+\tau(r))\bigr)
	\Bigr)\\
	&\quad
	+\bigl(\mathfrak{ad}_{\frak g}(x)\otimes\mathrm{id}_{\frak g}
	\otimes\mathrm{id}_{\frak g}\bigr){\bf C}(r)
	-\sum_{i,j}\lbrace x,[a_i,a_j]_{\frak g}\rbrace_{\frak g}
	\otimes b_i\otimes b_j
	\Bigr),\\
	& S (\sum_{i,j}
	-x\circ_{\frak g}a_i\otimes[b_i,a_j]_{\frak g}\otimes b_j
	-x\circ_{\frak g}a_i\otimes a_j\otimes[b_i,b_j]_{\frak g})\\
	=& S \Bigl(
	-(L_{\circ_{\frak g}}(x)\otimes\mathrm{id}_{\frak g}
	\otimes\mathrm{id}_{\frak g}){\bf C}(r)
	+\sum_{i,j}x\circ_{\frak g}[a_i,a_j]_{\frak g}
	\otimes b_i\otimes b_j
	\Bigr),\\
	& S (\sum_{i,j}
	[\lbrace x,b_i\rbrace_{\frak g},a_j]_{\frak g}
	\otimes b_j\otimes a_i)\\
	=& S \Bigl(
	\sum_j
	-(\mathrm{id}_{\frak g}\otimes\mathrm{id}_{\frak g}
	\otimes\mathrm{ad}_{\frak g}(a_j)\mathfrak{ad}_{\frak g}(x))
	\bigl(b_j\otimes(r+\tau(r))\bigr)
	-\sum_{i,j}[\lbrace x,a_j\rbrace_{\frak g},a_i]_{\frak g}
	\otimes b_i\otimes b_j
	\Bigr),\\
	& S (\sum_{i,j}
	[\lbrace x,b_i\rbrace_{\frak g},b_j]_{\frak g}
	\otimes a_i\otimes a_j)\\
	=& S \Bigl(
	\sum_j\Bigl(
	-\bigl(\mathrm{ad}_{\frak g}(b_j)\mathfrak{ad}_{\frak g}(x)
	\otimes\mathrm{id}_{\frak g}\otimes\mathrm{id}_{\frak g}\bigr)
	\bigl((r+\tau(r))\otimes a_j\bigr)\\
	&\quad
	-\bigl(\mathrm{id}_{\frak g}\otimes\mathrm{id}_{\frak g}
	\otimes\mathrm{ad}_{\frak g}(\lbrace x,a_j\rbrace_{\frak g})\bigr)
	\bigl(b_j\otimes(r+\tau(r))\bigr)
	\Bigr)
	+\sum_{i,j}[\lbrace x,a_i\rbrace_{\frak g},a_j]_{\frak g}
	\otimes b_i\otimes b_j
	\Bigr),\\
	& S (\sum_{i,j}
	-\lbrace[x,a_i]_{\frak g},b_j\rbrace_{\frak g}
	\otimes b_i\otimes a_j)\\
	=& S \Bigl(
	\sum_j
	-\bigl(\mathrm{id}_{\frak g}\otimes\mathrm{id}_{\frak g}
	\otimes\mathfrak{ad}_{\frak g}([x,a_j]_{\frak g})\bigr)
	\bigl(b_j\otimes(r+\tau(r))\bigr)
	+\sum_{i,j}\lbrace[x,a_i]_{\frak g},a_j\rbrace_{\frak g}
	\otimes b_i\otimes b_j
	\Bigr).\\
	& S (\sum_{i,j}
	\lbrace[x,a_i]_{\frak g},b_j\rbrace_{\frak g}
	\otimes a_j\otimes b_i)\\
	=& S \Bigl(
	\sum_j
	\bigl(\mathfrak{ad}_{\frak g}([x,a_j]_{\frak g})
	\otimes\mathrm{id}_{\frak g}\otimes\mathrm{id}_{\frak g}\bigr)
	\bigl((r+\tau(r))\otimes b_j\bigr)
	-\sum_{i,j}\lbrace[x,a_j]_{\frak g},a_i\rbrace_{\frak g}
	\otimes b_i\otimes b_j
	\Bigr),\\
	& S (\sum_{i,j}
	[x,b_i]_{\frak g}\otimes a_i\circ_{\frak g}a_j\otimes b_j
	-[x,b_i]_{\frak g}\otimes a_j\otimes
	\lbrace a_i,b_j\rbrace_{\frak g})\\
	=& S \Bigl(
	\sum_j\Bigl(
	\bigl(\mathrm{ad}_{\frak g}(x)\otimes R_{\circ_{\frak g}}(a_j)
	\otimes\mathrm{id}_{\frak g}\bigr)
	\bigl((r+\tau(r))\otimes b_j\bigr)
	+\bigl(\mathrm{id}_{\frak g}\otimes\mathfrak{ad}_{\frak g}(b_j)
	\otimes\mathrm{ad}_{\frak g}(x)\bigr)
	\bigl(a_j\otimes(r+\tau(r))\bigr)
	\Bigr)\\
	&\quad
	-\bigl(\mathrm{ad}_{\frak g}(x)\otimes\mathrm{id}_{\frak g}
	\otimes\mathrm{id}_{\frak g}\bigr){\bf T}(r)
	+\sum_{i,j}[x,a_i\circ_{\frak g}a_j]_{\frak g}
	\otimes b_i\otimes b_j
	\Bigr),\\
	& S (\sum_{i,j}
	-[x,b_i]_{\frak g}\otimes b_j\otimes a_i\circ_{\frak g}a_j
	+[x,b_i]_{\frak g}\otimes
	\lbrace a_i,b_j\rbrace_{\frak g}\otimes a_j)\\
	=& S \Bigl(
	-\sum_j\Bigl(
	(\mathrm{id}_{\frak g}\otimes\tau)
	\bigl(\mathrm{ad}_{\frak g}(x)\otimes R_{\circ_{\frak g}}(a_j)
	\otimes\mathrm{id}_{\frak g}\bigr)
	\bigl((r+\tau(r))\otimes b_j\bigr)\\
	&\quad
	+(\mathrm{id}_{\frak g}\otimes\tau)
	\bigl(\mathrm{id}_{\frak g}\otimes\mathfrak{ad}_{\frak g}(b_j)
	\otimes\mathrm{ad}_{\frak g}(x)\bigr)
	\bigl(a_j\otimes(r+\tau(r))\bigr)
	\Bigr)
	+(\mathrm{id}_{\frak g}\otimes\tau)
	\bigl(\mathrm{ad}_{\frak g}(x)\otimes\mathrm{id}_{\frak g}
	\otimes\mathrm{id}_{\frak g}\bigr){\bf T}(r)\\
	&\quad
	-\sum_{i,j}[x,a_j\circ_{\frak g}a_i]_{\frak g}
	\otimes b_i\otimes b_j
	\Bigr),
\end{align*}
and, for all $i,j$,
	\begin{align*}
		& [x\circ_{\frak g}a_{j}, a_{i}]_{\frak g}-[x\circ_{\frak g}a_{i}, a_{j}]_{\frak g}-\lbrace x, [a_{i}, a_{j}]_{\frak g}\rbrace_{\frak g}+x\circ_{\frak g}[a_{i}, a_{j}]_{\frak g}-[\lbrace x, a_{j}\rbrace_{\frak g}, a_{i}]_{\frak g} 
		+[\lbrace x, a_{i}\rbrace_{\frak g},a_{j}]_{\frak g}\\
		&\quad +[x,a_{j}]_{\frak g}\circ_{\frak g}a_{i}+\lbrace [x,a_{i}]_{\frak g}, a_{j}\rbrace_{\frak g}-[x,a_{i}]_{\frak g}\circ_{\frak g} a_{j}-\lbrace [x,a_{j}]_{\frak g}, a_{i}\rbrace_{\frak g}+[x,a_{i}\circ_{\frak g}a_{j}]_{\frak g}-[x,a_{j}\circ_{\frak g}a_{i}]_{\frak g}\\
		 &\overset{\eqref{292}}{=}  [x\circ_{\frak g}a_{j},a_{i}]_{\frak g}-[x\circ_{\frak g}a_{i},a_{j}]_{\frak g}+x\circ_{\frak g}[a_{i},a_{j}]_{\frak g}+[x,a_{j}]_{\frak g}\circ_{\frak g}a_{i}-[x,a_{i}]_{\frak g}\circ_{\frak g}a_{j}+2[x,\lbrace a_{i},a_{j}\rbrace_{\frak g}]_{\frak g}\\
		 &\overset{\eqref{886}}{=} 0,
	\end{align*}
	we finally have $S\big(
	(\delta_r\otimes\mathrm{id}_{\frak g})\Delta_r
	+(\Delta_r\otimes\mathrm{id}_{\frak g})\delta_r
	\big)(x)=$ the left hand side of \eqref{1421}.
	Therefore \eqref{1115} holds if and only if  \eqref{1421} holds. 
	Similarly, we can prove that \eqref{1098} holds if and only if \eqref{1706} holds.
\end{proof}

   \begin{pro}\label{1720}
   	Let $(\frak g,[-,-]_{\frak g},\circ_{\frak g})$ be a compatible LAPL algebra and $r=\sum\limits_{i}a_{i}\otimes b_{i}\in\frak g\otimes\frak g$. Let $\delta_{r}:\frak g\rightarrow
   	\frak g\otimes\frak g$ be a co-multiplication defined by 
   \eqref{1299} and $\theta_r:\frak g\rightarrow\frak g\otimes\frak g$ be a co-multiplication defined by \eqref{eq:aplCob}.
   Then \eqref{1179} holds automatically. Moreover, \eqref{1183} holds if and only if the following equation holds for all $x,y\in\frak g$:
  \begin{align}
 \bigg( \big( \mathrm{ad}_{\frak g}(x)\mathfrak{ad}_{\frak g}(y) +\mathrm{ad}_{\frak g}(y)\mathfrak{ad}_{\frak g}(x)-\mathfrak{ad}_{\frak g}([x,y]_{\frak g}) \big)\otimes\mathrm{id}_{\frak g} +\mathrm{id}_{\frak g}\otimes \big( L_{\circ_{\frak g}}([x,y]_{\frak g})-\mathrm{ad}_{\frak g}(x)L_{\circ_{\frak g}}(y)\big)&\nonumber\\
 +\big(\mathfrak{ad}_{\frak g}(y)\otimes\mathrm{ad}_{\frak g}(x)-\mathrm{ad}_{\frak g}(x)\otimes L_{\circ_{\frak g}}(y) -\mathrm{ad}_{\frak g}(y)\otimes L_{\circ_{\frak g}}(x)  \big) \bigg)\big(r+\tau(r)\big)&=0. \label{1384}
  \end{align}
   \end{pro}
  \begin{proof}
For all $x,y\in\frak g$, we have
\begin{align*}
 &\delta_{r}(x\circ_{\frak g}y)
-\big( L_{\circ_{\frak g}}(x)\otimes\mathrm{id}_{\frak g}+\mathrm{id}_{\frak g}\otimes L_{\circ_{\frak g}}(x)  \big)\delta_{r}(y)
-\big( R_{\circ_{\frak g}}(y)\otimes\mathrm{id}_{\frak g}+\mathrm{id}_{\frak g}\otimes R_{\circ_{\frak g}}(y)  \big)\delta_{r}(x) \\
&-\Delta_{r}([x,y]_{\frak g})+\big( \mathrm{id}_{\frak g}\otimes\mathrm{ad}_{\frak g}(x)+ \mathrm{ad}_{\frak g}(x)\otimes \mathrm{id}_{\frak g} \big)\Delta_{r}(y)+(\mathrm{id}^{\otimes 2}_{\frak g}-\tau)\big(\mathrm{id}_{\frak g}\otimes\mathrm{ad}_{\frak g}(y)\big)\theta_{r}(x) \\
=&\sum_{i} [x\circ_{\frak g}y, a_{i}]_{\frak g}\otimes b_{i}+a_{i}\otimes [x\circ_{\frak g}y, b_{i}]_{\frak g}
- \big( L_{\circ_{\frak g}}(x)\otimes\mathrm{id}_{\frak g}+\mathrm{id}_{\frak g}\otimes L_{\circ_{\frak g}}(x)  \big)( [y,a_{i}]_{\frak g}\otimes b_{i}+a_{i}\otimes [y, b_{i}]_{\frak g} )\\
&- \big( R_{\circ_{\frak g}}(y)\otimes\mathrm{id}_{\frak g}+\mathrm{id}_{\frak g}\otimes R_{\circ_{\frak g}}(y)  \big)( [x,a_{i}]_{\frak g}\otimes b_{i}+a_{i}\otimes [x, b_{i}]_{\frak g} )\\
&-[x,y]_{\frak g}\circ_{\frak g} a_{i}\otimes b_{i}+a_{i}\otimes \lbrace [x,y]_{\frak g}, b_{i}\rbrace_{\frak g}+b_{i}\otimes [x,y]_{\frak g}\circ_{\frak g} a_{i}-\lbrace [x,y]_{\frak g}, b_{i}\rbrace_{\frak g}\otimes a_{i}\\
&+\big( \mathrm{id}_{\frak g}\otimes\mathrm{ad}_{\frak g}(x)+\mathrm{ad}_{\frak g}(x)\otimes  \mathrm{id}_{\frak g} \big)( y\circ_{\frak g}a_{i}\otimes b_{i}-a_{i}\otimes \lbrace y , b_{i}\rbrace_{\frak g}-b_{i}\otimes y\circ_{\frak g}a_{i}+\lbrace y,b_{i}\rbrace_{\frak g}\otimes a_{i} )\\
&+(\mathrm{id}^{\otimes 2}_{\frak g}-\tau)\big(\mathrm{id}_{\frak g}\otimes\mathrm{ad}_{\frak g}(y)\big) (x\circ_{\frak g} a_{i}\otimes b_{i}-a_{i}\otimes\lbrace x, b_{i}\rbrace_{\frak g}) \\
=&\sum_{i} [x\circ_{\frak g}y, a_{i}]_{\frak g}\otimes b_{i}+a_{i}\otimes [x\circ_{\frak g}y, b_{i}]_{\frak g}-x\circ_{\frak g}[y, a_{i}]_{\frak g}\otimes b_{i}-x\circ_{\frak g}a_{i}\otimes [y, b_{i}]_{\frak g}-[y, a_{i}]_{\frak g}\otimes x\circ_{\frak g}b_{i}\\
&-a_{i}\otimes x\circ_{\frak g}[y, b_{i}]_{\frak g}-[x, a_{i}]_{\frak g}\circ_{\frak g}y\otimes b_{i}-a_{i}\circ_{\frak g} y\otimes [x,b_{i}]_{\frak g}-[x,a_{i}]_{\frak g}\otimes b_{i}\circ_{\frak g}y-a_{i}\otimes [x,b_{i}]_{\frak g}\circ_{\frak g}y\\
&-[x,y]_{\frak g}\circ_{\frak g} a_{i}\otimes b_{i}+a_{i}\otimes \lbrace [x,y]_{\frak g}, b_{i}\rbrace_{\frak g}+b_{i}\otimes [x,y]_{\frak g}\circ_{\frak g}a_{i}-\lbrace [x,y]_{\frak g}, b_{i}\rbrace_{\frak g}\otimes a_{i}\\
&+y\circ_{\frak g}a_{i}\otimes [x,b_{i}]_{\frak g}-a_{i}\otimes [x,\lbrace y, b_{i}\rbrace_{\frak g}]_{\frak g}-b_{i}\otimes [x,y\circ_{\frak g}a_{i}]_{\frak g}+\lbrace y, b_{i}\rbrace_{\frak g}\otimes [x, a_{i}]_{\frak g}\\
&+[x, y\circ_{\frak g} a_{i}]_{\frak g}\otimes b_{i}-[x, a_{i}]_{\frak g}\otimes \lbrace y, b_{i}\rbrace_{\frak g}-[x, b_{i}]_{\frak g}\otimes y\circ_{\frak g} a_{i}+[x, \lbrace y, b_{i}\rbrace_{\frak g}]_{\frak g}\otimes a_{i}\\
&+x\circ_{\frak g} a_{i}\otimes [y, b_{i}]_{\frak g}-a_{i}\otimes [y, \lbrace x, b_{i}\rbrace_{\frak g}]_{\frak g}-[y, b_{i}]_{\frak g}\otimes x\circ_{\frak g}a_{i}+[y, \lbrace x, b_{i}\rbrace_{\frak g}]_{\frak g}\otimes a_{i}\\ 
=&A(1)+A(2)+A(3)+A(4)+A(5)+A(6),
\end{align*}
where 
\begin{align*}
A(1) =&\sum_{i}-[y, a_{i}]_{\frak g}\otimes x\circ_{\frak g}b_{i}-[y, b_{i}]_{\frak g}\otimes x\circ_{\frak g}a_{i}=-\big(\mathrm{ad}_{\frak g}(y )\otimes L_{\circ_{\frak g}}(x)\big)\big( r+\tau(r) \big),\\
A(2) =&\sum_{i}-[x,a_{i}]_{\frak g}\otimes b_{i}\circ_{\frak g}y-[x, a_{i}]_{\frak g}\otimes \lbrace y, b_{i}\rbrace_{\frak g}-[x, b_{i}]_{\frak g}\otimes y\circ_{\frak g}a_{i}=-\big(\mathrm{ad}_{\frak g}(x )\otimes L_{\circ_{\frak g}}(y)\big)\big( r+\tau(r) \big),\\
A(3) =&\sum_{i} -x\circ_{\frak g}a_{i}\otimes [y, b_{i}]_{\frak g}+x\circ_{\frak g}a_{i}\otimes [y, b_{i}]_{\frak g}=0,\\
A(4) =&\sum_{i} -a_{i}\circ_{\frak g}y\otimes [x, b_{i}]_{\frak g}+y\circ_{\frak g}a_{i}\otimes [x, b_{i}]_{\frak g}+\lbrace y, b_{i}\rbrace_{\frak g}\otimes [x, a_{i}]_{\frak g}=\big( \mathfrak{ad}_{\frak g}(y)\otimes\mathrm{ad}_{\frak g}(x)  \big)\big(r+\tau(r)\big),\\
A(5) =&\sum_{i} [x\circ_{\frak g} y, a_{i}]_{\frak g}\otimes b_{i}-x\circ_{\frak g}[y, a_{i}]_{\frak g}\otimes b_{i}-[x, a_{i}]_{\frak g}\circ_{\frak g}y\otimes b_{i}-[x,y]_{\frak g}\circ_{\frak g}a_{i}\otimes b_{i}+[x, y\circ_{\frak g}a_{i}]_{\frak g}\otimes b_{i}\\
&-\lbrace [x,y]_{\frak g}, b_{i}\rbrace_{\frak g}\otimes a_{i}+[x, \lbrace y, b_{i}\rbrace_{\frak g}]_{\frak g}\otimes a_{i}+[y,\lbrace x, b_{i}\rbrace_{\frak g}]_{\frak g}\otimes a_{i}\\
=&\bigg(\big(\mathrm{ad}_{\frak g}(x)\mathfrak{ad}_{\frak g}(y) +\mathrm{ad}_{\frak g}(y)\mathfrak{ad}_{\frak g}(x)-\mathfrak{ad}_{\frak g}([x,y]_{\frak g}) \big)\otimes\mathrm{id}_{\frak g}\bigg)\big(  r+\tau(r) \big)\\
&+\sum_{i} [x\circ_{\frak g} y, a_{i}]_{\frak g}\otimes b_{i}-x\circ_{\frak g}[y, a_{i}]_{\frak g}\otimes b_{i}-[x, a_{i}]_{\frak g}\circ_{\frak g}y\otimes b_{i}-[x,y]_{\frak g}\circ_{\frak g}a_{i}\otimes b_{i}+[x, y\circ_{\frak g}a_{i}]_{\frak g}\otimes b_{i}\\
&+\lbrace [x,y]_{\frak g},a_{i}\rbrace_{\frak g}\otimes b_{i}-[x, \lbrace y, a_{i}\rbrace_{\frak g}]_{\frak g}\otimes b_{i}-[y,\lbrace x, a_{i}\rbrace_{\frak g}]_{\frak g}\otimes b_{i}\\
\overset{\eqref{646}}{=}&\bigg(\big(\mathrm{ad}_{\frak g}(x)\mathfrak{ad}_{\frak g}(y) +\mathrm{ad}_{\frak g}(y)\mathfrak{ad}_{\frak g}(x)-\mathfrak{ad}_{\frak g}([x,y]_{\frak g}) \big)\otimes\mathrm{id}_{\frak g}\bigg)\big(  r+\tau(r) \big),\\
A(6)=&\sum_{i} a_{i}\otimes [x\circ_{\frak g}y, b_{i}]_{\frak g}-a_{i}\otimes x\circ_{\frak g}[y, b_{i}]_{\frak g}-a_{i}\otimes [x, b_{i}]_{\frak g}\circ_{\frak g}y+a_{i}\otimes \lbrace [x,y]_{\frak g},b_{i}\rbrace_{\frak g}\\
&-a_{i}\otimes [x, \lbrace y, b_{i}\rbrace_{\frak g}]_{\frak g}-a_{i}\otimes [y, \lbrace x, b_{i}\rbrace_{\frak g}]_{\frak g}+b_{i}\otimes [x,y]_{\frak g}\circ_{\frak g}a_{i}-b_{i}\otimes [x, y\circ_{\frak g}a_{i}]_{\frak g}\\
=&\bigg(\mathrm{id}_{\frak g}\otimes \big( L_{\circ_{\frak g}}([x,y]_{\frak g})-\mathrm{ad}_{\frak g}(x)L_{\circ_{\frak g}}(y)\big)\bigg)\big( r+\tau(r) \big)\\
&+\sum_{i} a_{i}\otimes [x\circ_{\frak g}y, b_{i}]_{\frak g}-a_{i}\otimes x\circ_{\frak g}[y, b_{i}]_{\frak g}-a_{i}\otimes [x, b_{i}]_{\frak g}\circ_{\frak g}y+a_{i}\otimes \lbrace [x,y]_{\frak g},b_{i}\rbrace_{\frak g}\\
&-a_{i}\otimes [x, \lbrace y, b_{i}\rbrace_{\frak g}]_{\frak g}-a_{i}\otimes [y, \lbrace x, b_{i}\rbrace_{\frak g}]_{\frak g}-a_{i}\otimes [x,y]_{\frak g}\circ_{\frak g}b_{i}+a_{i}\otimes [x, y\circ_{\frak g}b_{i}]_{\frak g}\\
\overset{\eqref{646}}{=}&\bigg(\mathrm{id}_{\frak g}\otimes \big( L_{\circ_{\frak g}}([x,y]_{\frak g})-\mathrm{ad}_{\frak g}(x)L_{\circ_{\frak g}}(y)\big)\bigg)\big( r+\tau(r) \big).
\end{align*}
Therefore \eqref{1183} holds if and only if \eqref{1384} holds.
Similarly, we can prove that \eqref{1179} holds automatically.
  \end{proof}
  
  \begin{cor}\label{cor:triangular-compatible-lapl}
  	Let $(\frak g,[-,-]_{\frak g},\circ_{\frak g})$ be a compatible LAPL algebra and $r=\sum\limits_i a_i\otimes b_i\in\frak g\otimes\frak g$.
  	Let $\delta_r:\frak g\rightarrow\frak g\otimes\frak g$ be a co-multiplication defined by \eqref{1299} and $\theta_r:\frak g\rightarrow\frak g\otimes\frak g$ be a co-multiplication defined by \eqref{eq:aplCob}.
  	Then $(\frak g,[-,-]_{\frak g},\circ_{\frak g},\delta_r,\theta_r)$ is a compatible LAPL bialgebra if and only if
  	\eqref{eq:CYBE1}, \eqref{eq:CYBE2}, \eqref{eq:Apl 1}--\eqref{1399}
  	and  \eqref{1421}--\eqref{1384} hold.
  	In particular, if $r$ is an antisymmetric solution of the CLAPL-YBE in $(\frak g,[-,-]_{\frak g},\circ_{\frak g})$, then $(\frak g,[-,-]_{\frak g},\circ_{\frak g},\delta_r,\theta_r)$ is a compatible LAPL bialgebra, which is called a {\bf triangular compatible LAPL bialgebra}.
  \end{cor}
  \begin{proof}
  	It follows from Propositions \ref{1329}, \ref{1367}, \ref{1413} and \ref{1720}.
  \end{proof}
  
  \begin{pro}\label{1818}
  	Let $(\frak g,[-,-]_{\frak g},\circ_{\frak g})$ be a compatible LAPL algebra and $r\in\frak g\otimes\frak g$ be antisymmetric.
  	Then $r$ is a solution of the CLAPL-YBE in $(\frak g,[-,-]_{\frak g},\circ_{\frak g})$ if and only if the following equations hold:
  	\begin{eqnarray}
&& [r(a^{*}), r(b^{*})]_{\frak g}=r\Big( \mathrm{ad}^{*}_{\frak g}\big( r(a^{*}) \big)b^{*} -\mathrm{ad}^{*}_{\frak g}\big( r(b^{*}) \big)a^{*} \Big),\label{eq:of1}\\
&& r(a^{*})\circ_{\frak g} r(b^{*})=r\Big( -\mathfrak{ad}^{*}_{\frak g}\big( r(a^{*}) \big)b^{*}+R^{*}_{\circ_{\frak g}}\big( r(b^{*}) \big)a^{*}  \Big),\;\forall a^{*},b^{*}\in\frak g^{*}.\label{eq:of2}
  	\end{eqnarray}
  \end{pro}
  \begin{proof}
  	By \cite{Kup}, $\textbf{C}(r)=0$ if and only if \eqref{eq:of1} holds.
  	By \cite{TPA}, $\textbf{T}(r)=0$ if and only if \eqref{eq:of2} holds.
  Therefore the conclusion holds.
  \end{proof}
  
  \subsection{$\mathcal{O}$-operators of compatible LAPL algebras}\
  
  \begin{defi}
  	Let $(\rho,l_{\circ},r_{\circ},V)$ be a representation
  	of  a compatible LAPL algebra  $(\frak g,[-,-]_{\frak g},\circ_{\frak g})$.  If there is a
  	linear map $T:V\rightarrow\frak g$ satisfying the following conditions:
  	\begin{enumerate}
  		\item\cite{Kup} $T$ is an $\mathcal{O}$-operator of $(\frak g,[-,-]_{\frak g})$ associated to $(\rho,V)$, that is,
  		\begin{equation*}
  		 [T(u),T(v)]_{\frak g}=T\Big( \rho\big(T(u)\big)v-\rho\big(T(v)\big)u  \Big),\;\forall u,v\in V,
  		\end{equation*}
  		\item\cite{TPA} $T$ is an $\mathcal{O}$-operator of $(\frak g,\circ_{\frak g})$ with respect to $(l_{\circ} ,r_{\circ} ,V)$, that is,
  		\begin{equation*}
  			T(u)\circ_{\frak g} T(v)=T\Big( l_{\circ} \big(T(u)\big)v+r_{\circ} \big(T(v)\big)u     \Big),\;\forall u,v\in V,
  		\end{equation*}
  	\end{enumerate}
  	then we say $T$ is an {\bf $\mathcal{O}$-operator of $(\frak g,[-,-]_{\frak g},\circ_{\frak g})$ with respect to $(\rho,l_{\circ},r_{\circ},V)$}.
  \end{defi}
  
  Now we can rewrite Proposition \ref{1818} as follows.
  
  \begin{cor}\label{1853}
  	Let $(\frak g,[-,-]_{\frak g},\circ_{\frak g})$ be a compatible LAPL algebra and $r\in \frak g\otimes \frak g$ be antisymmetric.
  	Then $r$ is a solution of the CLAPL-YBE in $(\frak g,[-,-]_{\frak g},\circ_{\frak g})$ if and only if $r:\frak g^{*}\rightarrow \frak g$ is an $\mathcal{O}$-operator of $(\frak g,[-,-]_{\frak g},\circ_{\frak g})$ associated to $( \mathrm{ad}^{*}_{\frak g},-\mathfrak{ad}^{*}_{\frak g}, {R}^{*}_{\circ_{\frak g}},\frak g^{*})$.
  \end{cor}

  \begin{thm}\label{1859}
  		Let $(\rho,l_{\circ},r_{\circ},V)$ be a representation
  	of  a compatible LAPL algebra  $(\frak g,[-,-]_{\frak g},\circ_{\frak g})$.
  	Let $T:V\rightarrow \frak g$ be a linear map which is identified as an element in 
  	$$V^{*}\otimes \frak g\subset \frak g\ltimes_{\rho^{*},r_{\circ}^{*} -l_{\circ}^{*} , r_{\circ}^{*  }}V^{*}\otimes \frak g\ltimes_{\rho^{*},r_{\circ}^{*} -l_{\circ}^{*} , r_{\circ}^{*} }V^{*}.$$  Then $r=T-\tau(T)$ is a solution of the CLAPL-YBE in the compatible LAPL algebra $( \frak g\ltimes_{\rho^{*},r_{\circ}^{*} -l_{\circ}^{*} , r_{\circ}^{*} }V^{*} )$ if and only if
  	$T$ is an $\mathcal{O}$-operator of $(\frak g,[-,-]_{\frak g},\circ_{\frak g})$ with respect to $(\rho,l_{\circ},r_{\circ},V)$.
  \end{thm}
  \begin{proof}
  	By \cite{Bai2007}, $r=T-\tau(T)$ is a solution of the CYBE in  $\frak g\ltimes_{\rho^{*}} V^{*} $ if and only if
  	$T$ is an $\mathcal{O}$-operator of $(\frak g,[-,-]_{\frak g})$ with respect to $(\rho, V)$.
  	By \cite{TPA}, $r=T-\tau(T)$ is a solution of the APL-YBE in $\frak g\ltimes_{r_{\circ}^{*}-l_{\circ}^{*}, r_{\circ}^{*}} V^{*} $ if and only if
  	$T$ is an $\mathcal{O}$-operator of $(\frak g,\circ_{\frak g} )$ associated to $(l_{\circ} ,r_{\circ} ,V)$.
  Therefore the conclusion holds.
  \end{proof}

 Let $(\frak g,\diamond_{\frak g})$ be a pre-Lie algebra. Then there is a Lie algebra $(\frak g,[-,-]_{\frak g})$ with the Lie bracket $[-,-]_{\frak g}$ given by
 \begin{equation*}
 	[x,y]_{\frak g}=x\diamond_{\frak g} y-y\diamond_{\frak g} x,\;\forall x,y\in\frak g,
 \end{equation*}
which is called the {\bf sub-adjacent Lie algebra} of $(\frak g,\diamond_{\frak g})$.
 Moreover, by \cite{Bai2007}, $(L_{\diamond_{\frak g}},\frak g)$ is a representation of $(\frak g,[-,-]_{\frak g})$ and $T=\mathrm{id}_{\frak g}$ is an $\mathcal{O}$-operator of $(\frak g,[-,-]_{\frak g})$ with respect to $(L_{\diamond_{\frak g}},\frak g)$.
 
 On the other hand, recall \cite{TPA} that a {\bf pre-APL algebra}  is a triple $(\frak g,\succ_{\frak g}$,
 $\prec_{\frak g})$, where $\frak g$ is a vector space and $\succ_{\frak g},\prec_{\frak g}:\frak g\otimes \frak g\rightarrow \frak g$ are multiplications such that the following equations hold:
 	\begin{align*}
 		x\succ_{\frak g}(y\succ_{\frak g} z)-y\succ_{\frak g}(x\succ_{\frak g} z)
 		&= (y\succ_{\frak g} x+y\prec_{\frak g} x)\succ_{\frak g} z-(x\succ_{\frak g} y+x\prec_{\frak g} y)\succ_{\frak g} z, \ \ \ \ \  \\
 		z\prec_{\frak g}(x\succ_{\frak g} y+ x\prec_{\frak g} y)&= x\succ_{\frak g}(z\prec_{\frak g} y)+(x\succ_{\frak g} z)\prec_{\frak g} y-(z\prec_{\frak g} x)\prec_{\frak g} y, \ \ \ \ \  \\
 		x\succ_{\frak g}(y\succ_{\frak g} z)-y\succ_{\frak g}(x\succ_{\frak g} z)
 		&= (y\succ_{\frak g} z)\prec_{\frak g} x-(x\succ_{\frak g} z)\prec_{\frak g} y-(z\prec_{\frak g} y)\prec_{\frak g} x+(z\prec_{\frak g} x)\prec_{\frak g} y, 
 	\end{align*}
for all $x,y,z\in\frak g$. Let $(\frak g,\succ_{\frak g},\prec_{\frak g})$ be a pre-APL algebra. Then there is an anti-pre-Lie algebra $(\frak g,\circ_{\frak g})$ in which $\circ_{\frak g}$ is given by
 \begin{equation*}
 	x\circ_{\frak g} y=x\succ_{\frak g} y+ x\prec_{\frak g} y,\; \forall x,y\in \frak g,
 \end{equation*}
 which is called the {\bf sub-adjacent anti-pre-Lie algebra
 	of $(\frak g,\succ_{\frak g},\prec_{\frak g})$}.
 Moreover, $( {L}_{\succ_{\frak g}}, {R}_{\prec_{\frak g}},\frak g)$ is a representation of $(\frak g,\circ_{\frak g})$, and $T=\mathrm{id}_{\frak g}$ is an $\mathcal{O}$-operator of $(\frak g,\circ_{\frak g})$ associated to $( {L}_{\succ_{\frak g}}, {R}_{\prec_{\frak g}}, \frak g)$.

  \begin{defi}
  	Let $(\frak g,\diamond_{\frak g})$ be a pre-Lie algebra and
  	$(\frak g,[-,-]_{\frak g})$ be the sub-adjacent Lie algebra
  	of $(\frak g,\diamond_{\frak g})$.
  	Let $(\frak g,\succ_{\frak g},\prec_{\frak g})$ be a
  	pre-APL algebra \cite{TPA}, $(\frak g,\circ_{\frak g})$
  	be the sub-adjacent anti-pre-Lie algebra of
  	$(\frak g,\succ_{\frak g},\prec_{\frak g})$ and $(\frak g,\lbrace -,-\rbrace_{\frak g})$ be the sub-adjacent Lie algebra of $(\frak g,\circ_{\frak g})$.
  	If the following equations hold for all $x,y,z\in\frak g$:
  	\begin{align}
  		[x,y]_{\frak g}\succ_{\frak g}z-z\prec_{\frak g} [x,y]_{\frak g}-x\succ_{\frak g}(y\diamond_{\frak g}z)+(y\diamond_{\frak g} z)\prec_{\frak g} x+y\succ_{\frak g}(x\diamond_{\frak g} z)&\nonumber\\
  		-(x\diamond_{\frak g} z)\prec_{\frak g}y+\lbrace x,y\rbrace_{\frak g}\diamond_{\frak g}z-x\diamond_{\frak g}(y\succ_{\frak g}z)+x\diamond_{\frak g}(z\prec_{\frak g}y)+y\diamond_{\frak g}(x\succ_{\frak g} z)-y\diamond_{\frak g}(z\prec_{\frak g}x)&=0,\label{1906}\\
  		[x,y]_{\frak g}\succ_{\frak g}z-x\succ_{\frak g}(y\diamond_{\frak g}z) 
  		+y\succ_{\frak g}(x\diamond_{\frak g} z)-\lbrace x,y\rbrace_{\frak g}\diamond_{\frak g}z
  		-x\diamond_{\frak g}(y\succ_{\frak g}z)
  		+y\diamond_{\frak g}(x\succ_{\frak g} z)&=0,\label{1909}\\
  		(x\diamond_{\frak g} z)\prec_{\frak g}y+x\succ_{\frak g}(y\diamond_{\frak g}z)+z\prec_{\frak g}[x,y]_{\frak g} &\nonumber\\
  		+y\diamond_{\frak g}(x\succ_{\frak g}z  - z\prec_{\frak g}x)-x\diamond_{\frak g}(z\prec_{\frak g}y)-(x\circ_{\frak g}y)\diamond_{\frak g}z&=0,\label{1912}
  	\end{align}
  	then we say $(\frak g,\diamond_{\frak g},
  	\succ_{\frak g},\prec_{\frak g})$ is a
  	{\bf pre-compatible LAPL algebra}.
  \end{defi}

  \begin{pro}\label{1916}
  	Let $(\frak g,\diamond_{\frak g},\succ_{\frak g},\prec_{\frak g})$ be a  pre-compatible LAPL algebra, $(\frak g,[-,-]_{\frak g})$ be the sub-adjacent Lie algebra of $(\frak g,\diamond_{\frak g})$ and $(\frak g,\circ_{\frak g})$ be the sub-adjacent anti-pre-Lie algebra of $(\frak g,\succ_{\frak g},\prec_{\frak g})$.
  	Then the triple $(\frak g,[-,-]_{\frak g},\circ_{\frak g})$ is a compatible LAPL algebra, which is called the {\bf sub-adjacent compatible LAPL algebra} of $(\frak g,\diamond_{\frak g},\succ_{\frak g},\prec_{\frak g})$. 
  	Moreover, $(L_{\diamond_{\frak g}},L_{\succ_{\frak g}},R_{\prec_{\frak g}},\frak g)$ is a representation of $(\frak g,[-,-]_{\frak g},\circ_{\frak g})$ and $T=\mathrm{id}_{\frak g}$ is an $\mathcal{O}$-operator of $(\frak g,[-,-]_{\frak g},\circ_{\frak g})$ with respect to $(L_{\diamond_{\frak g}},L_{\succ_{\frak g}},R_{\prec_{\frak g}},\frak g)$.
  \end{pro}
  \begin{proof}
  	It is clear that $(\frak g,[-,-]_{\frak g})$ is a Lie algebra and $(\frak g,\circ_{\frak g})$ is an anti-pre-Lie algebra. Let $(\frak g,\lbrace -,-\rbrace_{\frak g})$ be the sub-adjacent Lie algebra of $(\frak g,\circ_{\frak g})$.
  	Then we have
  	\begin{align*}
  		&\lbrace[x,y]_{\frak g},z\rbrace_{\frak g}+ \lbrace[y,z]_{\frak g},x\rbrace_{\frak g}+ \lbrace[z,x]_{\frak g},y\rbrace_{\frak g}
  		+[\lbrace x,y\rbrace_{\frak g},z]_{\frak g}+ [\lbrace y,z\rbrace_{\frak g},x]_{\frak g}+ [\lbrace z,x\rbrace_{\frak g},y]_{\frak g}\\
  		&=[x,y]_{\frak g}\succ_{\frak g}z+(x\diamond_{\frak g} y-y\diamond_{\frak g} x)\prec_{\frak g}z-z\succ_{\frak g}(x\diamond_{\frak g}y-y\diamond_{\frak g} x)-z\prec_{\frak g}[x,y]_{\frak g}\\
  		&\quad +[y,z]_{\frak g}\succ_{\frak g}x
  		+(y\diamond_{\frak g}z-z\diamond_{\frak g}y)\prec_{\frak g}x -x\succ_{\frak g}(y\diamond_{\frak g}z-z\diamond_{\frak g}y)
  		-x\prec_{\frak g}[y,z]_{\frak g}\\
  		&\quad +[z,x]_{\frak g}\succ_{\frak g}y
  		+(z\diamond_{\frak g}x-x\diamond_{\frak g}z)\prec_{\frak g}y -y\succ_{\frak g}(z\diamond_{\frak g}x-x\diamond_{\frak g}z)
  		-y\prec_{\frak g}[z,x]_{\frak g}\\
  		&\quad +\lbrace x,y\rbrace_{\frak g}\diamond_{\frak g}z-z\diamond_{\frak g}(x\succ_{\frak g}y+x\prec_{\frak g}y)+z\diamond_{\frak g}(y\succ_{\frak g}x+y\prec_{\frak g}x)\\
  		& \quad+\lbrace y,z\rbrace_{\frak g}\diamond_{\frak g}x -x\diamond_{\frak g}(y\succ_{\frak g}z+y\prec_{\frak g}z)
  		+x\diamond_{\frak g}(z\succ_{\frak g}y+z\prec_{\frak g}y)\\
  		&\quad +\lbrace z,x\rbrace_{\frak g}\diamond_{\frak g}y -y\diamond_{\frak g}(z\succ_{\frak g}x+z\prec_{\frak g}x)
  		+y\diamond_{\frak g}(x\succ_{\frak g}z+x\prec_{\frak g}z)\\
  		&\overset{\eqref{1906}}{=}0,
  	\end{align*}
  	that is, $(\frak g,[-,-]_{\frak g},\lbrace-,-\rbrace_{\frak g})$ is a compatible Lie algebra. Similarly, by \eqref{1909} and \eqref{1912}, we obtain \eqref{646}.
  	Thus $(\frak g,[-,-]_{\frak g},\circ_{\frak g})$ is a compatible LAPL algebra.
  	Moreover, it is straightforward to check that the quadruple $(L_{\diamond_{\frak g}},L_{\succ_{\frak g}},R_{\prec_{\frak g}},\frak g)$ is a representation of $(\frak g,[-,-]_{\frak g},\circ_{\frak g})$ and $T=\mathrm{id}_{\frak g}$ is an $\mathcal{O}$-operator of $(\frak g,[-,-]_{\frak g},\circ_{\frak g})$ with respect to $(L_{\diamond_{\frak g}},L_{\succ_{\frak g}},R_{\prec_{\frak g}},\frak g)$.
  \end{proof}

  \begin{pro}\label{1946}
  	Let $(\frak g,\diamond_{\frak g},\succ_{\frak g},\prec_{\frak g})$ be a pre-compatible LAPL algebra and $(\frak g,[-,-]_{\frak g},\circ_{\frak g})$ be the sub-adjacent compatible LAPL algebra of $(\frak g,\diamond_{\frak g},\succ_{\frak g},\prec_{\frak g})$.
  	Let $\{ e_{1},\cdots,e_{n} \}$ be a basis of $\frak g$ and $\{ e^{*}_{1},\cdots,e^{*}_{n} \}$ be the dual basis.
  	Then
  	\begin{equation*}
  		r=\sum_{i=1}^{n}( e^{*}_{i}\otimes e_{i}-e_{i}\otimes e^{*}_{i} )
  	\end{equation*}
  	is an antisymmetric solution of the CLAPL-YBE in the compatible LAPL algebra
  	$\frak g\ltimes_{  {L}_{\diamond_{\frak g}}^{*}, {R}^{*}_{\prec_{\frak g}}- {L}^{*}_{\succ_{\frak g}}, {R}^{*}_{\prec_{\frak g}}}\frak g^{*}$.
  \end{pro}
  \begin{proof}
  	By Proposition \ref{1916}, $T=\mathrm{id}_{\frak g}$ is an $\mathcal{O}$-operator of $(\frak g,[-,-]_{\frak g},\circ_{\frak g})$ with respect to $(L_{\diamond_{\frak g}},L_{\succ_{\frak g}},R_{\prec_{\frak g}},\frak g)$. Note that $T=\mathrm{id}_{\frak g}$ is identified as
  	$\sum\limits_{i=1}^{n}e^{*}_{i}\otimes e_{i}$ in $\frak g^{*}\otimes
  	\frak g\subset
  	\frak g\ltimes_{  {L}_{\diamond_{\frak g}}^{*}, {R}^{*}_{\prec_{\frak g}}- {L}^{*}_{\succ_{\frak g}}, {R}^{*}_{\prec_{\frak g}}}\frak g^{*}\otimes
  	\frak g\ltimes_{  {L}_{\diamond_{\frak g}}^{*}, {R}^{*}_{\prec_{\frak g}}- {L}^{*}_{\succ_{\frak g}}, {R}^{*}_{\prec_{\frak g}}}\frak g^{*}$.
  	Hence the conclusion follows from Theorem \ref{1859}.
  \end{proof}

\section{From differential Novikov bialgebras to compatible LAPL bialgebras}

In this section, we extend to the bialgebraic level the construction of a compatible LAPL algebra from a differential Novikov algebra.
We show that a Manin triple of differential Novikov algebras induces a Manin triple of compatible Lie algebras with respect to an invariant-c$2$c form.
  Equivalently from the bialgebraic viewpoint,  a differential Novikov bialgebra gives rise to a compatible LAPL bialgebra, and this property is still available when both bialgebras are limited to the triangular cases.
  
  Recall \cite{HBG} that a {\bf Manin triple of Novikov algebras} is a collection $\big( (\frak g\oplus\frak g^{*},\star_{d},\mathcal{B}_{d}),\frak g,\frak g^{*}\big)$ such that $(\frak g\oplus\frak g^{*},\star_{d},\mathcal{B}_{d})$ is a quadratic Novikov algebra and $(\frak g,\star_{\frak g}),(\frak g^{*},\star_{\frak g^{*}})$ are Novikov subalgebras of $(\frak g\oplus\frak g^{*},\star_{d})$.  
  
  \begin{defi}
  	A {\bf  Manin triple of  differential Novikov algebras}  is a collection $\big( (\frak g\oplus\frak g^{*},\star_{d},\mathcal{B}_{d},D)$,
  	$\frak g,\frak g^{*}  \big)$ such that $\big( (\frak g\oplus\frak g^{*},\star_{d},\mathcal{B}_{d}), \frak g,\frak g^{*}  \big)$ is a Manin triple of Novikov algebras, 
  	$(\frak g\oplus\frak g^{*},\star_{d},\mathcal{B}_{d}, D)$ is a quadratic differential Novikov algebra 
  	 and $D(\frak g)\subset\frak g,\;
  	D(\frak g^{*})\subset\frak g^{*}$.
  \end{defi}

  \begin{lem}\label{lem:451}
  	Let
  	$\big((\frak g\oplus\frak g^{*},\star_d  ,\mathcal{B}_d),
  	\frak g,\frak g^{*}\big)$
  	be a Manin triple of Novikov algebras.
  	\begin{enumerate}
  		\item\label{641}
  		If there exists a linear map
  		$D:\frak g\oplus\frak g^{*}\rightarrow
  		\frak g\oplus\frak g^{*}$
  		such that
  		$\big((\frak g\oplus\frak g^{*},
  		\star_d ,\mathcal{B}_d,D),\frak g,\frak g^{*}\big)$
  		is a Manin triple of differential Novikov algebras,
  		then there exists a linear map
  		$P:\frak g\rightarrow\frak g$ such that $D=P-P^{t}$.
  		
  		\item\label{642}
  		Let $P:\frak g\rightarrow\frak g$ be a linear map.
  		Then
  		$\big((\frak g\oplus\frak g^{*},
  		\star_d ,\mathcal{B}_d,P-P^{t}),\frak g,\frak g^{*}\big)$
  		is a Manin triple of differential Novikov algebras
  		if and only if $P$ is a derivation of
  		$(\frak g,\star_{\frak g})$ and $P^{t}$ is a derivation of
  		$(\frak g^{*},\star_{\frak g^{*}})$.
  	\end{enumerate}
  \end{lem}
  
  \begin{proof}
  	To prove Item \eqref{641}, define
  	$P=D|_{\frak g}$ and $Q=D|_{\frak g^{*}}$.
  	For all $x\in\frak g$ and $a^{*}\in\frak g^{*}$, we have
  	\begin{align*}
  		\langle P^{t}(a^{*}),x\rangle
  		&=\mathcal{B}_d(D(x),a^{*})
  		=-\mathcal{B}_d(x,D(a^{*}))
  		=-\langle Q(a^{*}),x\rangle.
  	\end{align*}
  	Hence $Q=-P^{t}$, and thus $D=P-P^{t}$.
  	To prove Item \eqref{642}, suppose that $P$ and $P^{t}$
  	are derivations of $(\frak g,\star_{\frak g})$
  	and $(\frak g^{*},\star_{\frak g^{*}})$ respectively, and set
  	$D=P-P^{t}$. By \cite{HBG},
  	the product $\star_d $ on $\frak g\oplus\frak g^{*}$ is given by
  	\begin{equation}\label{1065}
  	(x+a^{*})\star_d (y+b^{*})
  	= x\star_{\frak g}y
  	+L^{*}_{\bullet_{\frak g^{*}}}(a^{*})y
  	-R_{\star_{\frak g^{*}}}^{*}(b^{*})x +a^{*}\star_{\frak g^{*}}b^{*}
  	+ L^{*}_{\bullet_{\frak g}}(x)b^{*}
  	-R_{\star_{\frak g}}^{*}(y)a^{*},
  	\end{equation}
  where 
  \begin{equation*}
  x\bullet_{\frak g}y=x\star_{\frak g}y+y\star_{\frak g}x,\;
  a^{*}\bullet_{\frak g^{*}}b^{*}= a^{*}\star_{\frak g^{*}}b^{*}+b^{*}\star_{\frak g^{*}} a^{*},\;\forall x,y\in\frak g, a^{*},b^{*}\in\frak g^{*}.
  \end{equation*}
  	Using the derivation identities of $P$ and $P^{t}$, we obtain
  	\begin{align*}
  		&D(x\star_d  a^{*})
  		-P(x)\star_d  a^{*}
  		+x\star_d  P^{t}(a^{*})\\
  		 &\overset{\eqref{1065}}{=} 
  		-PR_{\star_{\frak g^{*}}}^{*}(a^{*})x
  		+R_{\star_{\frak g^{*}}}^{*}(a^{*})P(x)
  		-R_{\star_{\frak g^{*}}}^{*}(P^{t}(a^{*}))x
  		-P^{t}L^{*}_{\bullet_{\frak g}}(x)a^{*}
  		-L^{*}_{\bullet_{\frak g}}(P(x))a^{*} 
  		+L^{*}_{\bullet_{\frak g}}(x)P^{t}(a^{*})\\
  		&= 0.
  	\end{align*}
  	Similarly we have
  	$
  	D(a^{*}\star_d  x)
  	+P^{t}(a^{*})\star_d  x
  	-a^{*}\star_d  P(x)=0.
  	$
  	Thus $D$ is a derivation of
  	$(\frak g\oplus\frak g^{*},\star_d )$.
  	Moreover, $D$ preserves both summands and 
  	is clearly skew-symmetric with respect to $\mathcal{B}_d$.
  	Therefore 	$\big((\frak g\oplus\frak g^{*},
  	\star_d ,\mathcal{B}_d,P-P^{t}),\frak g,\frak g^{*}\big)$
  	is a Manin triple of differential Novikov algebras. The converse side is obtained similarly.
  \end{proof}
  
  \begin{pro}\label{2013}
  	Let $\big( (\frak g\oplus\frak g^{*},\star_{d},\mathcal{B}_{d},D),\frak g,\frak g^{*}  \big)$ be a Manin triple of differential Novikov algebras. Then $\big( (\frak g\oplus\frak g^{*},[-,-]_{d},\lbrace-,-\rbrace_{d},\mathcal{B}_{d}),   \frak g,  \frak g^{*} \big)$
  	is a  Manin triple of compatible Lie algebras with respect to the invariant-c$2$c form, where the multiplications $[-,-]_{d}$ and $\lbrace-,-\rbrace_{d}$ are respectively given by
  	\begin{align}
  		[u,v]_{d}&=D(u)\star_{d} v-D(v)\star_{d} u,\label{1113}\\
  		\lbrace u,v\rbrace_{d}&=v\star_{d}u-u\star_{d}v,\;\forall u,v\in \frak g\oplus\frak g^{*}.\label{1114}
  	\end{align}
  	Moreover, the compatible LAPL algebra $(\frak g\oplus\frak g^{*},[-,-]_{d},\circ_{d})$ induced from $\mathcal{B}_{d}$ through Proposition \ref{673}
  	is given by \eqref{1113} and 
  	\begin{equation}\label{1961}
  		u\circ_{d}v=u\star_{d}v+2v\star_{d}u,\; \forall u,v\in \frak g\oplus\frak g^{*},
  	\end{equation}
  	which  contains $(\frak g,[-,-]_{\frak g},\circ_{\frak g})$ and $(\frak g^{*},[-,-]_{\frak g^{*}},\circ_{\frak g^{*}})$ as compatible LAPL subalgebras.
  \end{pro}
   \begin{proof}
   	By Proposition \ref{505}, the multiplications \eqref{1113}
   	and \eqref{1114} define a compatible Lie algebra structure $(\frak g\oplus\frak g^{*},[-,-]_{d},\lbrace-,-\rbrace_{d} )$
   	on $\frak g\oplus\frak g^{*}$ with invariant-c$2$c form
   	$\mathcal{B}_d$, and it is straightforward to check that the restrictions on $\frak g$ and $\frak g^{*}$ give rise to compatible Lie subalgebras. Therefore $\big( (\frak g\oplus\frak g^{*},[-,-]_{d},\lbrace-,-\rbrace_{d},\mathcal{B}_{d}),   \frak g,  \frak g^{*} \big)$
   	is a  Manin triple of compatible Lie algebras with respect to the invariant-c$2$c form.
   	By Proposition \ref{751}, there is a compatible LAPL algebra $(\frak g\oplus\frak g^{*},[-,-]_{d},\circ_{d})$ with the multiplication $\circ_{d}$ given by \eqref{1961}
    which is induced from $\mathcal{B}_{d}$ through Proposition \ref{673}.
  It is straightforward to check that the restrictions on $\frak g$ and $\frak g^{*}$ give rise to compatible LAPL subalgebras.
   \end{proof}

  \begin{defi}\cite{HBG}
  	A {\bf Novikov coalgebra} is a pair
  	$(\frak g,\alpha)$, where $\frak g$ is a vector space and
  	$\alpha:\frak g\rightarrow\frak g\otimes\frak g$
  	is a comultiplication satisfying the following equations:
  	\begin{align*}
  		(\mathrm{id}_{\frak g}\otimes\alpha)\alpha(x)
  		-(\tau\otimes\mathrm{id}_{\frak g})(\mathrm{id}_{\frak g}\otimes\alpha)\alpha(x)
  		&=(\alpha\otimes\mathrm{id}_{\frak g})\alpha(x)
  		-(\tau\otimes\mathrm{id}_{\frak g})(\alpha\otimes\mathrm{id}_{\frak g})\alpha(x),\\
  		(\tau\otimes\mathrm{id}_{\frak g})(\mathrm{id}_{\frak g}\otimes\alpha)\tau\alpha(x)
  		&=(\alpha\otimes\mathrm{id}_{\frak g})\alpha(x),
  		\; \forall x\in\frak g.
  	\end{align*}
  A {\bf Novikov bialgebra} is a triple
  $(\frak g,\star_{\frak g},\alpha)$ such that
  $(\frak g,\star_{\frak g})$ is a Novikov algebra,
  $(\frak g,\alpha)$ is a Novikov coalgebra,
  and the following equations hold for all $x,y\in\frak g$:
  \begin{align*}
  	\alpha(x\star_{\frak g}y)&=\big(R_{\star_{\frak g}}(y)\otimes\mathrm{id}_{\frak g}\big)\alpha(x)+\big(\mathrm{id}_{\frak g}\otimes L_{\bullet_{\frak g}}(x)\big)\big(\alpha(y)+\tau\alpha(y)\big),\\
  \big(L_{\bullet_{\frak g}}(x)\otimes\mathrm{id}_{\frak g}\big)\alpha(y)-\big(\mathrm{id}_{\frak g}\otimes L_{\bullet_{\frak g}}(x)\big)\tau\alpha(y)&=\big(L_{\bullet_{\frak g}}(y)\otimes\mathrm{id}_{\frak g}\big)\alpha(x)-\big(\mathrm{id}_{\frak g}\otimes L_{\bullet_{\frak g}}(y)\big)\tau\alpha(x),\\
  \big(\mathrm{id}_{\frak g}\otimes R_{\star_{\frak g}}(x)-R_{\star_{\frak g}}(x)\otimes\mathrm{id}_{\frak g}\big)\big(\alpha(y)+\tau\alpha(y)\big)&=\big(\mathrm{id}_{\frak g}\otimes R_{\star_{\frak g}}(y)-R_{\star_{\frak g}}(y)\otimes\mathrm{id}_{\frak g}\big)\big(\alpha(x)+\tau\alpha(x)\big).
  \end{align*}
  \end{defi}
  
  Now we introduce the notion of a derivation on a Novikov bialgebra.
  
  \begin{defi}
  	Let $(\frak g,\star_{\frak g},\alpha)$ be a Novikov bialgebra.
  	If there is a linear map $P:\frak g\rightarrow\frak g$ which is a derivation of $(\frak g,\star_{\frak g})$ and a coderivation of $(\frak g,\alpha)$ in the sense that
  	\begin{equation}\label{75}
  		\alpha\big(P(x)\big)
  		=\big(P\otimes\mathrm{id}_{\frak g}
  		+\mathrm{id}_{\frak g}\otimes P\big)\alpha(x),
  		\; \forall x\in\frak g,
  	\end{equation}
  	then we say $P$ is a {\bf derivation} of a Novikov bialgebra and $(\frak g,\star_{\frak g},\alpha,P)$ is a {\bf differential Novikov bialgebra}. 
  \end{defi} 
  
  \begin{pro}\label{2101}
  	Let $(\frak g,\star_{\frak g}), (\frak g^{*},\star_{\frak g^{*}})$ be Novikov algebras, $\alpha$ be the linear dual of $\star_{\frak g^{*}}$ and $P:\frak g\rightarrow\frak g$ be a linear map. Then $(\frak g,\star_{\frak g},\alpha,P)$ is a differential Novikov bialgebra if and only if there is a Manin triple of differential Novikov algebras $\big( (\frak g\oplus\frak g^{*},\star_{d},\mathcal{B}_{d},P-P^{t}),\frak g,\frak g^{*}  \big)$.
  \end{pro}
  \begin{proof}
  	By \cite{HBG}, $(\frak g,\star_{\frak g},\alpha )$ is a  Novikov bialgebra if and only if there is a Manin triple of   Novikov algebras $\big( (\frak g\oplus\frak g^{*},\star_{d},\mathcal{B}_{d} ),\frak g,\frak g^{*}  \big)$.
  	By a straightforward computation, $P$ is a coderivation of $(\frak g,\alpha)$ if and only if $P^{t}$ is a derivation of $(\frak g^{*},\star_{\frak g^{*}})$. Hence the conclusion follows from Lemma \ref{lem:451}.
  \end{proof}

  \begin{pro}\label{prop:differential-novikov-lapl}
  Let $(\frak g,\star_{\frak g},\alpha,P)$ be a differential Novikov bialgebra.
  Then there is
  	a compatible LAPL bialgebra
  	$(\frak g,[-,-]_{\frak g},\circ_{\frak g},\delta,\theta)$,
  	where
  	$[-,-]_{\frak g},\circ_{\frak g}:
  	\frak g\otimes\frak g\rightarrow\frak g$
  	are multiplications given by \eqref{331} and \eqref{408} respectively, and
  	$\delta,\theta:\frak g\rightarrow\frak g\otimes\frak g$
  	are respectively given by
  	\begin{align*}
  		\delta(x)
  		&=(\tau-\mathrm{id}^{\otimes 2}_{ \frak g})
  		(P\otimes\mathrm{id}_{\frak g})\alpha(x),
  		\\
  		\theta(x)
  		&=(\mathrm{id}^{\otimes 2}_{\frak g}+2\tau)\alpha(x),
  		\;\forall x\in\frak g.
  	\end{align*}
  \end{pro}
  \begin{proof}
  	The conclusion can be obtained by a direct verification.
  	Here we give a proof using Manin triples.
  	Let $(\frak g,\star_{\frak g},\alpha,P)$ be a
  	differential Novikov bialgebra.
  	By Proposition \ref{2101}, there is a Manin triple
  	of differential Novikov algebras
  $
  	\bigl(
  	(\frak g\oplus\frak g^{*},
  	\star_{d},\mathcal{B}_{d},P-P^{t}),
  	\frak g,\frak g^{*}
  	\bigr).
   $
  	By Proposition \ref{2013},
  	 $
  	\bigl(
  	(\frak g\oplus\frak g^{*},
  	[-,-]_{d},\lbrace-,-\rbrace_{d},\mathcal{B}_{d}),
  	\frak g,\frak g^{*}
  	\bigr)
   $
  	is a Manin triple of compatible Lie algebras with
  	respect to the invariant-c$2$c form $\mathcal{B}_{d}$,
  	where the multiplications $[-,-]_{d}$ and
  	$\lbrace-,-\rbrace_{d}$ are respectively given by
  	\eqref{1113} and \eqref{1114}.
  	Moreover, the compatible LAPL algebra
  	$(\frak g\oplus\frak g^{*},[-,-]_{d},\circ_{d})$
  	induced from $\mathcal{B}_{d}$ through
  	Proposition \ref{673} contains
  	$(\frak g,[-,-]_{\frak g},\circ_{\frak g})$ and
  	$(\frak g^{*},[-,-]_{\frak g^{*}},\circ_{\frak g^{*}})$
  	as compatible LAPL subalgebras.
  	Here $[-,-]_{\frak g}$ and $\circ_{\frak g}$ are
  	given by \eqref{331} and \eqref{408}, while
  	$[-,-]_{\frak g^{*}},\circ_{\frak g^{*}}:
  	\frak g^{*}\otimes\frak g^{*}\rightarrow\frak g^{*}$
  	are given by
  	\begin{align*}
  		[a^{*},b^{*}]_{\frak g^{*}}
  		&=-P^{t}(a^{*})\star_{\frak g^{*}}b^{*}
  		+P^{t}(b^{*})\star_{\frak g^{*}}a^{*},
  		\\
  		a^{*}\circ_{\frak g^{*}}b^{*}
  		&=a^{*}\star_{\frak g^{*}}b^{*}
  		+2b^{*}\star_{\frak g^{*}}a^{*},
  		\;\forall a^{*},b^{*}\in\frak g^{*}.
  	\end{align*}
  By Corollary \ref{1265}, there is a compatible LAPL bialgebra
  $(\frak g,[-,-]_{\frak g},\circ_{\frak g},\delta,\theta)$, in which $\delta$ is the linear dual of $[-,-]_{\frak g^{*}}$ and $\theta$ is the linear dual of $\circ_{\frak g^{*}}$. Hence the conclusion follows.
  \end{proof}

Let $(\frak g,\star_{\frak g})$ be a Novikov algebra and $r=\sum\limits_i a_i\otimes b_i\in\frak g\otimes\frak g$. Set
\begin{equation}\label{eq:novikov-ybe-tensor}
	\textbf{N}(r)=r_{13}\star_{\frak g}r_{23}+r_{12}\bullet_{\frak g}r_{23}+r_{13}\star_{\frak g}r_{12},
\end{equation}
where
\begin{equation*}
	r_{13}\star_{\frak g}r_{23}=\sum_{i,j}a_i\otimes a_j\otimes b_i\star_{\frak g}b_j,\;
	r_{12}\bullet_{\frak g}r_{23}=\sum_{i,j}a_i\otimes b_i\bullet_{\frak g}a_j\otimes b_j,\;
	r_{13}\star_{\frak g}r_{12}=\sum_{i,j}a_i\star_{\frak g}a_j\otimes b_j\otimes b_i.
\end{equation*}
If the following equation holds:
\begin{equation}\label{2265}
	\textbf{N}(r)=0,
\end{equation}
then we say $r$ is a solution of the {\bf Novikov Yang-Baxter equation} (NYBE in short) \cite{HBG}
in $(\frak g,\star_{\frak g})$.

\begin{lem}\cite{HBG}\label{2271}
	Let $(\frak g,\star_{\frak g})$ be a Novikov algebra, $r\in\frak g\otimes\frak g$
	and
	$\alpha_r:\frak g\rightarrow\frak g\otimes\frak g$ be given by
	\begin{equation}\label{2275}
		\alpha_r(x)=-\big(L_{\star_{\frak g}}(x)\otimes\mathrm{id}_{\frak g}+\mathrm{id}_{\frak g}\otimes L_{\bullet_{\frak g}}(x)\big)r,\;\forall x\in\frak g.
	\end{equation}
	If $r$ is an antisymmetric solution of the NYBE in $(\frak g,\star_{\frak g})$, then $(\frak g,\star_{\frak g},\alpha_r)$ is a Novikov bialgebra, which is called a {\bf triangular Novikov bialgebra}.
\end{lem}

\begin{rmk}
	The comultiplication $\alpha_r$ in \eqref{2275}
	differs by an overall minus sign from that used in
	\cite{HBG}.
	Replacing the comultiplication by its negative
	preserves the Novikov bialgebra axioms.
	This sign convention is chosen so that, in
	Theorem \ref{thm:triangular-novikov-lapl}, the induced
	compatible LAPL bialgebra has comultiplications
	$\delta_r$ and $\theta_r$ associated with the same
	tensor $r$.
\end{rmk}

\begin{lem}\label{2281}
	Let $(\frak g,\star_{\frak g},P)$ be a differential Novikov algebra and $r\in\frak g\otimes\frak g$.
	Let $\alpha_r:\frak g\rightarrow\frak g\otimes\frak g$ be given by \eqref{2275}.
	Then \eqref{75} holds if and only if the following equation holds:
	\begin{equation*}
		\big(L_{\star_{\frak g}}(x)\otimes\mathrm{id}_{\frak g}+\mathrm{id}_{\frak g}\otimes L_{\bullet_{\frak g}}(x)\big)\big(P\otimes\mathrm{id}_{\frak g}+\mathrm{id}_{\frak g}\otimes P\big)r=0,\;\forall x\in\frak g.
	\end{equation*}
	In particular, if the following equation holds:
	\begin{equation}\label{2289}
		\big(P\otimes\mathrm{id}_{\frak g}+\mathrm{id}_{\frak g}\otimes P\big)r=0,
	\end{equation}
	then \eqref{75} holds.
\end{lem}
\begin{proof}
	It follows from a straightforward verification.
\end{proof}

\begin{defi}
	Let $(\frak g,\star_{\frak g},P)$ be a differential Novikov algebra and $r\in\frak g\otimes\frak g$.
	If \eqref{2265} and \eqref{2289} hold, then we say $r$ is a solution of the {\bf differential NYBE} in $(\frak g,\star_{\frak g},P)$.
\end{defi}

\begin{pro}\label{2303}
	Let $(\frak g,\star_{\frak g},P)$ be a differential Novikov algebra, $r\in\frak g\otimes\frak g$
	and
	$\alpha_r:\frak g\rightarrow\frak g\otimes\frak g$ be given by \eqref{2275}.
	If $r$ is an antisymmetric solution of the differential NYBE in $(\frak g,\star_{\frak g},P)$, then $(\frak g,\star_{\frak g},\alpha_r,P)$ is a differential Novikov bialgebra, which is called a {\bf triangular differential Novikov bialgebra}.
\end{pro}
\begin{proof}
	It follows from Lemmas \ref{2271} and \ref{2281} directly.
\end{proof}

\begin{pro}\label{prop:novikov-clapl-ybe}
	Let $(\frak g,\star_{\frak g},P)$ be a differential Novikov algebra and let $(\frak g,[-,-]_{\frak g},\circ_{\frak g})$ be the compatible LAPL algebra
	with $[-,-]_{\frak g}$ and $\circ_{\frak g}$ given by \eqref{331} and \eqref{408}.
	Let $r=\sum\limits_{i}a_{i}\otimes b_{i}\in\frak g\otimes\frak g$. 
	Then the following identities hold:
\begin{align}
	\textbf{T}(r)
	&=\big(\mathrm{id}_{\frak g}^{\otimes3}
	-\tau\otimes\mathrm{id}_{\frak g}\big)\textbf{N}(r) +\sum_j\Big(
	\big(L_{\bullet_{\frak g}}(a_j)\otimes\mathrm{id}_{\frak g}
	+\mathrm{id}_{\frak g}\otimes L_{\star_{\frak g}}(a_j)\big)
	\big(r+\tau(r)\big)
	\Big)\otimes b_j,\label{eq:novikov-apl-ybe-bridge}\\
	\textbf{C}(r)
	&=\big(\tau\otimes\mathrm{id}_{\frak g}
	-\mathrm{id}_{\frak g}^{\otimes3}\big)
	\big(P\otimes\mathrm{id}_{\frak g}\otimes\mathrm{id}_{\frak g}\big)
	\textbf{N}(r)\notag\\
	&\quad+\sum_j\Big(
	\big(L_{\bullet_{\frak g}}(a_j)\otimes\mathrm{id}_{\frak g}
	+\mathrm{id}_{\frak g}\otimes L_{\bullet_{\frak g}}(a_j)\big)
	\big(P\otimes\mathrm{id}_{\frak g}
	+\mathrm{id}_{\frak g}\otimes P\big)r
	\Big)\otimes b_j\notag\\
	&\quad\quad- \Big(
	\big(L_{\bullet_{\frak g}}(a_j)\otimes P
	+\mathrm{id}_{\frak g}\otimes PL_{\star_{\frak g}}(a_j)\big)
	\big(r+\tau(r)\big)
	\Big)\otimes b_j\notag\\
	&\quad\quad+ \big(\mathrm{id}_{\frak g}\otimes\tau\big)
	\Big(
	\Big(
	\big(\mathrm{id}_{\frak g}\otimes R_{\star_{\frak g}}(b_j)\big)
	\big(P\otimes\mathrm{id}_{\frak g}
	+\mathrm{id}_{\frak g}\otimes P\big)r
	\Big)\otimes a_j
	\Big)\notag\\
	&\quad\quad-  a_j\otimes\Big(
	\big(\mathrm{id}_{\frak g}\otimes R_{\star_{\frak g}}(b_j)\big)
	\big(P\otimes\mathrm{id}_{\frak g}
	+\mathrm{id}_{\frak g}\otimes P\big)r
	\Big).\label{eq:novikov-cybe-bridge}
\end{align}
\end{pro}
\begin{proof}
	By \eqref{1318} and \eqref{eq:novikov-ybe-tensor}, we have
	\begin{align*}
		 \big(\tau\otimes\mathrm{id}_{\frak g}
		-\mathrm{id}_{\frak g}^{\otimes3}\big)\textbf{N}(r)+\textbf{T}(r) 
		&=\sum_{i,j}\Big(
		(a_i\bullet_{\frak g}a_j)\otimes b_i\otimes b_j
		+(b_i\bullet_{\frak g}a_j)\otimes a_i\otimes b_j\\
		&\;\;
		+a_i\otimes(a_j\star_{\frak g}b_i)\otimes b_j
		+b_i\otimes(a_j\star_{\frak g}a_i)\otimes b_j
		\Big)\\
		&=\sum_j\Big(
		\big(L_{\bullet_{\frak g}}(a_j)\otimes\mathrm{id}_{\frak g}
		+\mathrm{id}_{\frak g}\otimes L_{\star_{\frak g}}(a_j)\big)
		\big(r+\tau(r)\big)
		\Big)\otimes b_j.
	\end{align*}
	This proves \eqref{eq:novikov-apl-ybe-bridge}.
A similar but slightly more complicated computation yields \eqref{eq:novikov-cybe-bridge}.
\end{proof}
  
Combining Proposition \ref{prop:differential-novikov-lapl} with the above results,
we obtain the following construction of triangular compatible LAPL bialgebras.

\begin{thm}\label{thm:triangular-novikov-lapl}
	Let $(\frak g,\star_{\frak g},P)$ be a differential Novikov algebra,
	$r\in\frak g\otimes\frak g$ be an antisymmetric solution of the
	differential NYBE in $(\frak g,\star_{\frak g},P)$ and
	$(\frak g,\star_{\frak g},\alpha_r,P)$ be the triangular differential
	Novikov bialgebra given by Proposition \ref{2303}.
	Let $(\frak g,[-,-]_{\frak g},\circ_{\frak g},\delta,\theta)$
	be the corresponding compatible LAPL bialgebra given by
	Proposition \ref{prop:differential-novikov-lapl}.
	Then we have
	\begin{align}
		\delta(x)
		&=\big(
		\mathrm{ad}_{\frak g}(x)\otimes\mathrm{id}_{\frak g}
		+\mathrm{id}_{\frak g}\otimes\mathrm{ad}_{\frak g}(x)
		\big)r
		=\delta_r(x),\label{eq:triangular-induced-delta}\\
		\theta(x)
		&=\big(
		L_{\circ_{\frak g}}(x)\otimes\mathrm{id}_{\frak g}
		-\mathrm{id}_{\frak g}\otimes\mathfrak{ad}_{\frak g}(x)
		\big)r
		=\theta_r(x),\;\forall x\in\frak g.
		\label{eq:triangular-induced-theta}
	\end{align}
	Therefore
	$(\frak g,[-,-]_{\frak g},\circ_{\frak g},\delta,\theta)$
	is also triangular.
\end{thm}
\begin{proof}
	For all $x\in\frak g$, we have
	\begin{align*}
		\delta(x)
		&=(\tau-\mathrm{id}_{\frak g}^{\otimes 2})
		(P\otimes\mathrm{id}_{\frak g})\alpha_r(x)\overset{\eqref{2289}}{=}\Big(
		\big(
		L_{\star_{\frak g}}(P(x))-R_{\star_{\frak g}}(x)P
		\big)\otimes\mathrm{id}_{\frak g}
		+\mathrm{id}_{\frak g}\otimes
		\big(
		L_{\star_{\frak g}}(P(x))-R_{\star_{\frak g}}(x)P
		\big)
		\Big)r\\
		&=\big(
		\mathrm{ad}_{\frak g}(x)\otimes\mathrm{id}_{\frak g}
		+\mathrm{id}_{\frak g}\otimes\mathrm{ad}_{\frak g}(x)
		\big)r\\
		&=\delta_r(x).
	\end{align*}
	Therefore \eqref{eq:triangular-induced-delta} holds, and
	\eqref{eq:triangular-induced-theta} holds similarly.
	By Proposition \ref{prop:novikov-clapl-ybe},
	$r$ is an antisymmetric solution of the CLAPL-YBE.
	Therefore $(\frak g,[-,-]_{\frak g},\circ_{\frak g},\delta,\theta)$
	is triangular.
\end{proof}

Recall \cite{HBG} that a {\bf representation} of a Novikov algebra
$(\frak g,\star_{\frak g})$ is a triple $(l_{\star},r_{\star},V)$, where $V$ is a vector
space and $l_{\star},r_{\star}:\frak g\rightarrow\mathrm{End}_{\mathbb{K}}(V)$ are linear
maps satisfying the following equations:
\begin{align*}
	l_{\star}(x\star_{\frak g}y-y\star_{\frak g}x)
	&=l_{\star}(x)l_{\star}(y)-l_{\star}(y)l_{\star}(x),\\
	l_{\star}(x)r_{\star}(y)-r_{\star}(y)l_{\star}(x)
	&=r_{\star}(x\star_{\frak g}y)-r_{\star}(y)r_{\star}(x),\\
	l_{\star}(x\star_{\frak g}y)&=r_{\star}(y)l_{\star}(x),\\
	r_{\star}(x)r_{\star}(y)&=r_{\star}(y)r_{\star}(x),\;\forall x,y\in\frak g.
\end{align*}

\begin{defi}\label{def:differential-novikov-representation}
	Let $(\frak g,\star_{\frak g},P)$ be a differential Novikov algebra
	and $(l_{\star},r_{\star},V)$ be a representation of $(\frak g,\star_{\frak g})$.
	If there is a linear map $\alpha:V\rightarrow V$ satisfying the following
	equations:
	\begin{align}
		\alpha\big(l_{\star}(x)v\big)-l_{\star}\big(P(x)\big)v-l_{\star}(x)\alpha(v)&=0,
		\label{eq:differential-novikov-left-action}\\
		\alpha\big(r_{\star}(x)v\big)-r_{\star}\big(P(x)\big)v-r_{\star}(x)\alpha(v)&=0,
		\;\forall x\in\frak g, v\in V,
		\label{eq:differential-novikov-right-action}
	\end{align}
	then we say $(l_{\star},r_{\star},\alpha,V)$ is a {\bf representation} of
	$(\frak g,\star_{\frak g},P)$.
\end{defi}

\begin{ex}
	Let $(\frak g,\star_{\frak g},P)$ be a differential Novikov algebra.
	Then $(L_{\star_{\frak g}},R_{\star_{\frak g}},P,\frak g)$ is a
	representation of $(\frak g,\star_{\frak g},P)$, which is called the
	{\bf adjoint representation} of $(\frak g,\star_{\frak g},P)$.
\end{ex}

\begin{pro}\label{prop:differential-novikov-semidirect}
	Let $(\frak g,\star_{\frak g},P)$ be a differential Novikov algebra,
	$V$ be a vector space and
	$l_{\star},r_{\star}:\frak g\rightarrow\mathrm{End}_{\mathbb{K}}(V)$,
	$\alpha:V\rightarrow V$ be linear maps.
	Then $(l_{\star},r_{\star},\alpha,V)$ is a representation of
	$(\frak g,\star_{\frak g},P)$ if and only if there is a differential
	Novikov algebra structure $(\frak g\oplus V,\star_d  ,P+\alpha)$ on
	$\frak g\oplus V$ with $\star_d  $ given by
	\begin{equation*}
		(x+u)\star_d  (y+v)
		=x\star_{\frak g}y+l_{\star}(x)v+r_{\star}(y)u,
		\;\forall x,y\in\frak g,\ u,v\in V.
	\end{equation*}
	In this case, we denote this differential Novikov algebra structure by
	$(\frak g\ltimes_{l_{\star},r_{\star}}V,P+\alpha)$.
\end{pro}
\begin{proof}
	By \cite[Proposition 3.2]{HBG}, $(l_{\star},r_{\star},V)$ is a representation of
	$(\frak g,\star_{\frak g})$ if and only if $(\frak g\oplus V,\star_d  )$
	is a Novikov algebra.
	Moreover, since $P$ is a derivation of $(\frak g,\star_{\frak g})$,
	$P+\alpha$ is a derivation of $(\frak g\oplus V,\star_d  )$ if and only if
	\eqref{eq:differential-novikov-left-action} and
	\eqref{eq:differential-novikov-right-action} hold.
	The conclusion follows.
\end{proof}

\begin{pro}\label{prop:differential-novikov-lapl-representation}
	Let $(\frak g,\star_{\frak g},P)$ be a differential Novikov algebra
	and $(l_{\star},r_{\star},\alpha,V)$ be a representation of
	$(\frak g,\star_{\frak g},P)$.
	Let $(\frak g,[-,-]_{\frak g},\circ_{\frak g})$ be the compatible
	LAPL algebra given by Proposition \ref{751}.
	Define linear maps
	$\rho,l_{\circ},r_{\circ}:
	\frak g\rightarrow\mathrm{End}_{\mathbb{K}}(V)$ by
	\begin{equation}	\label{eq:novikov-lapl-representation-rho}
		\rho(x)v =l_{\star}(P(x))v-r_{\star}(x)\alpha(v),\;
	l_{\circ}(x)v =l_{\star}(x)v+2r_{\star}(x)v,\;
	r_{\circ}(x)v=r_{\star}(x)v+2l_{\star}(x)v, 
	\end{equation}
for all $x\in\frak g,v\in V$.
	Then $(\rho,l_{\circ},r_{\circ},V)$ is a representation of
	$(\frak g,[-,-]_{\frak g},\circ_{\frak g})$.
\end{pro}

\begin{proof}
	By Proposition \ref{prop:differential-novikov-semidirect}, there is a differential Novikov algebra
	$(\frak g\ltimes_{l_{\star},r_{\star}}V,P+\alpha)$.
	By Proposition \ref{751}, there is a compatible LAPL algebra
	$(\frak g\oplus V,[-,-]_d,\circ_d)$, where
	\begin{align*}
		[x+u,y+v]_d
		&=(P(x)+\alpha(u))\star_d(y+v)
		-(P(y)+\alpha(v))\star_d(x+u)\\
		&=P(x)\star_{\frak g}y-P(y)\star_{\frak g}x
		+l_{\star}(P(x))v+r_{\star}(y)\alpha(u) 
		-l_{\star}(P(y))u-r_{\star}(x)\alpha(v)\\
		&=[x,y]_{\frak g}+\rho(x)v-\rho(y)u,\\
		(x+u)\circ_d(y+v)
		&=(x+u)\star_d(y+v)+2(y+v)\star_d(x+u)\\
		&=x\star_{\frak g}y+2y\star_{\frak g}x
		+l_{\star}(x)v+r_{\star}(y)u
		+2l_{\star}(y)u+2r_{\star}(x)v\\
		&=x\circ_{\frak g}y+l_{\circ}(x)v+r_{\circ}(y)u,
	\end{align*}
	for all $x,y\in\frak g$ and $u,v\in V$.
	Therefore $(\rho,l_{\circ},r_{\circ},V)$ is a representation of
	$(\frak g,[-,-]_{\frak g},\circ_{\frak g})$.
\end{proof}
 
\begin{pro}\label{prop:differential-novikov-dual-representation}
	Let $(\frak g,\star_{\frak g},P)$ be a differential Novikov algebra
	and $(l_{\star},r_{\star},\alpha,V)$ be a representation of
	$(\frak g,\star_{\frak g},P)$.
	Then $(l_{\star}^*+r_{\star}^*,-r_{\star}^*,-\alpha^{t},V^*)$ is also a representation of
	$(\frak g,\star_{\frak g},P)$.
	In particular, the quadruple
	$(L_{\bullet_{\frak g}}^*,-R_{\star_{\frak g}}^*,-P^{t},\frak g^*)$
	is a representation of $(\frak g,\star_{\frak g},P)$.
\end{pro}
\begin{proof}
	By \cite[Proposition 3.3]{HBG},
	$(l_{\star}^*+r_{\star}^*,-r_{\star}^*,V^*)$
	is a representation of $(\frak g,\star_{\frak g})$.
	For all $x\in\frak g$, $v\in V$ and $v^*\in V^*$,
	we have
	\begin{align*}
		&\bigl\langle
		(-\alpha^t)(l_{\star}^*+r_{\star}^*)(x)v^*
		-(l_{\star}^*+r_{\star}^*)(P(x))v^*
		-(l_{\star}^*+r_{\star}^*)(x)(-\alpha^t)(v^*),
		v\bigr\rangle\\
		&=\bigl\langle v^*,
		(l_{\star}+r_{\star})(x)\alpha(v)
		+(l_{\star}+r_{\star})(P(x))v
		-\alpha\big((l_{\star}+r_{\star})(x)v\big)
		\bigr\rangle\\
		&\overset{\eqref{eq:differential-novikov-left-action},
			\eqref{eq:differential-novikov-right-action}}{=}0.
	\end{align*}
Therefore we have
\begin{equation*}
	(-\alpha^t)(l_{\star}^*+r_{\star}^*)(x)v^*
-(l_{\star}^*+r_{\star}^*)(P(x))v^*
-(l_{\star}^*+r_{\star}^*)(x)(-\alpha^t)(v^*)=0,
\end{equation*}
that is,
\eqref{eq:differential-novikov-left-action} holds for $(l_{\star}^*+r_{\star}^*,-r_{\star}^*,-\alpha^{t},V^*)$. Similarly, \eqref{eq:differential-novikov-right-action} holds for $(l_{\star}^*+r_{\star}^*,-r_{\star}^*,-\alpha^{t},V^*)$.
	Therefore
	$(l_{\star}^*+r_{\star}^*,-r_{\star}^*,-\alpha^t,V^*)$
	is a representation of $(\frak g,\star_{\frak g},P)$.
	Applying this result to the adjoint
	representation
	$(L_{\star_{\frak g}},R_{\star_{\frak g}},P,\frak g)$
	gives the remaining conclusion.
\end{proof}

\begin{defi}\label{def:differential-novikov-o-operator}
	Let $(\frak g,\star_{\frak g},P)$ be a differential Novikov algebra
	and $(l_{\star},r_{\star},\alpha,V)$ be a representation of
	$(\frak g,\star_{\frak g},P)$.
	Suppose that there is a linear map $T:V\rightarrow\frak g$ which is
	an $\mathcal{O}$-operator \cite{HBG} of $(\frak g,\star_{\frak g})$
	associated to $(l_{\star},r_{\star},V)$, that is, the following equation holds:
	\begin{equation*}
		T(u)\star_{\frak g}T(v)
	=T\big(l_{\star}(T(u))v+r_{\star}(T(v))u\big),
	\;\forall u,v\in V.
	\end{equation*}
	If in addition, the following equation holds:
	\begin{equation}\label{2632}
		PT=T\alpha,
	\end{equation}
	then we say $T$ is an {\bf $\mathcal{O}$-operator} of
	$(\frak g,\star_{\frak g},P)$ associated to $(l_{\star},r_{\star},\alpha,V)$.
\end{defi}

\begin{pro}\label{prop:differential-novikov-lapl-o-operator}
	Let $(\frak g,\star_{\frak g},P)$ be a differential Novikov algebra,
	$(l_{\star},r_{\star},\alpha,V)$ be a representation of
	$(\frak g,\star_{\frak g},P)$ and $T:V\rightarrow\frak g$ be an
	$\mathcal{O}$-operator of $(\frak g,\star_{\frak g},P)$ associated to
	$(l_{\star},r_{\star},\alpha,V)$.
		Let $(\frak g,[-,-]_{\frak g},\circ_{\frak g})$ be the compatible
	LAPL algebra given by Proposition \ref{751}.
	Then $T$ is an $\mathcal{O}$-operator of the compatible LAPL algebra
	$(\frak g,[-,-]_{\frak g},\circ_{\frak g})$ associated to
	$(\rho,l_{\circ},r_{\circ},V)$, where
 $\rho,l_{\circ},r_{\circ}$ are given by
	\eqref{eq:novikov-lapl-representation-rho}.
\end{pro}
\begin{proof}
It follows from a direct computation.
\end{proof}

\begin{pro}\label{prop:differential-nybe-operator-form}
	Let $(\frak g,\star_{\frak g},P)$ be a differential Novikov algebra
	and $r\in\frak g\otimes\frak g$ be antisymmetric.
	Then $r$ is a solution of the differential NYBE in
	$(\frak g,\star_{\frak g},P)$ if and only if the following equations hold:
	\begin{align}
		r(a^*)\star_{\frak g}r(b^*)
		&=r\Big(
		L_{\bullet_{\frak g}}^*\big(r(a^*)\big)b^*
		-R_{\star_{\frak g}}^*\big(r(b^*)\big)a^*
		\Big),
		\;\forall a^*,b^*\in\frak g^*,\nonumber\\
		Pr&=-rP^{t},\label{2693}
	\end{align}
	that is, $r$ is an $\mathcal{O}$-operator of
	$(\frak g,\star_{\frak g},P)$ associated to
	$(L_{\bullet_{\frak g}}^*,-R_{\star_{\frak g}}^*,-P^{t},\frak g^*)$.
\end{pro}
\begin{proof}
	By \cite[Theorem 3.27]{HBG}, the first equation is equivalent to
	$\textbf{N}(r)=0$.
	Moreover, it is straightforward to check that \eqref{2289} holds if and only if \eqref{2693} holds.
	The conclusion follows.
\end{proof}

\begin{thm}\label{thm:differential-novikov-o-operator-construction}
	Let $(\frak g,\star_{\frak g},P)$ be a differential Novikov algebra
	and $(l_{\star},r_{\star},\alpha,V)$ be a representation of
	$(\frak g,\star_{\frak g},P)$.
	Let $T:V\rightarrow\frak g$ be a linear map which is identified as an
	element in
	\[
	V^*\otimes\frak g
	\subset
	\big(\frak g\ltimes_{l_{\star}^*+r_{\star}^*,-r_{\star}^*}V^*\big)
	\otimes
	\big(\frak g\ltimes_{l_{\star}^*+r_{\star}^*,-r_{\star}^*}V^*\big).
	\]
	Then $T-\tau(T)$ is an antisymmetric solution of the differential NYBE in
	$(\frak g\ltimes_{l_{\star}^*+r_{\star}^*,-r_{\star}^*}V^*,P-\alpha^{t})$
	if and only if $T$ is an $\mathcal{O}$-operator of
	$(\frak g,\star_{\frak g},P)$ associated to $(l_{\star},r_{\star},\alpha,V)$.
\end{thm}
\begin{proof}
	Let $\{v_1,\ldots,v_n\}$ be a basis of $V$ and
	$\{v_1^*,\ldots,v_n^*\}$ be its dual basis.
	Then $T=\sum\limits_i v_i^*\otimes T(v_i)$.
	By \cite[Theorem 3.29]{HBG}, $T-\tau(T)$ is a solution of the NYBE in
	$\frak g\ltimes_{l_{\star}^*+r_{\star}^*,-r_{\star}^*}V^*$ if and only if $T$ is an
	$\mathcal{O}$-operator of $(\frak g,\star_{\frak g})$ associated to
	$(l_{\star},r_{\star},V)$.
	Moreover, we have
	\begin{align*}
		&\Big(
		(P-\alpha^{t})\otimes\mathrm{id}_{\frak g\oplus V^*}
		+\mathrm{id}_{\frak g\oplus V^*}\otimes(P-\alpha^{t})
		\Big)\big(T-\tau(T)\big)\\
		&=\sum_i v_i^*\otimes
		\big(P(T(v_i))-T(\alpha(v_i))\big)
		-\sum_i\big(P(T(v_i))-T(\alpha(v_i))\big)\otimes v_i^*,
	\end{align*}
which vanishes
	if and only if \eqref{2632} holds.
	The conclusion follows.
\end{proof}
 Recall \cite{HBG} that a {\bf pre-Novikov algebra} is a triple $(\frak g,\rhd_{\frak g},\lhd_{\frak g})$, where $\frak g$ is a vector space and $\rhd_{\frak g},\lhd_{\frak g}:\frak g\otimes\frak g\rightarrow\frak g$ are multiplications such that the following equations hold:
 \begin{align*}
 	 x\rhd_{\frak g}(y\rhd_{\frak g}z)&=(x\star_{\frak g} y)\rhd_{\frak g}z+y\rhd_{\frak g}(x\rhd_{\frak g}z)-(y\star_{\frak g} x)\rhd_{\frak g}z, \\
  x\rhd_{\frak g}(y\lhd_{\frak g}z)&=(x\rhd_{\frak g}y)\lhd_{\frak g}z+y\lhd_{\frak g}(x\star_{\frak g}z)-(y\lhd_{\frak g} x)\lhd_{\frak g}z,   \\
  (x\star_{\frak g} y)\rhd_{\frak g}z&=(x\rhd_{\frak g}z)\lhd_{\frak g}y, \\
  (x\lhd_{\frak g}y)\lhd_{\frak g}z&=(x\lhd_{\frak g}z)\lhd_{\frak g}y,\;\forall x,y ,z\in\frak g, 
 \end{align*}
 where $\star_{\frak g}:\frak g\otimes\frak g\rightarrow\frak g$ is given by
 \begin{align}\label{1954}
 	x\star_{\frak g} y=x\rhd_{\frak g}y+x\lhd_{\frak g}y,\;\forall x,y\in\frak g.
 \end{align}
 	Moreover, if $P:\frak g\rightarrow\frak g$ is a linear map which is a derivation of $(\frak g,\rhd_{\frak g},\lhd_{\frak g})$ in the sense that
 \begin{equation*}
 	P(x\rhd_{\frak g} y)=P(x)\rhd_{\frak g} y+x\rhd_{\frak g}P(y),\;
 	P(x\lhd_{\frak g} y)=P(x)\lhd_{\frak g} y+x\lhd_{\frak g}P(y),\;\forall x,y\in\frak g,
 \end{equation*}
then we say $(\frak g,\rhd_{\frak g},\lhd_{\frak g},P)$ is a {\bf differential pre-Novikov algebra}.

 \begin{pro}\label{1939}
 	Let $(\frak g,\rhd_{\frak g},\lhd_{\frak g})$ be a pre-Novikov algebra. Then $(\frak g,\star_{\frak g})$ with $\star_{\frak g}:\frak g\otimes\frak g\rightarrow\frak g$ given by
 	\eqref{1954} is a Novikov algebra.
 In particular, if there is a linear map $P:\frak g\rightarrow\frak g$ such that $(\frak g,\rhd_{\frak g},\lhd_{\frak g},P)$ is a  differential pre-Novikov algebra, then $(\frak g,\star_{\frak g},P)$ is a differential Novikov algebra.
 \end{pro}
 \begin{proof}
 	The first conclusion follows from \cite{HBG} and the second conclusion follows from a direct check.
 \end{proof}

 \begin{pro}\label{prop:differential-pre-novikov-bialgebra}
 	Let $(\frak g,\rhd_{\frak g},\lhd_{\frak g},P)$ be a differential
 	pre-Novikov algebra and $(\frak g,\star_{\frak g},P)$ be the
 	differential Novikov algebra given by Proposition \ref{1939}.
 	Then the identity map $\mathrm{id}_{\frak g}$ is an
 	$\mathcal{O}$-operator of $(\frak g,\star_{\frak g},P)$ associated to
 	$(L_{\rhd_{\frak g}},R_{\lhd_{\frak g}},P,\frak g)$.
 	Consequently,
 	\begin{equation*}
 		r=\mathrm{id}_{\frak g}-\tau(\mathrm{id}_{\frak g})
 		=\sum_{i=1}^{n}
 		\big(e_i^*\otimes e_i-e_i\otimes e_i^*\big)
 	\end{equation*}
 	is an antisymmetric solution of the differential NYBE in
 	\[
 	\big(
 	\frak g\ltimes_{L_{\rhd_{\frak g}}^*+R_{\lhd_{\frak g}}^*,
 		-R_{\lhd_{\frak g}}^*}\frak g^*,P-P^t
 	\big),
 	\]
 	where $\{e_1,\ldots,e_n\}$ is a basis of $\frak g$ and
 	$\{e_1^*,\ldots,e_n^*\}$ is its dual basis.
 	Further, there is a triangular differential Novikov bialgebra
 	$(\frak g\ltimes_{L_{\rhd_{\frak g}}^*+R_{\lhd_{\frak g}}^*,
 		-R_{\lhd_{\frak g}}^*}\frak g^*,\alpha_r,P-P^t)$,
 	where  $\alpha_r$ is defined by \eqref{2275} for
 	$(\frak g\oplus\frak g^*,\star_d  )$.
 \end{pro}
 \begin{proof}
 	It is straightforward to check that $(L_{\rhd_{\frak g}},R_{\lhd_{\frak g}},P,\frak g)$ is a representation
 	of $(\frak g,\star_{\frak g},P)$, and 
 	  $\mathrm{id}_{\frak g}$ is an $\mathcal{O}$-operator of
 	$(\frak g,\star_{\frak g},P)$ associated to
 	$(L_{\rhd_{\frak g}},R_{\lhd_{\frak g}},P,\frak g)$.
 	The remaining conclusions follow from
 	Theorem \ref{thm:differential-novikov-o-operator-construction}
 	and Proposition \ref{2303}.
 \end{proof}

\begin{pro}\label{1953}
Suppose that $(\frak g,\rhd_{\frak g},\lhd_{\frak g},P)$ is a differential
pre-Novikov algebra.
Let $\diamond_{\frak g}:\frak g\otimes\frak g\rightarrow\frak g$ be a multiplication given by
	\begin{equation*}
		x\diamond_{\frak g} y=P(x)\rhd_{\frak g} y-P(y)\lhd_{\frak g} x,\;\forall x,y\in\frak g.
	\end{equation*} 
	Then  $(\frak g,\diamond_{\frak g})$ is a pre-Lie algebra. 
	Moreover, define multiplications $\succ_{\frak g},\prec_{\frak g}:\frak g\otimes\frak g\rightarrow\frak g$ by
	\begin{equation*}
		x\succ_{\frak g}y=x\rhd_{\frak g}y+2y\lhd_{\frak g}x,\;
		x\prec_{\frak g}y=x\lhd_{\frak g}y+2y\rhd_{\frak g}x,\;
		\;\forall x,y\in\frak g.
	\end{equation*}
	Then $(\frak g, \succ_{\frak g},\prec_{\frak g})$ is a pre-APL algebra.
	Moreover,
	 \eqref{1906}-\eqref{1912} hold such that $(\frak g,\diamond_{\frak g},\succ_{\frak g},\prec_{\frak g})$ is a pre-compatible LAPL algebra.
\end{pro}
\begin{proof}
	It follows from a direct computation.
\end{proof}

Combining the preceding results, we obtain the following commutative
diagram.

\begin{equation*}
	\xymatrix@C=2.2cm@R=1.0cm{
		{\begin{gathered}
				\text{Manin triples of differential}\\
				\text{Novikov algebras}
		\end{gathered}}
		\ar[r]^{\text{Proposition }\ref{2013}}
		&
		{\begin{gathered}
				\text{Manin triples of compatible}\\
				\text{Lie algebras with respect to}\\
				\text{the invariant-c}2\text{c form}
		\end{gathered}}
		\\
		{\begin{gathered}
				\text{differential Novikov}\\
				\text{bialgebras}
		\end{gathered}}
		\ar[r]^{\text{Proposition }\ref{prop:differential-novikov-lapl}}
		\ar@{<->}[u]^{\text{Proposition }\ref{2101}}
		&
		{\begin{gathered}
				\text{compatible LAPL}\\
				\text{bialgebras}
		\end{gathered}}
		\ar@{<->}[u]_{\text{Corollary }\ref{1265}}
		\\
		{\begin{gathered}
				\text{antisymmetric solutions}\\
				\text{of the differential NYBE}
		\end{gathered}}
		\ar[r]^{\text{Proposition }\ref{prop:novikov-clapl-ybe}}
		\ar[u]^{\text{Proposition }\ref{2303}}
		&
		{\begin{gathered}
				\text{antisymmetric solutions}\\
				\text{of the CLAPL-YBE}
		\end{gathered}}
		\ar[u]_{\text{Corollary }\ref{cor:triangular-compatible-lapl}}
		\\
		{\begin{gathered}
				\mathcal{O}\text{-operators of differential}\\
				\text{Novikov algebras}
		\end{gathered}}
		\ar[r]^{\text{Proposition }\ref{prop:differential-novikov-lapl-o-operator}}
		\ar[u]^{\text{Theorem }\ref{thm:differential-novikov-o-operator-construction}}
		&
		{\begin{gathered}
				\mathcal{O}\text{-operators of}\\
				\text{compatible LAPL algebras}
		\end{gathered}}
		\ar[u]_{\text{Theorem }\ref{1859}}
		\\
		{\begin{gathered}
				\text{differential pre-Novikov}\\
				\text{algebras}
		\end{gathered}}
		\ar[r]^{\text{Proposition }\ref{1953}}
		\ar[u]^{\text{Proposition }\ref{prop:differential-pre-novikov-bialgebra}}
		&
		{\begin{gathered}
				\text{pre-compatible LAPL}\\
				\text{algebras}
		\end{gathered}}
		\ar[u]_{\text{Proposition }\ref{1916}}
	}
\end{equation*}
 
\noindent{\bf Acknowledgements.}
This work is supported by NSFC (W2412041),
the China Postdoctoral Science Foundation (2024T005TJ,
2024M761507), the Postdoctoral Fellowship Program of CPSF (GZC20240755), and Nankai Zhide Foundation.

\noindent{\bf Declaration of interests.}
The author has no conflicts of interest to disclose.

\noindent{\bf Data availability.}
Data sharing is not applicable as no data were created or analyzed.

\end{document}